\documentclass[11pt]{article}
\usepackage[top=25mm,bottom=25mm,left=25mm,right=25mm]{geometry}

\usepackage{amsmath,amssymb,amsthm,latexsym,amscd,stmaryrd,mathrsfs,longtable,color,url}
\usepackage[dvipdfmx]{hyperref}

\newtheoremstyle{standard}%
{9pt}%
{9pt}%
{\it}
{}%
{\bfseries}%
{}
{ }%
{#3}%

\newcommand{\db}[1]{(\!({#1})\!)}

\newcommand{\wak}{k}
\newcommand{\centralc}{C}
\newcommand{\subL}{P}

\numberwithin{equation}{section}

\newcommand{\N}{{\mathbb N}}
\newcommand{\Z}{{\mathbb Z}}
\newcommand{\Q}{{\mathbb Q}}

\newcommand{\C}{{\mathbb C}}

\newcommand{\wi}{i}

\newcommand{\wf}{f}
\newcommand{\Har}{H}

\newcommand{\ws}{s}
\newcommand{\wl}{l}

\newcommand{\rankL}{d}
\newcommand{\mN}{N}
\newcommand{\mW}{W}

\newcommand{\ewo}{u}
\newcommand{\ewt}{v}

\newcommand{\module}{M}

\newcommand{\hei}{{\mathfrak h}}

\newcommand{\sU}{{\mathscr U}}

\newcommand{\mmbox}[1]{\mbox{\scriptsize #1}}

\newcommand{\nor}{\begin{subarray}{c}\circ\\\circ\end{subarray}}

\newcommand{\fg}{{\mathfrak g}}
\newcommand{\fh}{{\mathfrak h}}

\newcommand{\ul}[1]{{#1}}
\newcommand{\lu}{u}
\newcommand{\lw}{w}
\newcommand{\lv}{v}

\newcommand{\mK}{K}

\newcommand{\lE}{t}
\newcommand{\ExB}{E}

\newcommand{\sv}{P}
\newcommand{\vac}{{\mathbf 1}}
\newcommand{\lattice}{L}

\newcommand{\nS}{S}

\newcommand{\eS}{\epsilon(S)}

\DeclareMathOperator{\Hom}{Hom}

\DeclareMathOperator{\Ker}{Ker}

\DeclareMathOperator{\Ext}{Ext} 
\DeclareMathOperator{\id}{id}

\DeclareMathOperator{\Span}{Span} 
\DeclareMathOperator{\wt}{wt}

\DeclareMathOperator{\rank}{rank} 

\DeclareMathOperator{\tw}{tw}

\newtheorem{lemma}{Lemma}[section]
\newtheorem{theorem}[lemma]{Theorem}
\newtheorem{proposition}[lemma]{Proposition}
\newtheorem{corollary}[lemma]{Corollary}
\theoremstyle{definition}
\newtheorem{definition}[lemma]{Definition}

\newtheorem{remark}[lemma]{Remark}

\theoremstyle{standard}

\title{The irreducible weak modules for the fixed point subalgebra of the vertex algebra associated to
a non-degenerate even lattice by an automorphism of order $2$ (Part $2$)}
\author{Kenichiro Tanabe\footnote{Research was partially supported by the Grant-in-aid
(No. 21K03172) for Scientific Research, JSPS.}\\\\
Faculty of Liberal Arts and Sciences\\
Tokyo City University\\
1-28-1 Tamazutsumi,Setagaya-ku,
Tokyo 158-8557\\
Japan\\
ktanabe@tcu.ac.jp}
\date{}
\begin{document}
\maketitle

\begin{abstract}
Let $V_{\lattice}$ be the vertex algebra  associated to a non-degenerate even lattice $\lattice$,
$\theta$ the automorphism of $V_{\lattice}$ induced from the $-1$ symmetry of $\lattice$, and
$V_{\lattice}^{+}$ the fixed point subalgebra of $V_{\lattice}$ under the action of $\theta$.
In this series of papers, we classify the irreducible weak $V_{\lattice}^{+}$-modules
and show that any irreducible weak $V_{\lattice}^{+}$-module 
is isomorphic to a weak submodule of some irreducible weak $V_{\lattice}$-module or 
to a submodule of some irreducible $\theta$-twisted $V_{\lattice}$-module.
Let $M(1)^{+}$ be the fixed point subalgebra of the Heisenberg vertex operator algebra $M(1)$ under the action of $\theta$.
In this paper (Part $2$), we show that there exists an irreducible  $M(1)^{+}$-submodule
in any non-zero weak $V_{\lattice}^{+}$-module and we compute extension groups for $M(1)^{+}$.
\end{abstract}



\tableofcontents
\section{\label{section:introduction}Introduction}
Let $L$ be a non-degenerate even lattice of finite rank $d$, $V_{\lattice}$ the vertex algebra associated to $\lattice$,
$\theta$ the automorphism of $V_{\lattice}$ induced from the $-1$ symmetry of $\lattice$, and
$V_{\lattice}^{+}$ the fixed point subalgebra of $V_{\lattice}$ under the action of $\theta$.
The  fixed point subalgebras play an important role in the study of vertex algebras.
For example, the moonshine vertex algebra 
$V^{\natural}$
  is constructed 
as a direct sum of $V_{\Lambda}^{+}$ and some irreducible $V_{\Lambda}^{+}$-module
in \cite{FLM} where $\Lambda$ is the Leech lattice.
The moonshine conjecture \cite{CN1979}, which is an unexpected connection between the monster group and modular functions,
was proved by Borcherds using  $V^{\natural}$ in \cite{B1992}.
The aim of this series of papers is to classify the irreducible weak $V_{\lattice}^{+}$-modules (see Definition \ref{definition:weak-module}
for the definition).
Because of the large number of pages in the original paper \cite{Tanabe2019}, we divide the paper into 3 parts  in a series for publication. 
This paper is Part $2$ and  a continuation of Part $1$ \cite{Tanabe2021-1}.
I will write the main result here again, which is stated in 
\cite[Theorem 1.1]{Tanabe2021-1}:
\begin{theorem}
\label{theorem:classification-weak-module}
Let $\lattice$ be a non-degenerate even lattice of finite rank with a bilinear form $\langle\ ,\ \rangle$.
The following is a complete set of representatives of equivalence classes of the irreducible weak $V_{\lattice}^{+}$-modules:
\begin{enumerate}
\item
$V_{\lambda+\lattice}^{\pm}$, $\lambda+\lattice\in \lattice^{\perp}/\lattice$ with $2\lambda\in \lattice$.   
\item
$V_{\lambda+\lattice}\cong V_{-\lambda+\lattice}$, $\lambda+\lattice\in \lattice^{\perp}/\lattice$ with $2\lambda\not\in \lattice$.   
\item
$V_{\lattice}^{T_{\chi},\pm}$ for any irreducible $\hat{\lattice}/P$-module $T_{\chi}$ with central character $\chi$.
\end{enumerate}
\end{theorem}
\noindent{}In the theorem, $L^{\perp}$ is the dual lattice of $\lattice$, 
$V_{\lambda+\lattice}^{\pm}=\{\lu\in V_{\lambda+\lattice}\ |\ \theta(\lu)=\pm \lu\}$
for $\lambda+\lattice\in \lattice^{\perp}/\lattice$ with $2\lambda\in \lattice$,
$\hat{\lattice}$ is the canonical central extension of $\lattice$
by the cyclic group $\langle\kappa\rangle$ of order $2$ with  the commutator map
$c(\alpha,\beta)=\kappa^{\langle\alpha,\beta\rangle}$ for $\alpha,\beta\in\lattice$,
$\subL=\{\theta(a) a^{-1}\ |\ a\in\hat{\lattice}\}$,
$V_{\lattice}^{T_{\chi}}$ is an irreducible $\theta$-twisted $V_{\lattice}$-module,
and $V_{\lattice}^{T_{\chi},\pm}=\{\lu\in V_{\lattice}^{T_{\chi}}\ |\ \theta(\lu)=\pm \lu\}$.
Note that in Theorem \ref{theorem:classification-weak-module}, 
$V_{\lattice}^{T_{\chi},\pm}$ in (3) are $V_{\lattice}^{+}$-modules,
however, if $\lattice$ is not positive definite,
then $V_{\lambda+\lattice}^{\pm}$ in (1) and $V_{\lambda+\lattice}$ in (2) are not  
$V_{\lattice}^{+}$-modules (cf. \cite[(2.18)]{Tanabe2021-1}).
See Section 1 of Part $1$ \cite{Tanabe2021-1} for the background and the detailed introduction to Theorem \ref{theorem:classification-weak-module}.
\textcolor{black}{
We note that since we do not assume any grading in the definition of a weak module, Theorem \ref{theorem:classification-weak-module} does not follow from \cite[Theorem 8.1]{Huang2020}, which deals with lower-bounded generalized modules for fixed point vertex algebras with some conditions.}
In fact, if $\lattice$ is not positive definite,
then the weak modules listed in (1) and (2) in Theorem \ref{theorem:classification-weak-module}
are not lower-bounded generalized modules.

Let  $M(1)$ be the Heisenberg vertex operator algebra associated to $\fh=\C\otimes_{\Z}\lattice$ (see the explanation under \eqref{eq:untwist-induced}
for the definition)
and $M(1)^{+}$ the fixed point subalgebra of $M(1)$ under the action of $\theta$.
The vertex operator algebra $M(1)^{+}$ is a subalgebra of $V_{\lattice}^{+}$
and, as stated in  \cite[Section 1]{Tanabe2021-1},
representations of $M(1)^{+}$ play a crucial role in the proof of Theorem \ref{theorem:classification-weak-module}.
The irreducible $M(1)^{+}$-modules 
are classified in \cite[Theorem 4.5]{DN1999-1} for the case of $\dim_{\C}\fh=1$ 
and \cite[Theorem 6.2.2]{DN2001} for the general case as follows: any irreducible $M(1)^{+}$-module is isomorphic to one of 
\begin{align}
\label{eqn:classificationM(1)p}
M(1)^{\pm}, M(1)(\theta)^{\pm}, \mbox{ or } M(1,\lambda)\cong M(1,-\lambda)\ (0\neq \lambda\in \fh).
\end{align}
Here $M(1)(\theta)$ is the irreducible $\theta$-twisted $M(1)$-module,
$M(1)^{\pm}=\{u\in M(1)\ |\ \theta \lu=\pm \lu\},
M(1)(\theta)^{\pm}=\{u\in M(1)(\theta)\ |\ \theta \lu=\pm \lu\}$,
and $M(1,\lambda)$ is the irreducible $M(1)$-module generated by
the vector $e^{\lambda}$ such that $(\alpha(-1)\vac)_{0}e^{\lambda}=\langle\alpha,\lambda\rangle e^{\lambda}$
and $(\alpha(-1)\vac)_{n}e^{\lambda}=0$ for  all $\alpha\in\fh$ and $n\in\Z_{>0}$
(See \eqref{eqn:hathei=heiotimes}--\eqref{eqn:M(1)(theta)pm=u} in Section \ref{section:preliminary} for the precise definitions of these symbols).
In the previous paper (Part 1), we showed that when the rank of $L$ is $1$,
for any non-zero weak $V_{\lattice}^{+}$-module $\module$ 
there exists a non-zero $M(1)^{+}$-submodule in $\module$.
In this paper  (Part $2$), we first strengthen and generalize this result to $\lattice$ of  an  arbitrary rank.
Precisely, we show that  for any non-zero weak $V_{\lattice}^{+}$-module $\module$,
there exists an irreducible $M(1)^{+}$-submodule in $\module$ 
(Proposition \ref{proposition:M(1)plusmoduleinweak} and Corollary \ref{corollary:M(1)plusirreduciblemoduleinweak}).
We next study extension groups and generalized Verma modules (see \cite[Theorem 6.2]{DLM1998t} for the definition) for $M(1)^{+}$ (Proposition \ref{proposition:Ext-total},
Corollary \ref{corollary:verma-irreducible}, and Lemma \ref{lemma:generalizedVermaModule-M(1)-}).
We shall explain how to use these results in Part $3$.
 Let $\module$ be an irreducible weak $V_{\lattice}^{+}$-module.
By Corollary \ref{corollary:M(1)plusirreduciblemoduleinweak}, there exists an irreducible $M(1)^{+}$-submodule $\mK$ of $\module$.
If $\mK\cong M(\theta)^{\pm}$, then the same argument as in Section \ref{section:Modules for the Zhu algera of general} of the present paper
shows that 
$\module$ is a $V_{\lattice}^{+}$-module.
In this case, \cite[Proposition 4.15]{Yamsk2008} shows that
$M$ is isomorphic to one of the irreducible weak $V_{\lattice}^{+}$-modules in Theorem \ref{theorem:classification-weak-module} (3).
Assume $\mK\cong M(1)^{\pm}$ or $M(1,\lambda)$ with $0\neq \lambda\in \fh$.
Since $V_{\lattice}^{+}$ is a direct sum of irreducible
$M(1)^{+}$-modules, for any irreducible $M(1)^{+}$-submodule $\mN$ of $V_{\lattice}^{+}$,
the $V_{\lattice}^{+}$-module structure of $\module$
induces an intertwining operator $I(\mbox{ },x) : \mN\times K\rightarrow\module\db{x}$ for weak $M(1)^{+}$-modules (see Definition \ref{definition:intertwining-operator} for the definition).
We denote by $Q$ the weak $M(1)^{+}$-submodule of $\module$ that is the image of $I(\mbox{ },x)$.
The same argument as in Section \ref{section:Modules for the Zhu algera of general} of the present paper
shows that there exists an irreducible $M(1)^{+}$-submodule $R$ of $Q$.
Moreover, If $R\neq Q$, then $Q/R$ is an irreducible $M(1)^{+}$-module
and by Proposition \ref{proposition:Ext-total} the exact sequence 
$0\rightarrow R\rightarrow Q\rightarrow Q/R\rightarrow 0$ splits, namely
$Q\cong R\oplus Q/R$ as  $M(1)^{+}$-modules.
Since $V_{\lattice}^{+}$ is a direct sum of irreducible
$M(1)^{+}$-modules,
this leads to the result that $\module$ is a direct sum of irreducible $M(1)^{+}$-modules.
Moreover, we find that the irreducible $M(1)^{+}$-modules in the direct sum
are pairwise non-isomorphic.
Using \textcolor{black}{fusion rules} (see the explanation under Definition \ref{definition:intertwining-operator}) for the irreducible $M(1)^{+}$-modules obtained in \cite[Theorem 5.13]{Abe2001} and \cite[Theorem 7.7]{ADL2005},
we can determine the weak $V_{\lattice}^{+}$-module $\module$ with such an $M(1)^{+}$-module structure and thus $\module$ is one of the irreducible weak  $V_{\lattice}^{+}$-modules in Theorem \ref{theorem:classification-weak-module} (1) and (2).

Let us briefly explain the basic idea to show Proposition \ref{proposition:M(1)plusmoduleinweak}, the main result in Section \ref{section:Modules for the Zhu algera of general}.
Let $V$ be a vertex algebra and  $\module$ a weak $V$-module.
For $a\in V$ and $\lu\in \module$, we define  $\epsilon(a,\lu)\in\Z\cup\{-\infty\}$ by
\begin{align}
\label{eqn:max-vanish0}
a_{\epsilon(a,\lu)}\lu&\neq 0\mbox{ and }a_{i}\lu
=0\mbox{ for all }i>\epsilon(a,\lu)
\end{align}
if $Y_{\module}(a,x)\lu\neq 0$ and $\epsilon(a,\lu)=-\infty$
if $Y_{\module}(a,x)\lu= 0$.
It is well-known that the vertex operator algebra $M(1)^{+}$ is generated by
homogeneous elements
$\omega^{[i]}$ of weight $2$, $J^{[i]}$( or $H^{[i]}$) of weight $4$, and $S_{lm}(1,r)$ 
 of weight $r+1$ $(1\leq i\leq \rankL, 1\leq m<l\leq \rankL, r=1,2,3)$
such that $[\omega^{[i]}_{k}, \omega^{[j]}_{l}]=[\omega^{[i]}_{k}, H^{[j]}_{l}]=[H^{[i]}_{k}, H^{[j]}_{l}]=0$
for any $k,l\in\Z$ and any pair of distinct elements $i,j\in\{1,\ldots,\rankL\}$
(see \eqref{eq:def-oega-i-H-i} and \eqref{eq:def-nS_ij(l,m)} 
for these symbols). 
Hence for any non-zero  weak $V_{\lattice}^{+}$-module $\module$,
it follows from \cite[Lemma 3.7]{Tanabe2021-1} that there exists a simultaneous eigenvector $\lu$ of
 $\{\omega^{[i]}_{1},\Har^{[i]}_{3}\}_{i=1}^{\rankL}$ in $\module$
such that $\epsilon(\omega^{[i]},\lu)\leq 1$ and
$\epsilon(\Har^{[i]},\lu)\leq 3$ for all $i=1,\ldots,\rankL$.
By induction on $\max\{\epsilon(S_{ij}(1,1),\lu)\ |\ i>j\}$, we get a simultaneous eigenvector of
 $\{\omega^{[i]}_{1},\Har^{[i]}_{3}\}_{i=1}^{\rankL}$, which we denote by the same symbol $u$, 
such that  $\epsilon(\omega^{[i]},\lu)\leq 1$,
$\epsilon(\Har^{[i]},\lu)\leq 3$,  and $\epsilon(S_{lm}(1,r),\lu)\leq r$
for all $i=1,\ldots,\rankL, 1\leq m<l\leq \rankL$, and $r=1,2,3$.
Namely, $\lu\in\Omega_{M(1)^{+}}(\module)$ (see \eqref{eq:OmegaV(U)=BigluinU} for the definition). Since $\lu$ is a simultaneous eigenvector of
 $\{\omega^{[i]}_{1},\Har^{[i]}_{3}\}_{i=1}^{\rankL}$, 
$A(M(1)^{+})\lu$ is of \textcolor{black}{finite dimension}, where $A(M(1)^{+})$ is the Zhu algebra 
for $M(1)^{+}$  (see \eqref{eq:zhu-ideal-multi}--\eqref{eq:zhu-bimodule} for the definition).
Hence, by \cite[Theorem 6.2]{DLM1998t} we have the result.

We next explain the basic idea to show Proposition \ref{proposition:Ext-total}, the main result in Section \ref{section:Extension groups for M(1)+}.
The result shows that in most cases the exact sequence 
\begin{align}
\label{eq:intro-exactsequence}
0\rightarrow W\overset{}{\rightarrow} N\overset{\pi}{\rightarrow} M\rightarrow 0
\end{align}
splits for two irreducible $M(1)^{+}$-modules $\mW=\oplus_{i\in\delta+\Z_{\geq 0}}\mW_{i}, \module=\oplus_{i\in\gamma+\Z_{\geq 0}}M_{i}$ and a weak $M(1)^{+}$-module $\mN$,
where $\omega=\sum_{j=1}^{\rankL}\omega^{[j]}$ is the conformal vector (Virasoro element) of $M(1)^{+}$
and $M_{i}:=\{u\in M\ |\ \omega_1u=i u\}$ for $i\in\C$.
Precisely, in Section \ref{section:Extension groups for M(1)+}
we deal with the case where $\module$ is a general $M(1)^{+}$-module in the exact sequence \eqref{eq:intro-exactsequence}
in order to show Corollary \ref{corollary:verma-irreducible}, however, here we assume $\module$ is irreducible to simplify the argument. 
As in Part $1$ \cite{Tanabe2021-1},
we first find some relations for $\omega^{[i]},\Har^{[i]}\ (i=1,\ldots,\rankL), S_{lm}(1,r)\ (1\leq m<l\leq \rankL,r=1,2,3)$ in $M(1)^{+}$ 
with the help of computer algebra system 
Risa/Asir\cite{Risa/Asir} (\eqref{eq:s11-3}--\eqref{eq:s11-4-3}).
For $\zeta=(\zeta^{[1]},\ldots,\zeta^{[\rankL]}), \xi=(\xi^{[1]},\ldots,\xi^{[\rankL]})\in \C^{\rankL}$, 
let $\lv\in \module_{\gamma}$ such that 
$(\omega^{[i]}_{1}-\zeta^{[i]})\lv=(\Har^{[i]}_{3}-\xi^{[i]})\lv=0$
for all $i=1,\ldots,\rankL$.
Assume $(W,M_{\gamma})\not\cong (M(1)^{+},M(1)^{-}_{1}), (M(1)^{-},M(1)^{+}_{0})$.
Using the relations \eqref{eq:(epsilonS-1)big((18zeta[i]+3)epsilonS5-0}--\eqref{eq:s11-zeta-2} obtained by \eqref{eq:s11-3}--\eqref{eq:s11-4-3}, we can take $\lu\in N_{\gamma}$ such that $\pi(\lu)=\lv$,
$(\omega^{[i]}_{1}-\zeta^{[i]})\lu=(\Har^{[i]}_{3}-\xi^{[i]})\lu=0$ for all $i=1,\ldots,\rankL$
and $\lu\in \Omega_{M(1)^{+}}(N_{\gamma})$ (Lemmas \ref{lemma:M1lambda-submodule} and \ref{lemma:M1lambda-submodule-2}).
After studying the following two cases, we have the result:
\begin{enumerate}
\item The case that $\module\not\cong \mW$.
Assume that the exact sequence \eqref{eq:intro-exactsequence} does not split.
Since the intersection of $\mW$ and the $M(1)^{+}$-submodule of $\mN$ generated by $\lu$
is not trivial, we have $\delta\in\gamma+\Z_{\geq 0}$. 
By taking the restricted dual of \eqref{eq:intro-exactsequence},
the same argument shows that $\gamma\in\delta+\Z_{\geq 0}$ and hence $\delta=\gamma$.
Since $\module\not\cong\mW$, we have $\mN_{\gamma}\cong \mW_{\gamma}\oplus \module_{\gamma}$ as $A(M(1)^{+})$-modules
and hence the exact sequence \eqref{eq:intro-exactsequence} splits. This is a contradiction (Lemma \ref{lemma:Ext-split-lambda-mu}).
\item
The case that $\module=\mW$ and $M\in \{M(1)^{\pm}, M(1)(\theta)^{\pm}\}$.
Using the relations \eqref{eq:(epsilonS-1)big((18zeta[i]+3)epsilonS5-0}--\eqref{eq:s11-zeta-2} again, we have $N_{\gamma}\cong M_{\gamma}\oplus W_{\gamma}$ as  $A(M(1)^{+})$-modules
and hence the exact sequence \eqref{eq:intro-exactsequence} splits (Lemma \ref{lemma:Ext-M-M}).
\end{enumerate}

Complicated computation has been done by a computer algebra system Risa/Asir\cite{Risa/Asir}.
Throughout this paper, the word \lq\lq a direct computation\rq\rq \ often means a direct computation with the help of Risa/Asir.
Details of computer calculations such as \eqref
{eq:Sij(2,1)=omega_0S_ij(1,1)-S_ij(1,2)}, \eqref{eq:s11-3}, \eqref{eq:(epsilonS-1)big((18zeta[i]+3)epsilonS5-0}, 
\eqref{eqn:epsilon(Si1)5epsilon(Si1-1)4}, etc.,
and \eqref{eq:omega[j]0Sij(1,1)=Sij(1,2)-1vac}--\eqref{eq:last-appendix} in \ref{section:appendix}  
can be found on the internet at \cite{computer-nondeg-va-part2}.

The organization of the paper is as follows.
\textcolor{black}{In Section \ref{section:preliminary}} we recall some basic properties of 
weak modules for a vertex algebra. We also recall the Heisenberg algebra $M(1)$
and its fixed point algebra $M(1)^{+}$.
In Section \ref{section:Modules for the Zhu algera of general} we show that for any non-zero weak $V_{\lattice}^{+}$-module $\module$
there exists a non-zero submodule for $M(1)^{+}$ in $\module$.
In Section \ref{section:Extension groups for M(1)+} we study extension groups and generalized Verma modules for $M(1)^{+}$.
In \ref{section:appendix} we put 
computations of $a_{k}b$ for some $a,b\in V_{\lattice}^{+}$
and $k=0,1,\ldots$ to find the commutation relation $[a_{i},b_{j}]=\sum_{k=0}^{\infty}\binom{i}{k}(a_{k}b)_{i+j-k}$. 
In \ref{section:notation} we list some notation.
\section{\label{section:preliminary}Preliminary}
We assume that the reader is familiar with the basic knowledge on
vertex algebras as presented in \cite{B1986,FLM,LL,Li1996}. 

Throughout this paper, $V$ is a vertex algebra and 
we always assume that $V$ has an element $\omega$ such that $\omega_{0}a=a_{-2}\vac$ for all $a\in V$.
For a vertex operator algebra $V$, this condition automatically holds since $V$ has the conformal vector (Virasoro element).
Throughout this paper, we follow the notation and terminology of \cite{Tanabe2021-1}.
We will explain some of them. 
We note that  if $V$ is a vertex operator algebra, then 
the notion of a module for $V$ viewed as a vertex algebra is different from the notion of a module for $V$ viewed as a vertex operator algebra (cf. \cite[Definitions 4.1.1 and 4.1.6]{LL}).
To avoid confusion, throughout this paper, we refer to a module for a vertex algebra defined in \cite[Definition 4.1.1]{LL} as a {\it weak module}.
Here we write down the definition of a weak $V$-module:
\begin{definition}
\label{definition:weak-module}
A {\it weak $V$-module} $\module$ is a vector space over $\C$ equipped with a linear map
\begin{align}
\label{eq:inter-form}
Y_{\module}(\ , x) : V\otimes_{\C}\module&\rightarrow \module\db{x}\nonumber\\
a\otimes u&\mapsto  Y_{\module}(a, x)\lu=\sum_{n\in\Z}a_{n}\lu x^{-n-1}
\end{align}
such that the following conditions are satisfied:
\begin{enumerate}
\item $Y_{\module}(\vac,x)=\id_{\module}$.
\item
For $a,b\in V$ and $\lu\in \module$,
\begin{align}
\label{eq:inter-borcherds}
&x_0^{-1}\delta(\dfrac{x_1-x_2}{x_0})Y_{\module}(a,x_1)Y_{\module}(b,x_2)\lu-
x_0^{-1}\delta(\dfrac{x_2-x_1}{-x_0})Y_{\module}(b,x_2)Y_{\module}(a,x_1)\lu\nonumber\\
&=x_1^{-1}\delta(\dfrac{x_2+x_0}{x_1})Y_{\module}(Y(a,x_0)b,x_2)\lu.
\end{align}
\end{enumerate}
\end{definition}
For $n\in\C$ and a weak $V$-module $\module$, we define $M_{n}=\{\lu\in V\ |\ \omega_1 \lu=n \lu\}$.
For $a\in V_{n}\ (n\in\C)$, $\wt a$ denotes $n$.
For a vertex algebra $V$ which admits a decomposition $V=\oplus_{n\in\Z}V_n$ and a subset $U$ of a weak $V$-module, we 
define
\begin{align}
\label{eq:OmegaV(U)=BigluinU}
\Omega_{V}(U)&=\Big\{\lu\in U\ \Big|\ 
\begin{array}{l}
a_{i}\lu=0\ \mbox{for all homogeneous elements }a\in V\\
\mbox{and }i>\wt a-1.
\end{array}\Big\}.
\end{align}
For a vertex algebra $V$ which admits a decomposition $V=\oplus_{n\in\Z}V_n$, a weak $V$-module $\mN$
 is called {\it $\N$-graded} if $N$ admits a decomposition $N=\oplus_{n=0}^{\infty}N(n)$
such that $a_{i}N(n)\subset N(\wt a-i-1+n)$ for all homogeneous elements $a\in V$, $i\in\Z$, and $n\in\Z_{\geq 0}$, where 
we define $N(n)=0$ for all $n<0$.
For a vertex algebra $V$ which admits a decomposition $V=\oplus_{n\in\Z}V_n$, a weak $V$-module $\mN$
 is called a {\it $V$-module} if $\mN$ admits a decomposition $\mN=\oplus_{n\in\C}\mN_{n}$
such that $\dim_{\C}\mN_{n}<\infty$ for all $n\in\C$ and $\mN_{n}=0$ for $n$ whose real part is sufficiently negative.
We recall the definition of an intertwining operator
from \cite[Definition 5.4.1]{FHL}.
\begin{definition}
	\label{definition:intertwining-operator}
Let $V$ be a vertex algebra and let $\module, \mW$, and $\mN$ be three weak $V$-modules. 
	An {\em intertwining operator} of type $\binom{\mN}{\module\ \mW}$ is a linear map 
	\begin{align}
		\label{eq:inter-form}
		I(\ , x) : \module\otimes_{\C}\mW&\rightarrow \mN\{x\}\nonumber\\
		I(\ewo, x)\ewt&=\sum_{\alpha\in\C}\ewo_{\alpha}\ewt x^{-\alpha-1},\nonumber\\
		&\ewo\in \module,\ewt\in \mW, \mbox{ and }\ewo_{\alpha}\in \Hom_{\C}(\mW,\mN),
	\end{align}
	such that the following conditions are satisfied:
	\begin{enumerate}
		\item
		For $\ewo\in \module,\ewt\in \mW$, and $\alpha\in\C$,
		\begin{align}
			\label{eq:inter-truncation}
			\ewo_{\alpha+m}\ewt&=0\mbox{ for sufficiently large $m\in\N$.}
		\end{align}
		
		\item
		For $\ewo\in \module$ and $a\in V$,
		\begin{align}
			\label{eq:inter-borcherds}
			&x_0^{-1}\delta(\dfrac{x_1-x_2}{x_0})Y(a,x_1)I(\ewo,x_2)-
			x_0^{-1}\delta(\dfrac{x_2-x_1}{-x_0})I(\ewo,x_2)Y(a,x_1)\nonumber\\
			&=x_1^{-1}\delta(\dfrac{x_2+x_0}{x_1})I(Y(a,x_0)\ewo,x_2).
		\end{align}
		\item
		For $\ewo\in \module$,
		\begin{align}
			\label{eq:inter-derivative}
		I(\omega_{0}\ewo,x)&=\dfrac{d}{dx}I(\ewo,x).
	\end{align}\end{enumerate}
\end{definition}
For irreducible weak $V$-modules $\module,\mW$, and $\mN$, 
 $I\binom{\mN}{\module\ \mW}$ denotes
the space of all intertwining operators of type $\binom{\mN}{\module\ \mW}$
and we call its dimension the {\it fusion rule} of type $\binom{\mN}{\module\ \mW}$. 
In this paper, for an intertwining operator 
$I(\mbox{ },x)$ from $\module\times \mW$ to $\mN$,
we consider only the case that the image of $I(\mbox{ },x)$ is contained in $\mN\db{x}$.
For $A\subset \module$ and $B\subset\mW$,
\begin{align}
A\cdot B\mbox{ denotes }\Span_{\C}\{a_{i}b\ |\ a\in A, i\in\Z, b\in B\}\subset N. 
\end{align}
For an intertwining operator $I(\mbox{ },x) : \module\times \mW\rightarrow \mN\db{x}$,
$\lu\in\module$, and $\lv\in \mW$, we define  $\epsilon_{I}(\lu,\lv)=\epsilon(\lu,\lv)\in\Z\cup\{-\infty\}$ by
\begin{align}
\label{eqn:max-vanish}
\lu_{\epsilon_{I}(\lu,\lv)}\lv&\neq 0\mbox{ and }\lu_{i}\lv
=0\mbox{ for all }i>\epsilon_{I}(\lu,\lv)
\end{align}
if $I(\lu,x)\lv\neq 0$ and $\epsilon_{I}(\lu,\lv)=-\infty$
if $I(\lu,x)\lv= 0$.
For a subset $A$ of $V$ and a subset $B$ of a weak $V$-module $\module$,
let 
\begin{align}
\label{eqn:A-B:=Span}
A_{-}B&:=\Span_{\C}\{a_{-i}b\ |\ a\in A, b\in B,\mbox{ and }i\in\Z_{>0}\}\subset \module
\end{align}
and 
\begin{align}
\label{eq:a(1)-i1cdotsa(n)-inb}
\langle A_{-}\rangle B&:=\Span_{\C}\Big\{a^{(1)}_{-i_1}\cdots a^{(n)}_{-i_n}b\ \Big|\ 
\begin{array}{l}
n\in\Z_{\geq 0}, a^{(1)},\ldots,a^{(n)}\in A, b\in B,\\
i_1,\ldots,i_n\in\Z_{>0}
\end{array}\Big\}\subset \module.
\end{align}
When $B=\{b\}$,  we will simply write  
$A_{-}B$ 
and 
$\langle A_{-}\rangle B$ as $A_{-}b$
and $\langle A_{-}\rangle b$, respectively.

\begin{lemma}
\label{lemma:ksumklimitsm=1n(kepsilon(a(m))+im)kk+(kepsilon(b)+j)-1+kr}
Let $V=\oplus_{n=0}^{\infty}V_n$ be a vertex operator algebra,
$A$ and $B$ finite subsets of $\bigcup_{n=1}^{\infty}V_n$, 
$\module$ a weak $V$-module, and \textcolor{black}{$U$ a finite dimensional subspace} of $\module$.
Assume $A\cdot \langle A_{-}\rangle\vac\subset \langle A_{-}\rangle\vac$.
For each $a\in A$ and $b\in B$, we choose $\epsilon(a), \epsilon(b)\in \Z$
so that $\epsilon(a)\geq \max\{\{\epsilon(a,u)\ |\ u\in U\}\cup\{-1\}\}$ and $\epsilon(b)\geq \max\{\{\epsilon(b,u)\ |\ u\in U\}\cup\{-1\}\}$.
For 
$j\in\Z_{\geq 0}$, we define
\begin{align}
\label{eq:gamma(k):=maxBiggsum}
\delta(j)&:=
\max\Bigg\{
\begin{array}{l}
\sum\limits_{i=1}^{n}(\epsilon(a^{(i)})-\wt(a^{(i)})+1)\\
\end{array}\ \Bigg|
\
\begin{array}{l}
n\in\Z_{\geq 0},
a^{(1)},\ldots,a^{(n-1)}\in A, a^{(n)}\in B,\\
\sum\limits_{i=1}^{n}\wt(a^{(i)})\leq j.
\end{array}
\Bigg\}-1+j,
\end{align}
where $\sum_{i=1}^{n}(\epsilon(a^{(i)})-\wt(a^{(i)})+1):=0$ for $n=0$.
\begin{enumerate}
\item
Assume $A\cdot (\langle A_{-}\rangle B_{-}\vac)\subset \langle A_{-}\rangle B_{-}\vac$.
For each $a\in A$ and $b\in B$, we choose $\gamma(a),\gamma(b)\in\Z_{\geq -1}$.
Then, for a homogeneous element $c\in \langle A_{-}\rangle B_{-}\vac$, $u\in U$, and $k\in \Z$, $c_ku$ 
is a linear combination of elements of the form
\begin{align}
\label{eq:p(1)epsilon(p(1))cdots p(m)epsilon(p(m)}
p^{(1)}_{\sigma_1}\cdots p^{(l)}_{\sigma_l}
q_{\tau}p^{(l+1)}_{\sigma_{l+1}}\cdots p^{(m)}_{\sigma_m}
\lu
\end{align}
where $l,m\in\Z_{\geq 0}$ with $0\leq l\leq m$, $p^{(1)},\ldots,p^{(m)}\in A$, 
$\sigma_i\in \Z_{\leq \gamma(p^{(i)})}\ (i=1,\ldots,l)$,
$\sigma_i\in \Z_{\geq \gamma(p^{(i)})+1}\ (i=l+1,\ldots,m)$,
$q\in B, \tau\in\Z$ such that 
\begin{align}
\label{eq:sumi=1mwt(p(i))+wtq}
\sum_{i=1}^{m}\wt (p^{(i)})+\wt (q)&\leq \wt(c)\mbox{ and}\\
\label{eq:wt(c)-k-1=sum_i=1m}\wt (c)-k-1&=\sum_{i=1}^{m}(\wt (p^{(i)})-\sigma_i-1)+\wt (q)-\tau-1.
\end{align}
In particular, $c_{\delta(\wt(c))}u$ 
is a linear combination of elements of the form
\begin{align}
\label{eq:p(1)epsilon(p(1))cdots p(m)epsilon(p(m)-1}
p^{(1)}_{\epsilon(p^{(1)})}\cdots p^{(m)}_{\epsilon(p^{(m)})}q_{\epsilon(q)}\lu
\end{align}
where $m\in\Z_{\geq 0}$, $p^{(1)},\ldots,p^{(m)}\in A$, and $q\in B$.
Moreover, for $k>\delta(\wt(c))$,
\begin{align}
\label{eq:cru=0}
c_{k}u&=0.
\end{align}
\item 
Assume  $B\cdot \langle A_{-}\rangle\vac \subset B_{-}\langle A_{-}\rangle\vac$.
Then, for a homogeneous element $c\in B_{-}\langle A_{-}\rangle\vac$ and $u\in U$, $c_{\delta(\wt(c))}u$ is a linear combination of elements of the form
\begin{align}
\label{eq:p(1)epsilon(p(1))cdots p(m)epsilon(p(m)-2}
q_{\epsilon(q)}p^{(1)}_{\epsilon(p^{(1)})}\cdots p^{(m)}_{\epsilon(p^{(m)})}\lu
\end{align}
where $m\in\Z_{\geq 0}$, $p^{(1)},\ldots,p^{(m)}\in A$, and $q\in B$. Moreover,
for $k>\delta(\wt(c))$,
\begin{align}
\label{eq:cru=0-1}
c_{k}u&=0.
\end{align}
\item 
Assume $B\cdot \langle A_{-}\rangle\vac \subset B_{-}\langle A_{-}\rangle\vac$
and $A\cdot (B_{-}\langle A_{-}\rangle\vac) \subset B_{-}\langle A_{-}\rangle\vac$.
For $a\in A$, we define
\begin{align}
\label{eq:zeta(a):=maxepsilon(a)}
\zeta(a):=\max\{\{\epsilon(a)\}\cup \{\delta(\wt a+\wt b-1)-\epsilon(b)\ |\ b\in B\}\}.
\end{align}
If $a_{\epsilon(a)}\lu\in U$ for all $a\in A$ and $u\in U$, 
then the subspace $W:=\Span_{\C}\{b_{\epsilon(b)}\lu\ |\ b\in B, u\in U\}$ of $\module$
is stable under the action of $a_{\zeta(a)}$ for all $a\in A$. Moreover, for $a\in A$ and $k>\zeta(a)$,
$a_{k}\mW=0$.
\item 
Assume $B\cdot \langle A_{-}\rangle\vac \subset B_{-}\langle A_{-}\rangle\vac$, $A\cdot (B_{-}\langle A_{-}\rangle\vac) \subset B_{-}\langle A_{-}\rangle\vac$, and
$B\cdot (B_{-}\langle A_{-}\rangle\vac) \subset B_{-}\langle A_{-}\rangle\vac$.
For any $a\in A$, $b\in B$, and $u\in U$, we assume $\epsilon(a,\lu)\geq \wt(a)-1$
and the value $\epsilon(b)-\wt(b)+1$ is a constant independent of $b\in B$, which we denote by $\rho$.
We define $W:=\Span_{\C}\{b_{\epsilon(b)}\lu\ |\ b\in B, u\in U\}$.
If $a_{\wt(a)-1}\lu\in U$ for all $a\in A$ and $\lu\in U$, then for any homogeneous element $c\in B_{-}\langle A_{-}\rangle\vac$,
$c_{\wt(c)-1}W\subset W$ and $c_{k}W=0$ for all $k>\wt(c)-1$. 
\end{enumerate}

\end{lemma}

\begin{proof}
\begin{enumerate}
\item
Let 
$c:=a^{(1)}_{-i_1}\cdots a^{(n)}_{-i_n}b_{-i_{n+1}}\vac\in \langle A_{-}\rangle B_{-}\vac$
where $n\in\Z_{\geq 0}$, $a^{(1)},\ldots,a^{(n)}\in A$, $b\in B$, and $i_1,\ldots,i_n, i_{n+1}\in\Z_{>0}$.
We shall show \eqref{eq:p(1)epsilon(p(1))cdots p(m)epsilon(p(m)}--\eqref{eq:wt(c)-k-1=sum_i=1m} by induction on $\wt c$.
If $c=\textbf{0}$ or $c=b_{-i}\vac$ with $b\in B$ and $i\in\Z_{>0}$, then the results hold.
Let $n\geq 1$.
We define $\wf:=a^{(2)}_{-i_2}\cdots a^{(n)}_{-i_n}b_{-i_{n+1}}\vac$
and note that $\wt(\wf)<\wt(c)$. 
For $k\in\Z$, using \cite[Lemma 2.2]{Tanabe2021-1}, we have
\begin{align}
\label{eq:(a(1)_-i_1lv)_klu-0}
c_{k}\lu&=(a^{(1)}_{-i_1}\wf)_{k}\lu\nonumber\\
&=\sum_{\begin{subarray}{l}\ws\leq \gamma(a^{(1)})\\\ws+t+i_1=k\end{subarray}}
\binom{-\ws-1}{i_1-1}a^{(1)}_{\ws}\wf_{t}\lu+\sum_{\begin{subarray}{l}\ws\geq \gamma(a^{(1)})+1\\\ws+t+i_1=k\end{subarray}}\binom{-\ws-1}{i_1-1}\wf_{t}a^{(1)}_{\ws}\lu\nonumber\\
&\quad{}+(-1)^{i_1}\sum_{\wl=0}^{\infty}\binom{\wl+i_1-1}{i_1-1}
\binom{\gamma(a^{(1)})+i_1}{\wl+i_1}(a^{(1)}_{\wl}\wf)_{k-i_1-\wl}\lu.
 \end{align}
In the second term in \eqref{eq:(a(1)_-i_1lv)_klu-0},
by the induction hypothesis in the setting $c$, $k$, and $\lu$ are replaced by 
$f$, $t$, and $a^{(1)}_{\ws}\lu$, respectively, we find that $\wf_{t}a^{(1)}_{\ws}\lu$
is a linear combination of elements of the form
\eqref{eq:p(1)epsilon(p(1))cdots p(m)epsilon(p(m)}--\eqref{eq:wt(c)-k-1=sum_i=1m}.
In the third term in \eqref{eq:(a(1)_-i_1lv)_klu-0}, since $\wt (a^{(1)}_{\wl}\wf)=\wt(c)-i_1-l<\wt(c)$ for $l\in\Z_{\geq 0}$, by the induction hypothesis
in the setting $c$ and $k$ are replaced by $a^{(1)}_{\wl}\wf$ and $k-i_1-\wl$, respectively, we find that
$(a^{(1)}_{\wl}\wf)_{k-i_1-\wl}\lu$ 
is a linear combination of elements of the form
\eqref{eq:p(1)epsilon(p(1))cdots p(m)epsilon(p(m)}--\eqref{eq:wt(c)-k-1=sum_i=1m}.
Hence \eqref{eq:p(1)epsilon(p(1))cdots p(m)epsilon(p(m)}--\eqref{eq:wt(c)-k-1=sum_i=1m} hold.

We shall show \eqref{eq:p(1)epsilon(p(1))cdots p(m)epsilon(p(m)-1} and \eqref{eq:cru=0}. 
Let $k\geq \delta(\wt c))$.
We set 
$\gamma(a)=\epsilon(a)$ for $a\in A$ and 
$\gamma(b)=\epsilon(b)$ for $b\in B$.
Assume that in the expansion of $c_k\lu$, the coefficient of an element 
of the form 
\eqref{eq:p(1)epsilon(p(1))cdots p(m)epsilon(p(m)} is not zero.
Since $p_j\lu=0$ for $p\in A$ and $j>\epsilon(p)$, we have $l=m$ and hence $\tau\leq \epsilon(q)$. 
By  \eqref{eq:wt(c)-k-1=sum_i=1m},
\begin{align*}
0&
=-\wt (c)+k+1+\sum_{i=1}^{m}(\wt (p^{(i)})-\sigma_i-1)+\wt (q)-\tau-1\\
&\geq -\wt (c)+\delta(\wt(c))+1+\sum_{i=1}^{m}(\wt (p^{(i)})-\sigma_i-1)+\wt (q)-\tau-1\\
&\geq -\wt (c)+\big(\sum\limits_{i=1}^{m}(\epsilon(p^{(i)})-\wt(p^{(i)})+1)+(\epsilon(q)-\wt(q)+1)-1+\wt (c)\big)+1\\
&\quad{}+\sum_{i=1}^{m}(\wt (p^{(i)})-\sigma_i-1)+\wt (q)-\tau-1\\
&=\sum\limits_{i=1}^{m}(\epsilon(p^{(i)})-\sigma_i)+(\epsilon(q)-\tau)\geq 0
\end{align*}
and hence $k=\delta(\wt (c))$, 
$\epsilon(p^{(i)})=\sigma_i$ for all $i=1,\ldots,m$ and $\epsilon(b)=\tau$.
Here we have used \eqref{eq:sumi=1mwt(p(i))+wtq} and  the definition \eqref{eq:gamma(k):=maxBiggsum} of $\delta(\wt(c))$.
Thus, \eqref{eq:p(1)epsilon(p(1))cdots p(m)epsilon(p(m)-1} and \eqref{eq:cru=0} hold.

\item
The same argument as in (1) shows the results.

\item
For $j\in\Z_{>0}$, by the definition of $\delta(j)$
we have
$\delta(j)-1\geq \delta(j-1)$
and 
hence 
$\delta(j)-i\geq \delta(j-i)$ for all $i=0,\ldots,j$.
Let $\lu\in U$ and $k\in\Z_{\geq \zeta(a)}$.
For $a\in A$ and $b\in B$,
\begin{align}
a_{k}b_{\epsilon(b)}\lu&=b_{\epsilon(b)}a_{k}\lu+[a_{k},b_{\epsilon(b)}]\lu
=b_{\epsilon(b)}a_{k}\lu+\sum_{i=0}^{\infty}\binom{k}{i}(a_{i}b)_{k+\epsilon(b)-i}\lu.
\end{align}
For $i\in \Z_{\geq 0}$, since 
\begin{align}
k+\epsilon(b)-i&\geq 
\zeta(a)+\epsilon(b)-i\geq \delta(\wt a+\wt b-1)-i\nonumber\\
&\geq \delta(\wt a+\wt b-1-i)=\delta(\wt(a_ib)),
\end{align}
we have the results by (2).
\item 
For any $a\in A$ and $u\in U$, since $\epsilon(a,\lu)\geq \wt(a)-1$, we choose $\epsilon(a)=\wt(a)-1$.
For a homogeneous element $c\in B_{-}\langle A_{-}\rangle\vac$, by the definition \eqref{eq:gamma(k):=maxBiggsum} of $\delta$,
\begin{align}
\label{eq:beginarrayl-1+wt(c)endarrayBigg|}
\delta(\wt(c))&
=\max\{\epsilon(b)-\wt(b)\ |\ b\in B\}+\wt(c)=\rho+\wt(c)-1.
\end{align}
For $a\in A$ and $b\in B$, by \eqref{eq:beginarrayl-1+wt(c)endarrayBigg|},
\begin{align}
\delta(\wt(a)+\wt(b)-1)-\epsilon(b)&
=\rho+\wt(a)+\wt(b)-2-\epsilon(b)\nonumber\\
&=\epsilon(b)-\wt(b)+\wt(a)+\wt(b)-1-\epsilon(b)\nonumber\\
&=\wt(a)-1.
\end{align}
Since we have chosen $\epsilon(a)=\wt(a)-1$, $\zeta(a)=\wt(a)-1$ in (2).
By (2), the results hold for $a\in A$.
For $b, b^{\prime}\in B$, $\lu\in U$, and $j\geq \wt(b)-1$, by the Borcherds identity (cf. \cite[(3.1.7)]{LL}
putting 
$u=b, v=b^{\prime}, l=-1, m=j+1$, and $n=\epsilon(b^{\prime})$ in the symbol used there) we have
\begin{align}
b_{j}b^{\prime}_{\epsilon(b^{\prime})}\lu
& = \sum\limits_{i=0}^{\infty}
\binom{j+1}{i}(b_{-1+i}b^{\prime})_{j+1+\epsilon(b^{\prime})-i}\lu.
\end{align}
Since
\begin{align}
\label{eq:delta(wt(c))&=max}
\delta(\wt(b_{-1+i}b^{\prime}))&=\rho +\wt(b)+\wt(b^{\prime})-i-1\nonumber\\
&=\epsilon(b^{\prime})-\wt(b^{\prime})+\wt(b)+\wt(b^{\prime})-i\nonumber\\
&=(\wt(b)-1)+1+\epsilon(b^{\prime})-i\nonumber\\
&\leq j+1+\epsilon(b^{\prime})-i,
\end{align}
the results hold for $b\in B$ by (2).
Thus, for any homogeneous element $c\in B_{-}\langle A_{-}\rangle\vac$, an inductive argument on $\wt c$ shows the results.
\end{enumerate}

\end{proof}

We recall the {\it Zhu algebra} $A(V)$ of a vertex operator algebra $V$ from \cite[Section 2]{Z1996}.
For homogeneous $a\in V$ and $b\in V$, we define
\begin{align}
\label{eq:zhu-ideal-multi}
a\circ b&=\sum_{i=0}^{\infty}\binom{\wt a}{i}a_{i-2}b\in V
\end{align}
and 
\begin{align}
\label{eq:zhu-bimodule-left}
a*b&=\sum_{i=0}^{\infty}\binom{\wt a}{i}a_{i-1}b\in V.
\end{align}
We extend \eqref{eq:zhu-ideal-multi} and \eqref{eq:zhu-bimodule-left} for an arbitrary $a\in V$ by linearity.
We also define
$O(V)=\Span_{\C}\{a\circ b\ |\ a,b\in V\}$.
Then, the quotient space
\begin{align}
\label{eq:zhu-bimodule}
A(V)&=M/O(V)\textcolor{black}{,}
\end{align}
called the {\em Zhu algebra} of $V$, is an associative $\C$-algebra with multiplication  
\eqref{eq:zhu-bimodule-left} by \cite[Theorem 2.1.1]{Z1996}.
It is shown in \cite[Theorem 2.2.1]{Z1996} that
for a vertex operator algebra $V$
there is a one to one correspondence between the set of all isomorphism classes of irreducible $\N$-graded weak $V$-modules 
and that of irreducible $A(V)$-modules.

We recall the vertex operator algebra $M(1)$ associated to the Heisenberg algebra and 
the vertex algebra $V_{\lattice}$ associated to a non-degenerate even lattice $\lattice$
from \cite[Sections 6.3--6.5]{LL} and \cite[Section 2.2]{DN2001}.
Let $\hei$ be a finite dimensional vector space equipped with a non-degenerate symmetric bilinear form
$\langle \mbox{ }, \mbox{ }\rangle$.
Set a Lie algebra
\begin{align}
\label{eqn:hathei=heiotimes}
\hat{\hei}&=\hei\otimes \C[t,t^{-1}]\oplus \C \centralc
\end{align} 
with the Lie bracket relations 
\begin{align}
[\beta\otimes t^{m},\gamma\otimes t^{n}]&=m\langle \beta,\gamma\rangle\delta_{m+n,0}\centralc,&
[\centralc,\hat{\hei}]&=0
\end{align}
for $\beta,\gamma\in \hei$ and $m,n\in\Z$.
For $\beta\in \hei$ and $n\in\Z$, $\beta(n)$ denotes $\beta\otimes t^{n}\in\widehat{\fh}$. 
Set two Lie subalgebras of $\fh$:
\begin{align}
\widehat{\hei}_{{\geq 0}}&=\bigoplus_{n\geq 0}\hei \otimes t^n\oplus \C \centralc&\mbox{ and }&&
\widehat{\hei}_{<0}&=\bigoplus_{n\leq -1}\hei\otimes t^n.
\end{align}
For $\beta\in\fh$,
$\C e^{\beta}$ denotes the one dimensional $\widehat{\fh}_{\geq 0}$-module uniquely determined
by the condition that for $\gamma\in\fh$
\begin{align}
\gamma(i)\cdot e^{\beta}&
=\left\{
\begin{array}{ll}
\langle\gamma,\beta\rangle e^{\beta}&\mbox{ for }i=0\\
0&\mbox{ for }i>0
\end{array}
\right.&&\mbox{ and }&
\centralc\cdot e^{\beta}&=e^{\beta}.
\end{align}
We take an $\widehat{\fh}$-module 
\begin{align}
\label{eq:untwist-induced}
M(1,\beta)&=\sU (\widehat{\fh})\otimes_{\sU (\widehat{\fh}_{\geq 0})}\C e^{\beta}
\cong \sU(\widehat{\fh}_{<0})\otimes_{\C}\C e^{\beta}
\end{align}
where $\sU(\fg)$ is the universal enveloping algebra of a Lie algebra $\fg$.
Then, $M(1)=M(1,0)$ has a vertex operator algebra structure with
the conformal vector
\begin{align}
\label{eq:conformal-vector}
\omega&=\dfrac{1}{2}\sum_{i=1}^{\dim\fh}h_i(-1)h_i^{\prime}(-1)\vac
\end{align}
where $\{h_1,\ldots,h_{\dim\fh}\}$ is a basis of $\fh$ and
$\{h_1^{\prime},\ldots,h_{\dim\fh}^{\prime}\}$ is its dual basis.
Moreover, $M(1,\beta)$ is an irreducible $M(1)$-module for any $\beta\in\fh$. 
The vertex operator algebra $M(1)$ is called the {\it vertex operator algebra associated to
 the Heisenberg algebra} $\oplus_{0\neq n\in\Z}\fh\otimes t^{n}\oplus \C \centralc$. 

Let $\lattice$ be a non-degenerate even lattice.
We define $\fh=\C\otimes_{\Z}\lattice$
and denote by $\lattice^{\perp}$ the dual of $\lattice$: $\lattice^{\perp}=\{\gamma\in\fh\ |\ \langle\beta,\gamma\rangle\in\Z\mbox{ for all }\beta\in\lattice\}$. 
Taking $M(1)$ for $\fh$,
we define $V_{\lambda+\lattice}=\oplus_{\beta\in\lambda+\lattice}M(1,\beta)
\cong \sU(\widehat{\hei}_{<0})\otimes_{\C}(\oplus_{\beta\in\lambda+\lattice}\C e^{\beta})$
for $\lambda+\lattice\in \lattice^{\perp}/\lattice$.
Then, $V_{\lattice}=V_{0+\lattice}$ admits a unique vertex algebra structure compatible with the action of $M(1)$ and is called 
 the {\it vertex algebra associated to $\lattice$} (cf. \cite[Section 6.5]{LL}). 
Moreover, for each $\lambda+\lattice\in\lattice^{\perp}/\lattice$
the vector space
$V_{\lambda+\lattice}$ is an irreducible weak $V_{\lattice}$-module which admits the following decomposition:
\begin{align}
V_{\lambda+\lattice}&=\bigoplus_{n\in\langle\lambda,\lambda\rangle/2+\Z}(V_{\lambda+\lattice})_{n} \mbox{ where }
(V_{\lambda+\lattice})_{n}=\{a\in V_{\lambda+\lattice}\ |\ \omega_{1}a=na\}.
\end{align}

Let  $\hat{\lattice}$ be the canonical central extension of $\lattice$
by the cyclic group $\langle\kappa\rangle$ of order $2$ with  the commutator map
$c(\alpha,\beta)=\kappa^{\langle\alpha,\beta\rangle}$ for $\alpha,\beta\in\lattice$:
\begin{align}
	0\rightarrow \langle\kappa\rangle\overset{}{\rightarrow} \hat{\lattice}\overset{-}{\rightarrow} \lattice\rightarrow 0.
\end{align}
Then, the $-1$-isometry of $\lattice$ induces an automorphism $\theta$ of $\hat{\lattice}$ of order $2$
and an automorphism, by abuse of notation we also denote by $\theta$, of $V_{\lattice}$ of order $2$ (see \cite[(8.9.22)]{FLM}).
In $M(1)$, we have
\begin{align}
\label{eq:theta}
\theta(h^1(-i_1)\cdots h^n(-i_n)\vac)&=(-1)^{n}h^1(-i_1)\cdots h^n(-i_n)\vac
\end{align}
for $n\in\Z_{\geq 0}$, $h^1,\ldots,h^{n}\in\fh$, and $i_1,\ldots,i_n\in\Z_{>0}$.
We set
\begin{align}
V_{L}^{+}&=\{a\in V_{L}\ |\ \theta(a)=a\}
\mbox{ and  }M(1)^{+}=\{a\in M(1)\ |\ \theta(a)=a\}.
\end{align}
For a weak $V_{\lattice}$-module $\module$,
we define a weak $V_{\lattice}$-module $(\module\circ \theta,Y_{\module\circ \theta})$
by $\module\circ\theta=\module$ and 
\begin{align}
Y_{\module\circ \theta}(a,x)&=Y_{\module}(\theta(a),x)
\end{align}
for  $a\in V_{\lattice}$.
Then
$V_{\lambda+\lattice}\circ\theta\cong V_{-\lambda+\lattice}$
for $\lambda\in \lattice^{\perp}$.
Thus, for $\lambda\in \lattice^{\perp}$ with $2\lambda\in\lattice$
we define
\begin{align}
V_{\lambda+\lattice}^{\pm}&=\{u\in V_{\lambda+\lattice}\ |\ \theta(u)=\pm u\}.
\end{align}
Next, we recall the construction of $\theta$-twisted modules for $M(1)$ and $V_{\lattice}$ from \cite[Section 9]{FLM}.
Set a Lie algebra
\begin{align}
\hat{\hei}[-1]&=\hei\otimes t^{1/2}\C[t,t^{-1}]\oplus \C \centralc
\end{align} 
with the Lie bracket relations 
\begin{align}
[\centralc,\hat{\hei}[-1]]&=0&\mbox{and}&&
[\alpha\otimes t^{m},\beta\otimes t^{n}]&=m\langle\alpha,\beta\rangle\delta_{m+n,0}\centralc
\end{align}
for $\alpha,\beta\in \hei$ and $m,n\in1/2+\Z$.
For $\alpha\in \hei$ and $n\in1/2+\Z$, $\alpha(n)$ denotes $\alpha\otimes t^{n}\in\widehat{\hei}$. 
Set two Lie subalgebras of $\hat{\hei}[-1]$:
\begin{align}
\widehat{\hei}[-1]_{{\geq 0}}&=\bigoplus_{n\in 1/2+\N}\hei\otimes t^n\oplus \C \centralc&\mbox{ and }&&
\widehat{\hei}[-1]_{{<0}}&=\bigoplus_{n\in 1/2+\N}\hei\otimes t^{-n}.
\end{align}
Let $\C \vac_{\tw}$ denote a unique one dimensional $\widehat{\hei}[-1]_{{\geq 0}}$-module 
such that 
\begin{align}
h(i)\cdot \vac_{\tw}&
=0\quad\mbox{ for }h\in\fh\mbox{ and }i\in \frac{1}{2}+\N,\nonumber\\
\centralc\cdot \vac_{\tw}&=\vac_{\tw}.
\end{align}
We take an $\widehat{\hei}[-1]$-module 
\begin{align}
\label{eq:twist-induced}
M(1)(\theta)
&=\sU (\widehat{\hei}[-1])\otimes_{\sU (\widehat{\hei}[-1]_{\geq 0})}\C u_{\ul{\zeta}}
\cong\sU (\widehat{\hei}[-1]_{<0})\otimes_{\C}\C u_{\ul{\zeta}}.
\end{align}
We define for $\alpha\in \hei$, 
\begin{align}
	\alpha(x)&=\sum_{i\in 1/2+\Z}\alpha(i)x^{-i-1}
	\end{align}
and for $u=\alpha_1(-\wi_1)\cdots \alpha_k(-\wi_k)\vac\in M(1)$, 
\begin{align}
	Y_{0}(u,x)&=\nor
\dfrac{1}{(\wi_1-1)!}	(\dfrac{d^{\wi_1-1}}{dx^{\wi_1-1}}\alpha_1(x))
	\cdots
\dfrac{1}{(\wi_k-1)!}	(\dfrac{d^{\wi_k-1}}{dx^{\wi_k-1}}\alpha_k(x))\nor.
\end{align}
Here, for $\beta_1,\ldots,\beta_{n}\in \fh$ and $i_1,\ldots,i_n\in1/2+\Z$, we define 
$\nor \beta_1(i_1)\cdots\beta_{n}(i_n)\nor$ inductively by
\begin{align}
\label{eq:nomal-ordering}
\nor \beta_1(i_1)\nor&=\beta_1(i_1)\qquad\mbox{ and}\nonumber\\
\nor \beta_1(i_1)\cdots\beta_{n}(i_n)\nor&=
\left\{
\begin{array}{ll}
\nor \beta_{2}(i_2)\cdots\beta_{n}(i_n)\nor \beta_1(i_1)&\mbox{if }i_1\geq 0,\\
\beta_{1}(i_1)\nor \beta_{2}(i_2)\cdots\beta_{n}(i_n)\nor &\mbox{if }i_1<0.
\end{array}\right.
\end{align}
Let $h^{[1]},\ldots,h^{[\dim\fh]}$ be an orthonormal basis of $\fh$.
We define $c_{mn}\in\Q$ for $ m,n\in \Z_{\geq 0}$ by
\begin{align}
	\sum_{m,n=0}^{\infty}c_{mn}x^{m}y^{n}&=-\log (\dfrac{(1+x)^{1/2}+(1+y)^{1/2}}{2})
	\end{align}
and
	\begin{align}
		\Delta_{x}&=\sum_{m,n=0}^{\infty}c_{mn}\sum_{i=1}^{\dim\fh}h^{[i]}(m)h^{[i]}(n)x^{-m-n}.
		\end{align}
Then, for $u\in M(1)$ we define a vertex operator $Y_{M(1)(\theta)}$ by
\begin{align}
Y_{M(1)(\theta)}(u,x)&=Y_{0}(e^{\Delta_{x}}u,x).
\end{align}
Then, \cite[Theorem 9.3.1]{FLM}
shows that 
$(M(1)(\theta),Y_{M(1)(\theta)})$ is an irreducible $\theta$-twisted $M(1)$-module.
We define the action of $\theta$ on $M(1)(\theta)$ by
\begin{align}
\theta(h^1(-i_1)\cdots h^n(-i_n)\vac_{\tw}&=(-1)^{n}h^1(-i_1)\cdots h^n(-i_n)\vac_{\tw}
\end{align}
for $n\in\Z_{\geq 0}$, $h^1,\ldots,h^{n}\in\fh$, $i_1,\ldots,i_n\in 1/2+\Z_{>0}$
and set
\begin{align}
\label{eqn:M(1)(theta)pm=u}
M(1)(\theta)^{\pm}&=\{u\in M(1)(\theta)\ |\ \theta \lu=\pm \lu\}.
\end{align}
Set a submodule $\subL=\{\theta(a) a^{-1}\ |\ a\in\hat{\lattice}\}$ of $\hat{\lattice}$.
Let $T_{\chi}$ be the irreducible $\hat{L}/\subL$-module associated to a
central character $\chi$ such that $\chi(\kappa)=-1$.
We set
\begin{align}
V_{\lattice}^{T_{\chi}}&=M(1)(\theta)\otimes T_{\chi}.
\end{align} 
Then, \cite[Theorem 9.5.3]{FLM} shows that 
$V_{\lattice}^{T_{\chi}}$ admits an irreducible $\theta$-twisted $V_{\lattice}$-module structure compatible with the action of $M(1)$.
We define the action of $\theta$ on $V_{\lattice}^{T_{\chi}}$ by
\begin{align}
\theta(h^1(-i_1)\cdots h^n(-i_n)\lu)&=(-1)^{n}h^1(-i_1)\cdots h^n(-i_n)\lu
\end{align}
for $n\in\Z_{\geq 0}$, $h^1,\ldots,h^{n}\in\fh$, $i_1,\ldots,i_n\in 1/2+\Z_{>0}$, and $\lu\in T_{\chi}$
and set
\begin{align}
V_{\lattice}^{T_{\chi},\pm}&=\{\lu\in V_{\lattice}^{T_{\chi}}\ |\ \theta(u)=\pm u\}.
\end{align}

\textcolor{black}{
Let $h^{[1]},\ldots,h^{[\rankL]}$ be an orthonormal basis of $\fh$.
For $i=1,\ldots, \rankL$,
we define the following elements in $M(1)^{+}$:
\begin{align}
\label{eq:def-oega-i-H-i}
\omega^{[i]}&=\dfrac{1}{2}h^{[i]}(-1)^2\vac,\nonumber\\
\omega&=\omega^{[1]}+\cdots+\omega^{[\rankL]},\nonumber\\
\Har^{[i]}&=\dfrac{1}{3}h^{[i]}(-3)h^{[i]}(-1)\vac-\dfrac{1}{3}h^{[i]}(-2)^2\vac,\nonumber\\
J^{[i]}&=h^{[i]}(-1)^4\vac-2h^{[i]}(-3)h^{[i]}(-1)\vac+\dfrac{3}{2}h^{[i]}(-2)^2\vac\nonumber\\
&=-9\Har^{[i]}+4(\omega^{[i]}_{-1})^2\vac-3\omega^{[i]}_{-3}\vac.
\end{align}
}
For $\alpha\in \fh$, we define
\begin{align}
E(\alpha)&=e^{\alpha}+\theta(e^{\alpha}).
\end{align}
We recall the following notation and some results from \textcolor{black}{\cite[Sections 4 and 5]{DN2001}}:
for any pair of distinct elements $i,j\in\{1,\ldots,\rankL\}$ and $r,s\in\Z_{>0}$,
\begin{align}
\label{eq:def-nS_ij(l,m)}
\nS_{ij}(r,s)&=h^{[i]}(-r)h^{[j]}(-s)\vac,\nonumber\\
E^{u}_{ij}&=5\nS_{ij}(1,2)+25\nS_{ij}(1,3)+36\nS_{ij}(1,4)+16\nS_{ij}(1,5),\nonumber\\
E^{t}_{ij}&=-16\nS_{ij}(1,2)+145\nS_{ij}(1,3)+19\nS_{ij}(1,4)+8\nS_{ij}(1,5),\nonumber\\
\Lambda_{ij}
&=45\nS_{ij}(1,2)+190\nS_{ij}(1,3)+240\nS_{ij}(1,4)+96\nS_{ij}(1,5).
\end{align}
It follows from \cite[Proposition 5.3.14]{DN2001} that in the Zhu algebra $(A(M(1)^{+}),*)$, $A^{u}=\oplus_{i,j}\C E^{u}_{ij}$ and 
$A^{t}=\oplus_{i,j}\C E^{t}_{ij}$ are two-sided ideals, each of  which is isomorphic to the $\rankL\times \rankL$ matrix algebra
and $A^{u}*A^{t}=A^{t}*A^{u}=0$. 
By \cite[Proposition 5.3.12]{DN2001}, for any pair of distinct elements $i,j\in\{1,\ldots,\rankL\}$,
we have $A^{u}*\Lambda_{ij}=\Lambda_{ij}*A^{u}=A^{t}*\Lambda_{ij}=\Lambda_{ij}*A^{t}=0$.
By \cite[Proposition 5.3.15]{DN2001},
$A(M(1)^{+})/(A^{u}+A^{t})$ is a commutative algebra generated by
the images of $\omega^{[i]},\Har^{[i]}$ and $\Lambda_{jk}$ where
$i=1,\ldots,\rankL$ and $j,k\in \{1,\ldots,\rankL\}$ with $j\neq k$.

For $\lambda\in\fh$, $k=1,\ldots,\rankL$, and any pair of distinct elements $i,j\in\{1,\ldots,\rankL\}$,
\begin{align}
\label{eq:norm2Si1(1,1)1lu=-Si1(1,2)2lu}
\omega^{[k]}_{1}e^{\lambda}
&=\frac{\langle\lambda,h^{[k]}\rangle^2}{2}e^{\lambda},\nonumber\\
\Har^{[k]}_{3}e^{\lambda}&=0,\nonumber\\
S_{ij}(1,1)_{1}e^{\lambda}&=-S_{ij}(1,2)_{2}e^{\lambda}
=S_{ij}(1,3)_{3}e^{\lambda}=\langle\lambda,h^{[i]}\rangle\langle\lambda,h^{[j]}\rangle e^{\lambda},
\end{align}
\begin{align}
\label{eq:norm2Si1(1,1)1lu=-Si1(1,2)2lu-2}
\omega^{[k]}_{1}h^{[j]}(-1)\vac
&=\delta_{jk}h^{[j]}(-1)\vac,\nonumber\\
\Har^{[k]}_{3}h^{[j]}(-1)\vac&=\delta_{jk}h^{[j]}(-1)\vac,\nonumber\\
S_{ij}(1,1)_{1}h^{[j]}(-1)\vac&=h^{[i]}(-1)\vac,\nonumber\\
S_{ij}(1,2)_{2}h^{[j]}(-1)\vac&=-2h^{[i]}(-1)\vac,\nonumber\\
S_{ij}(1,3)_{3}h^{[j]}(-1)\vac&=3h^{[i]}(-1)\vac,
\end{align}
and 
\begin{align}
\label{eq:norm2Si1(1,1)1lu=-Si1(1,2)2lu-twist-0}
\textcolor{black}{\omega^{[k]}_{1}\vac_{\tw}}&=\frac{1}{16} \vac_{\tw},\nonumber\\
\Har^{[k]}_{3} \vac_{\tw}&=\frac{-1}{128} \vac_{\tw},\nonumber\\
S_{ij}(1,1)_{1} \vac_{\tw}&=
S_{ij}(1,2)_{2} \vac_{\tw}=
S_{ij}(1,3)_{3} \vac_{\tw}=0,\nonumber\\
\omega^{[k]}_{1}h^{[j]}(-\frac{1}{2})\vac_{\tw}&=\delta_{jk}\frac{9}{16}h^{[j]}(-\frac{1}{2})\vac_{\tw},\nonumber\\
\Har^{[k]}_{3}h^{[j]}(-\frac{1}{2})\vac_{\tw}&=\delta_{jk}\frac{15}{128}h^{[j]}(-\frac{1}{2})\vac_{\tw},\nonumber\\
S_{ij}(1,1)_{1}h^{[j]}(-\frac{1}{2})\vac_{\tw}&=\frac{1}{2}h^{[i]}(-\frac{1}{2})\vac_{\tw},\nonumber\\
S_{ij}(1,2)_{2}h^{[j]}(-\frac{1}{2})\vac_{\tw}&=\frac{-3}{4}h^{[i]}(-\frac{1}{2})\vac_{\tw},\nonumber\\
S_{ij}(1,3)_{3}h^{[j]}(-\frac{1}{2})\vac_{\tw}&=\frac{15}{16}h^{[i]}(-\frac{1}{2})\vac_{\tw}.
\end{align}
For any pair of distinct elements $i,j,k\in\{1,\ldots,\rankL\}$, $l,m\in\Z$, and $r,s\in\Z_{>0}$,
a direct computation shows that
\begin{align}
	\label{eq:[h[j](l),Sij(1,1)m]}
	\textcolor{black}{[h^{[j]}(l),S_{ij}(r,s)_{m}]}&=s\binom{l}{s}\binom{-l-m+r+s-2}{r-1}h^{[i]}(l+m-r-s+1),\\
[\omega^{[j]}_l,S_{ij}(1,r)_{m}]
&=rS_{ij}(1,r+1)_{m+l}+lrS_{ij}(1,r)_{m+l-1}, \label{eq:omega-s11-2}\\
[S_{kj}(1,1)_{l},S_{ij}(1,r)_{m}]
&=r\sum_{t=1}^{r+1}\binom{l}{t}S_{ik}(1,t)_{l+m-r-1+t}.
\label{eq:omega-s11}
\end{align}
We also have
\begin{align}
\label{eq:Sij(2,1)=omega_0S_ij(1,1)-S_ij(1,2)}
S_{ij}(2,1)&=
\omega_{0 } S_{ij}(1,1)
-S_{ij}(1,2),
\nonumber\\
S_{ij}(3,1)
&=\frac{1}{2}\omega_{0 }^{2}S_{ij}(1,1)
-\omega_{0 } S_{ij}(1,2)
+S_{ij}(1,3),\nonumber\\
S_{ij}(2,2)
&=
\omega_{0 } S_{ij}(1,2)
-2S_{ij}(1,3),\nonumber\\
S_{ij}(3,2)
&=-\omega^{[j]}_{-2 } S_{ij}(1,1)
+2\omega^{[j]}_{-1 } S_{ij}(1,2)
+\frac{1}{2}\omega_{0 }^{2}S_{ij}(1,2)
-2\omega_{0 } S_{ij}(1,3),\nonumber\\
S_{ij}(3,3)
&=
\frac{-1}{2}\omega_{0 } \omega^{[j]}_{-2 } S_{ij}(1,1)
+\frac{3}{2}\omega^{[i]}_{-2 } S_{ij}(1,2)
\nonumber\\&\quad{}
-\omega_{0 } \omega^{[i]}_{-1 } S_{ij}(1,2)
+\omega_{0 } \omega^{[j]}_{-1 } S_{ij}(1,2)
\nonumber\\&\quad{}
+\frac{1}{4}\omega_{0 }^{3}S_{ij}(1,2)
+2\omega^{[i]}_{-1 } S_{ij}(1,3)
-\omega_{0 }^{2}S_{ij}(1,3).
\end{align}
For $P\subset \{1,\ldots,\rankL\}$,
we define the subspace
\begin{align}
\label{eq:definitionM1P}
M(1)_{P}&:=\Span_{\C}\Big\{
h^{[i_1]}(-j_1)\ldots h^{[i_n]}(-j_n)\vac\ \Big|\ 
\begin{array}{l}
n\in \Z_{\geq 0}, i_1,\ldots,i_n\in\{1,\ldots,\rankL\}, j_1,\ldots,j_n\in\Z_{>0},\\ 
\big\{l\in\{1,\ldots,n\}\ \big|\ |\{k\ |\ i_k=l\}|\mbox{ is  odd}\big\}=P
\end{array}
\Big\}
\end{align}
of $M(1)$ and the subspace 
\begin{align}
\label{eq:definitionM1P-2}
M(1)^{+}_{P}&:=M(1)_{P}\cap M(1)^{+}
\end{align}
of $M(1)^{+}$.
Note that if $|P|$ is odd, then $M(1)^{+}_{P}=\{0\}$.

For $P, P^{\prime}\subset \{1,\ldots,\rankL\}$, we define 
$P\ominus P^{\prime}:=(P\cup P^{\prime})\setminus(P\cap P^{\prime})\subset\{1,\ldots,\rankL\}$.
\begin{lemma}
\mbox{}
\label{lemma:(M(1)GammaPM(1)Gammaprime}
\begin{enumerate}
\item
We have
$M(1)^{+}=\bigcup_{\begin{subarray}{l}P\subset\{1,\ldots,\rankL\}\\|P|\mmbox{ is even}\end{subarray}}M(1)^{+}_{P}$.
\item
For $P,P^{\prime}\subset \{1,\ldots,\rankL\}$,
\begin{align}
\label{eq:M(1)GammaPM(1)Gammaprime}
(M(1)_{P})\cdot M(1)_{P^{\prime}}
&\subset M(1)_{P\ominus P^{\prime}}
\mbox{ and }\nonumber\\
(M(1)^{+}_{P})\cdot M(1)^{+}_{P^{\prime}}
&\subset M(1)^{+}_{P\ominus P^{\prime}}.
\end{align}
\end{enumerate}
\end{lemma}
\begin{proof}
The result (1) follows from the definition of $M(1)^{+}_{P}$.
The result (2) follows from the fact that
for $i\in\{1,\ldots,\rankL\}$ and $j\in\Z$, 
$h^{[i]}(j)M(1)_{P^{\prime}}\subset 
M(1)_{\{i\}\ominus P^{\prime}}$.
\end{proof}
For $P=\{p_1,\ldots,p_{2t}\}\subset \{1,\ldots,\rankL\}$ with $p_1>\cdots>p_{2t}$, we define
\begin{align}
B_{P}&:=\{
h^{[p_1]}(-1)h^{[p_2]}(-r_1)\cdots h^{[p_{2t-1}]}(-1)h^{[p_{2t}]}(-r_t)\vac\ |\ 
r_1,\ldots,r_t\in\{1,2,3\}
\}\nonumber\\
&=\{
S_{p_1,p_2}(1,r_1)_{-1}\cdots S_{p_{2t-1},p_{2t}}(1,r_{t})_{-1}\vac\ |\ 
r_1,\ldots,r_t\in\{1,2,3\}
\}.
\end{align}
By \eqref{eq:norm2Si1(1,1)1lu=-Si1(1,2)2lu-2}, \eqref{eq:[h[j](l),Sij(1,1)m]}, 
\eqref{eq:Sij(2,1)=omega_0S_ij(1,1)-S_ij(1,2)}, we know that
\begin{align}
\label{eq]M(1)+P=M(1)+}
M(1)^{+}_{P}&=M(1)^{+}_{\varnothing}\cdot (B_{P})_{-}\vac
\end{align}
and hence by \eqref{eq:[h[j](l),Sij(1,1)m]} again,
$M(1)^{+}_{P}$ is also spanned by the elements of the form
\begin{align}
\label{eq:Sp1,p2(1,r1)-s1cdotsSp}
S_{p_1,p_2}(1,r_1)_{-s_1}\cdots S_{p_{2t-1},p_{2t}}(1,r_{t})_{-s_{t}}a,
\end{align}
where  $a\in M(1)^{+}_{\varnothing}$, $s_1,\ldots,s_{t}\in \Z_{>0}$, and $r_1,\ldots,r_t\in\{1,2,3\}$.
By \eqref{eq:norm2Si1(1,1)1lu=-Si1(1,2)2lu-2}, \eqref{eq:[h[j](l),Sij(1,1)m]}, 
\eqref{eq:Sij(2,1)=omega_0S_ij(1,1)-S_ij(1,2)},  and \cite[(3.4), (3.5), (3.7), (3.9), (3.10), (3.11)]{Tanabe2021-1},
$M(1)^{+}_{P}$ is spanned by the elements of the form
\begin{align}
\label{eq:a(1)(-l1)ldots a(m)lm}
a^{(1)}_{-l_1}\ldots a^{(m)}_{-l_m}S_{p_1,p_2}(1,r_1)_{-s_1}\cdots S_{p_{2t-1},p_{2t}}(1,r_t)_{-s_t}\vac,
\end{align}
where $m\in\Z_{\geq 0}$, $a^{(1)},\ldots,a^{(m)}\in\{\omega^{[j]},\Har^{[j]}\ |\ j\in\{1,\ldots,\rankL\}\}$,
$l_1,\ldots,l_m,s_1,\ldots,s_t\in\Z_{>0}$, and $r_1,\ldots,r_t\in\{1,2,3\}$.
By \eqref{eq:Sp1,p2(1,r1)-s1cdotsSp}, in the same way, we know that
$M(1)^{+}_{P}$ is spanned by the elements of the form
\begin{align}
\label{eq:Sp1p2(1r1)-s1cdots}
S_{p_1,p_2}(1,r_1)_{-s_1}\cdots S_{p_{2t-1},p_{2t}}(1,r_t)_{-s_t}a^{(1)}_{-l_1}\ldots a^{(m)}_{-l_m}\vac,
\end{align}
where $m\in\Z_{\geq 0}$, $a^{(1)},\ldots,a^{(m)}\in\{\omega^{[j]},\Har^{[j]}\ |\ j\in\{1,\ldots,\rankL\}\}$,
$l_1,\ldots,l_m,s_1,\ldots,s_t\in\Z_{>0}$, and $r_1,\ldots,r_t\in\{1,2,3\}$.

\begin{remark}
\label{remark:com-appendix}
\textcolor{black}{Let $u,v\in M(1)^{+}$, $i,j\in\Z,\rho,\sigma\in\Z_{\geq -1}$, and $p,q$ a pair} of distinct elements in $\{1,\ldots,\rankL\}$.
Throughout this paper, 
if $[u_{i},v_{j}]=\sum_{k=0}^{\infty}\binom{i}{k}(u_{k}v)_{i+j-k}\in M(1)^{+}_{\{p,q\}}$,
then we frequently express this element 
as a linear combination of elements of the form
\begin{align}
a^{(1)}_{i_1}\cdots a^{(k)}_{i_{k}}S_{pq}(1,r)_{t}b^{(1)}_{j_1}\cdots b^{(l)}_{j_{l}}
	\end{align}
where $k,l\in\Z_{\geq 0}$, $r=1,2,3$, $t\in\Z$, and 
\begin{align}
(a^{(1)},i_1),\ldots,(a^{(k)},i_{k})&\in\{(\omega^{[k]},m)\ |\ m\leq \rho\}_{k=1}^{\rankL}\cup\{(\Har^{[k]},n)\ |\ n\leq \sigma\}_{k=1}^{\rankL},\nonumber\\
(b^{(1)},j_1),\ldots,(b^{(l)},j_{l})&\in\{(\omega^{[k]},m)\ |\ m\geq \rho+1\}_{k=1}^{\rankL}\cup\{(\Har^{[k]},n)\ |\ n\geq \sigma+1\}_{k=1}^{\rankL}.
\end{align}
For the calculation, we use \cite[Lemma 2.2]{Tanabe2021-1} and the data $a_{k}b\ (k=0,1,\ldots)$ in  \ref{section:appendix}.
In most cases, we obtain the explicit expressions of the results by using computer algebra system Risa/Asir\cite{Risa/Asir}.

For example, we shall compute $[	\Har^{[j]}_{l}, S_{ij}(1,1)_{n}]$ for a pair of distinct elements $i,j\in\{1,\ldots,\rankL\}$
and $l,n\in\Z$.
For $m\in\Z_{\geq -1}$, by \cite[Lemma 2.2]{Tanabe2021-1}, 
we have
\begin{align}
\label{eq:(omega[j]_-1S_ij(1,1))_n}
(\omega^{[j]}_{-2 } S_{ij}(1,1))_{n}
&=
\sum_{\begin{subarray}{l}r\leq m\\r+s+2=n\end{subarray}}(-r-1)\omega^{[j]}_{r} S_{ij}(1,1)_{s}+
\sum_{\begin{subarray}{l}r\geq m+1\\r+s+2=n\end{subarray}}(-r-1)S_{ij}(1,1)_{s}\omega^{[j]}_{r} \nonumber\\
&\quad{}+
\sum_{t=0}^{1}(t+1)\binom{m+2}{t+2}(\omega^{[j]}_{t}S_{ij}(1,1))_{n-2-t}\nonumber\\
&=
\sum_{\begin{subarray}{l}r\leq m\\r+s+2=n\end{subarray}}(-r-1)\omega^{[j]}_{r} S_{ij}(1,1)_{s}+
\sum_{\begin{subarray}{l}r\geq m+1\\r+s+2=n\end{subarray}}(-r-1)S_{ij}(1,1)_{s}\omega^{[j]}_{r} \nonumber\\
&\quad{}+
\binom{m+2}{2}S_{ij}(1,2)_{n-2}
+2\binom{m+2}{3}S_{ij}(1,1)_{n-3},\nonumber\\
(\omega^{[j]}_{-1} S_{ij}(1,2))_{n}
&=
\sum_{\begin{subarray}{l}r\leq m\\r+s+1=n\end{subarray}}\omega^{[j]}_{r} S_{ij}(1,2)_{s}+
\sum_{\begin{subarray}{l}r\geq m+1\\r+s+1=n\end{subarray}}S_{ij}(1,2)_{s}\omega^{[j]}_{r} \nonumber\\
&\quad{}-
\sum_{t=0}^{1}\binom{t}{0}\binom{m+1}{t+1}(\omega^{[j]}_{t}S_{ij}(1,2))_{n-1-t}\nonumber\\
&=
\sum_{\begin{subarray}{l}r\leq m\\r+s+1=n\end{subarray}}\omega^{[j]}_{r} S_{ij}(1,2)_{s}+
\sum_{\begin{subarray}{l}r\geq m+1\\r+s+1=n\end{subarray}}S_{ij}(1,2)_{s}\omega^{[j]}_{r} \nonumber\\
&\quad{}-2\binom{m+1}{1}S_{ij}(1,3)_{n-1}
-2\binom{m+1}{2}S_{ij}(1,2)_{n-2}.
\end{align}
By \eqref{eq:Har[j]0Sij(1,1)=-2omega[j]-2Sij(1,1)-1vac} and \eqref{eq:(omega[j]_-1S_ij(1,1))_n}, we have
\begin{align}
&[	\Har^{[j]}_{l}, S_{ij}(1,1)_{n}]
=\sum_{k=0}^{3}\binom{l}{k}(\Har^{[j]}_{k}S_{ij}(1,1))_{l+n-k}\nonumber\\
&=
	(-2\omega^{[j]}_{-2 } S_{ij}(1,1)	+4\omega^{[j]}_{-1 } S_{ij}(1,2))_{l+n}\nonumber\\
&\quad{}+4lS_{ij}(1,3)_{l+n-1}+\binom{l}{2}	\frac{7}{3}S_{ij}(1,2)_{l+n-2}
+\binom{l}{3}S_{ij}(1,1)_{l+n-3}\nonumber\\
&=-2\big(\sum_{\begin{subarray}{l}r\leq m\\r+s+2=l+n\end{subarray}}(-r-1)\omega^{[j]}_{r} S_{ij}(1,1)_{s}+
\sum_{\begin{subarray}{l}r\geq m+1\\r+s+2=l+n\end{subarray}}(-r-1)S_{ij}(1,1)_{s}\omega^{[j]}_{r} \nonumber\\
&\qquad{}+
\binom{m+2}{2}S_{ij}(1,2)_{l+n-2}
+2\binom{m+2}{3}S_{ij}(1,1)_{l+n-3}\big)\nonumber\\
&\quad{}+4\big(\sum_{\begin{subarray}{l}r\leq m\\r+s+1=l+n\end{subarray}}\omega^{[j]}_{r} S_{ij}(1,2)_{s}+
\sum_{\begin{subarray}{l}r\geq m+1\\r+s+1=l+n\end{subarray}}S_{ij}(1,2)_{s}\omega^{[j]}_{r}\nonumber\\
&\qquad{}-2\binom{m+1}{1}S_{ij}(1,3)_{l+n-1}-2\binom{m+1}{2}S_{ij}(1,2)_{l+n-2}\big)\nonumber\\
&\quad{}+4lS_{ij}(1,3)_{l+n-1}+\binom{l}{2}	\frac{7}{3}S_{ij}(1,2)_{l+n-2}
+\binom{l}{3}S_{ij}(1,1)_{l+n-3}.
\end{align}
\end{remark}

\section{Modules for the Zhu algebra of $M(1)^{+}$ in a weak $V_{\lattice}^{+}$-module: the general case}
\label{section:Modules for the Zhu algera of general}
Let $L$ be a non-degenerate even lattice of finite rank $d$ and $\fh:=\C\otimes_{\Z}L$.
In this section, we shall show that there exists an irreducible \textcolor{black}{$M(1)^{+}$}-module in any non-zero 
weak $V_{\lattice}^{+}$-module (Proposition \ref{proposition:M(1)plusmoduleinweak}).
Throughout this section $M$ is a weak $V_{\lattice}^{+}$-module.

\begin{lemma}
For a non-degenerate even lattice $\lattice$ of finite rank $\rankL$,
there exists a sequence 
	of elements $\beta_1,\ldots,\beta_{\rankL}\in \lattice$ such that
	$\langle\beta_i,\beta_i\rangle\neq 0$ for $i=1,\ldots,\rankL$ and 
	$\langle\beta_j,\beta_k\rangle=0$ 
for any pair of distinct elements $j,k\in\{1,\ldots,\rankL\}$. 
\end{lemma}
\begin{proof}
	Let $\gamma_1,\ldots,\gamma_d$ be a basis of $\Q\otimes_{\Z}\lattice$
	such that $\langle\gamma_i,\gamma_i\rangle\neq 0$ and 
	$\langle\gamma_j,\gamma_k\rangle=0$ for all $i\in\{1,\ldots,\rankL\}$
	and $j,k\in\{1,\ldots,\rankL\}$ with $j\neq k$. . 
	Since $\gamma_1,\ldots,\gamma_d\in \Q\otimes_{\Z}\lattice$,
	there exists a non-zero integer $m_i$ such that 
	$m\gamma_i\in \lattice$ for all $i=1,\ldots,d$.
	Then, the elements $\beta_i=m\gamma_i\ (i=1,\ldots,d)$ satisfy the condition.
\end{proof}

Let $\Lambda=\oplus_{i=1}^{d}\Z\beta_i$ be a sublattice of $\lattice$
such that $\langle\beta_i,\beta_i\rangle\neq 0$ for $i=1,\ldots,\rankL$ and 
$\langle\beta_j,\beta_k\rangle=0$ for any pair of distinct elements $j,k\in\{1,\ldots,\rankL\}$. 
We have
\begin{align}
	V_{\Z\beta_1}^{+}\otimes \cdots\otimes V_{\Z\beta_d}^{+} \subset V_{\lattice}^{+}
\end{align}
and take the orthonormal basis $h^{[1]},\ldots,h^{[\rankL]}$ of $\fh$ defined by
\begin{align}
\label{eq:h[i]=frac1sqrtlanglebeta[i]}
h^{[i]}&=\frac{1}{\sqrt{\langle\beta^{[i]},\beta^{[i]}\rangle}}\beta^{[i]}\quad (i=1,\ldots,\rankL).
\end{align}
Since $[h^{[i]}_{l},h^{[j]}_{m}]=0$ for any pair of distinct  elements $i,j\in\{1,\ldots,\rankL\}$ and $l,m\in\Z$, 
it follows by induction on $d$ \textcolor{black}{using} \cite[Lemma 3.7]{Tanabe2021-1} that there exists a simultaneous eigenvector $\lu$ of
 $\{\omega^{[i]}_{1},\Har^{[i]}_{3}\}_{i=1}^{\rankL}$ in a weak $V_{\lattice}^{+}$-module $\module$
such that $\epsilon(\omega^{[i]},\lu)\leq 1$ and
$\epsilon(\Har^{[i]},\lu)\leq 3$ for all $i=1,\ldots,\rankL$.

\begin{lemma}
	\label{lemma:bound-index-S}
	Let $U$ be a subspace of a weak $M(1)^{+}$-module.
	\begin{enumerate}
		\item
		Let $i,j\in\{1,\ldots,\rankL\}$ with $i\neq j$ and
		$\wak\in\Z$ such that
		$\wak\geq\epsilon(S_{ij}(1,1),\lu)$ for all $\lu\in U$.
		If $U$ is stable under the action of $\omega_{1}^{[j]}$, then
		\begin{align}
			\label{eq:sab-bound}
			\epsilon(S_{ij}(1,r+1),\lu)\leq \wak+r
		\end{align}
		for all $r\in\Z_{\geq 0}$.
		\item
		Assume $\alpha\in\C h^{[1]}$.
		Let $i\in\{2,\ldots,\rankL\}$ and
		$\lE\in\Z$ such that
		$\lE\geq\epsilon(\ExB(\alpha),\lu)$ for all $\lu\in U$.
		If $U$ is stable under the action of $S_{i1}(1,1)_{1}$, then
		\begin{align}
			\label{eq:sab-bound-1}
			\epsilon(S_{i1}(1,1)_{0}\ExB(\alpha),\lu)\leq \lE+1.
		\end{align}
		
	\end{enumerate}
	
\end{lemma}
\begin{proof}
	\begin{enumerate}
		\item
		For $l,m\in\Z$ and $r\in\Z_{>0}$, by \eqref{eq:omega-s11-2}
		\begin{align}
			\label{eq:s12s13-bound}
			S_{ij}(1,r+1)_{l+m}&=\dfrac{1}{m}[\omega^{[j]}_l,S_{ij}(1,r)_{m}]-lS_{ij}(1,r)_{l+m-1},
		\end{align}
		which implies \eqref{eq:sab-bound}.
		\item
		Since $S_{i1}(1,1)_{m}\ExB(\alpha)=0$ for all $m\in\Z_{>0}$,
		the same argument as above shows the result.
	\end{enumerate}
\end{proof}

\begin{lemma}
\label{lemma:epsilon(omegailuleq1epsilonHariluleq3}
Let $U$ be a finite dimensional subspace of a weak $M(1)^{+}$-module $\module$ such that
for all $k=1,\ldots,\rankL$ and $\lu\in U$,
$\epsilon(\omega^{[k]},\lu)\leq 1, \epsilon(\Har^{[k]},\lu)\leq 3$,
and $\omega^{[k]}_{1}\lu\in U$, $\Har^{[k]}_{3}\lu\in U$.
For any pair of distinct elements $i,j\in\{1,\ldots,\rankL\}$ we denote 
$\max\{\{\epsilon(S_{ij}(1,1),\lu)\ |\ \lu\in U\}\cup\{-1\}\}$ by $\epsilon(S_{ij})$.
\begin{enumerate}
\item
Let $i,j$ be a pair of distinct elements in $\{1,\ldots,\rankL\}$.
We define $\mW:=\Span_{\C}\{S_{ij}(1,r)_{\epsilon(S_{ij})+r-1}\lu\ |\ r=1,2,3\}$.
For any $\lw\in \mW$ and $k=1,\ldots,\rankL$, we have
$\epsilon(\omega^{[k]},\lw)\leq 1, \epsilon(\Har^{[k]},\lw)\leq 3$
and  $\omega^{[k]}_{1}\lw, \Har^{[k]}_{3}\lw\in \mW$.
\item
Assume $\epsilon(S_{ij})\leq 1$ for any pair of distinct elements $i,j\in\{1,\ldots,\rankL\}$, namely
$U\subset\Omega_{M(1)^{+}}(\module)$.
For \textcolor{black}{$P=\{p_1,\ldots, p_{2t}\}\subset \{1,\ldots,\rankL\}$ with $p_1>\cdots>p_{2t}$,}
we define
\begin{align}
S_{P}U&:=\Span_{\C}\{S_{p_1,p_2}(1,r_1)_{r_1}
\cdots
S_{p_{2t-1},p_{t}}(1,r_{t})_{r_{t}}\lu
\ |\  \lu\in U\mbox{ and }r_1,\ldots,r_{t}\in \{1,2,3\}\}
\end{align}
and  $SU:=\sum_{\begin{subarray}{l}P\subset \{1,\ldots,\rankL\}, \\|P|\mmbox{ is even}\end{subarray}}S_{P}U$.
Then, $SU$ is an $A(M(1)^{+})$-submodule of $\Omega_{M(1)^{+}}(\module)$.
\end{enumerate}
\end{lemma}
\begin{proof}
\begin{enumerate}
\item
We may take $(i,j)$ to be $(2,1)$. We define
 ${\eS}:=\epsilon(S_{21}(1,1))$,
$\epsilon(S_{21}(1,2)):={\eS}+1$, and $\epsilon(S_{21}(1,3)):={\eS}+2$.
By Lemma \ref{lemma:bound-index-S} (1), we have $\epsilon(S_{21}(1,i))\geq \epsilon(S_{21}(1,i),\lu)$ for all $i=1,2,3$.
We define $A:=\{\omega^{[i]},\Har^{[i]}\}_{i=1}^{\rankL}$, $B:=\{S_{21}(1,i)\ |\ i=1,2,3\}$, 
$\epsilon(\omega^{[i]}):=1=\wt(\omega^{[i]})-1$ and $\epsilon(\Har^{[i]}):=3=\wt(\Har^{[i]})-1$
for $i=1,\ldots,\rankL$.
In order to apply Lemma \ref{lemma:ksumklimitsm=1n(kepsilon(a(m))+im)kk+(kepsilon(b)+j)-1+kr} (3) to $\lu$,
for $a\in A$ we shall compute $\zeta(a)$ defined  in  \eqref{eq:zeta(a):=maxepsilon(a)}. 
Note that $M(1)^{+}_{\varnothing}=\langle A_{-}\rangle\vac$ and $M(1)^{+}_{\{2,1\}}
=B_{-}\langle A_{-}\rangle\vac
=\langle A_{-}\rangle B_{-}\vac$ (see \eqref{eq:definitionM1P-2} for the definition of $M(1)^{+}_{P}$). 
By Lemma \ref{lemma:(M(1)GammaPM(1)Gammaprime} (2), we have
$
A\cdot \langle A_{-}\rangle\vac \subset \langle A_{-}\rangle\vac,
A\cdot (\langle A_{-}\rangle B_{-}\vac)\subset \langle A_{-}\rangle B_{-}\vac,
B\cdot \langle A_{-}\rangle\vac \subset B_{-}\langle A_{-}\rangle\vac$,
and 
$A\cdot (B_{-}\langle A_{-}\rangle\vac) \subset B_{-}\langle A_{-}\rangle\vac$,
 where the symbols $\langle A_{-}\rangle\vac, \langle A_{-}\rangle B_{-}\vac$, and $ B_{-}\langle A_{-}\rangle\vac$
are defined in \eqref{eqn:A-B:=Span} and \eqref{eq:a(1)-i1cdotsa(n)-inb}.
For $j\in\Z$ with $j\geq \min\{\wt(b)\ |\ b\in B\}=2$, by the definition \eqref{eq:gamma(k):=maxBiggsum} of $\delta$ , 
\begin{align}
\delta(j)
&=\epsilon(S_{21}(1,1))-2+j.
\end{align}
For $a\in\{\omega^{[i]}_{1},\Har^{[i]}_{3}\}_{i=1}^{\rankL}$ and $j=1,2,3$,
since $\wt(a)\geq 2$, we have
\begin{align}
\label{eq:delta(wt(a)+wt(S211,j))-1)-epsilon(S21(1,j))}
&\delta(\wt(a)+\wt (S_{21}(1,j))-1)-\epsilon(S_{21}(1,j))
=\delta(\wt(a)+j)-(\epsilon(S_{21}(1,1))+j-1)\nonumber\\
&=\epsilon(S_{21}(1,1))-2+(\wt(a)+j)-(\epsilon(S_{21}(1,1))+j-1)=\wt(a)-1
\end{align}
and hence $\zeta(a)=\wt(a)-1$.
Applying Lemma \ref{lemma:ksumklimitsm=1n(kepsilon(a(m))+im)kk+(kepsilon(b)+j)-1+kr} (3) to $\lu$,
we have $a_{\wt(a)-1}U\subset U$ and $a_{k}U=0$ for $a\in A$ and $k>\wt(a)-1$.
\item
For $P\subset\{1,\ldots,\rankL\}$ such that $|P|$ is even,
by using (1), an inductive argument on $|P|$ shows that 
$\epsilon(\omega^{[k]},\lu)\leq 1, \epsilon(\Har^{[k]},\lu)\leq 3$,
and $\omega^{[k]}_{1}\lu\in S_{P}U$, $\Har^{[k]}_{3}\lu\in S_{P}U$
for all $k=1,\ldots,\rankL$ and $\lu\in S_{P}U$.
Let 
$A:=\{\omega^{[i]},\Har^{[i]}\}_{i=1}^{\rankL}$
and 
\begin{align}
B&:=\Big\{
h^{[p_1]}(-1)h^{[p_2]}(-r_1)\cdots h^{[p_{2t-1}]}(-1)h^{[p_{2t}]}(-r_n)\vac\ \Big|\ 
\begin{array}{l}
t\in\Z_{\geq 0}, r_1,\ldots,r_t\in\{1,2,3\}\\
p_1,\ldots,p_{2t}\in\{1,\ldots,\rankL\}\mbox{ such that}\\
p_1>\cdots>p_{2t} 
\end{array}
\Big\}\nonumber\\
&=\Big\{
S_{p_1,p_2}(1,r_1)_{-1}\cdots S_{p_{2t-1},p_{2t}}(1,r_{t})_{-1}\vac\ \Big|\ 
\begin{array}{l}
t\in\Z_{\geq 0}, r_1,\ldots,r_t\in\{1,2,3\}\\
p_1,\ldots,p_{2t}\in\{1,\ldots,\rankL\}\mbox{ such that}\\
p_1>\cdots>p_{2t} 
\end{array}
\Big\}.
\end{align}
By Lemma \ref{lemma:(M(1)GammaPM(1)Gammaprime}, 
\eqref{eq:a(1)(-l1)ldots a(m)lm}, and \eqref{eq:Sp1p2(1r1)-s1cdots}, we have $A\cdot (B_{-}\vac)=\langle A_{-}\rangle B_{-}\vac=B_{-}\langle A_{-}\rangle\vac=M(1)^{+}$.
We define $\epsilon(u):=\wt(u)-1$ for a homogeneous element  $u\in M(1)^{+}$.
Note that $\Span_{\C}\{b_{\epsilon(b)}\lu\ |\ b\in B, \lu\in U\}=SU$.
Since 
$
A\cdot \langle A_{-}\rangle\vac \subset \langle A_{-}\rangle\vac$ and 
$A\cdot M(1)^{+}=B\cdot M(1)^{+}=M(1)^{+}$,
the result follows from Lemma \ref{lemma:ksumklimitsm=1n(kepsilon(a(m))+im)kk+(kepsilon(b)+j)-1+kr} (4).
\end{enumerate}
\end{proof}
\begin{lemma}
\label{lemma:Zhu-Omega}
For a non-zero weak $V_{\lattice}^{+}$-module $\module$,
there exists an irreducible $A(M(1)^{+})$-submodule 
of
 $\Omega_{M(1)^{+}}(M)$.
\end{lemma}
\begin{proof}
Let $h^{[1]},\ldots,h^{[\rankL]}$ be the orthonormal basis of $\fh$
defined by \eqref{eq:h[i]=frac1sqrtlanglebeta[i]}.
By the argument just after \eqref{eq:h[i]=frac1sqrtlanglebeta[i]}, we can take a simultaneous eigenvector $\lu$ of 
	$\{\omega^{[i]}_{1},\Har^{[i]}_{3}\}_{i=1}^{\rankL}$
	such that $\epsilon(\omega^{[i]},\lu)\leq 1$ and
	$\epsilon(\Har^{[i]},\lu)\leq 3$ for all $i=1,\ldots,\rankL$.
We take a pair of distinct elements $i,j$ so that $\epsilon(S_{ij}(1,1),\lu)\geq \epsilon(S_{lm}(1,1),\lu)$
	for any pair of distinct elements $l,m\in\{1,\ldots,\rankL\}$.
We may assume $(i,j)=(2,1)$. 
We define ${\eS}=\epsilon(S_{21}(1,1)):=\epsilon(S_{21}(1,1),\lu)$,
$\epsilon(S_{21}(1,2)):={\eS}+1$, and $\epsilon(S_{21}(1,3)):={\eS}+2$.
By Lemma \ref{lemma:bound-index-S} (1), we have $\epsilon(S_{21}(1,i))\geq \epsilon(S_{21}(1,i),\lu)$ for all $i=1,2,3$.
Hence if ${\eS}\leq 0$, then $u\in \Omega_{M(1)^{+}}(\module)$
and $\C u$ is an irreducible $A(M(1)^{+})$-module.

{From }now, we assume \textcolor{black}{${\eS}\geq 1$}.
We define the subspace $\mW:=\sum_{r=1}^{3}\C S_{21}(1,r)_{{\eS}+r-1}\lu$ of $\module$.
Let $\lw\in\mW$ and $(j,k)$ a pair of distinct elements in $\{1,\ldots,\rankL\}$.
We shall investigate $\epsilon(S_{jk}(1,1),\lw)$.
We note that $\epsilon(S_{jk}(1,r),\lu)\leq \eS+r-1$ for all $r\geq 1$ by the definition of $\eS$ and Lemma \ref{lemma:bound-index-S}.
If $\{j,k\}\cap \{1,2\}=\varnothing$,  then
\begin{align}
\label{eq:epsilon(Sjk(1,1),lw)&leqepsilon(Sjk(1,1),lu)}
\epsilon(S_{jk}(1,1),\lw)&\leq \epsilon(S_{jk}(1,1),\lu)
\end{align}
since $S_{jk}(1,1)_{l}S_{21}(1,r)_{m}=S_{21}(1,r)_{m} S_{jk}(1,1)_{l}$
for all $l,m\in\Z$ and $r=1,2,3$.
For $j\in\{3,\ldots,\rankL\}$, $r\in\{1,2,3\}$, and $l\in\Z$, by \eqref{eq:omega-s11}
\begin{align}
\label{eq:Sj1(1,1)lS21(1,r)eS+r-1lu}
S_{j1}(1,1)_{l}S_{21}(1,r)_{\eS+r-1}\lu&=S_{21}(1,r)_{\eS+r-1}S_{j1}(1,1)_{l}\lu
+r\sum_{s=1}^{r+1}\binom{l}{s}S_{2j}(1,s)_{\eS-2+l+s}\lu.
\end{align}
Thus, 
$\epsilon(S_{j1}(1,1), S_{21}(1,r)_{\eS+r-1}\lu)\leq \max\{\epsilon(S_{j1}(1,1),\lu),1\}$
for $r=1,2,3$ and hence 
\begin{align}
\label{eq:epsilon(Sjk(1,1),lw)&leqepsilon(Sjk(1,1),lu)-2}
\epsilon(S_{j1}(1,1),\lw)\leq \max\{\epsilon(S_{j1}(1,1),\lu),1\}.
\end{align}
The same argument shows that for $j\in\{3,\ldots,\rankL\}$,
\begin{align}
\label{eq:epsilon(Sjk(1,1),lw)&leqepsilon(Sjk(1,1),lu)-3}
\epsilon(S_{j2}(1,1),\lw)\leq \max\{\epsilon(S_{j2}(1,1),\lu),1\}.
\end{align}
We set $(j,k)=(2,1)$. For $i\in\Z$ and $r=1,2,3$, by \cite[Proposition 4.5.7]{LL} putting 
$u=\nS_{21}(1,1), v=\nS_{21}(1,r), p=i, q=\eS+r-1, l=\eS+1$, and $m=0$ in the symbol used there,
we have
\begin{align}
\label{eq:nS21(1,1)inS21(1,r)}
\nS_{21}(1,1)_{i}\nS_{21}(1,r)_{\eS+r-1}\lu&=\sum_{j=0}^{\eS+1}\binom{\eS+1}{j}(\nS_{21}(1,1)_{i-(\eS+1)+j}\nS_{21}(1,r))_{2\eS+r-j}\lu.
\end{align}
By Lemma \ref{lemma:(M(1)GammaPM(1)Gammaprime} (2), $\nS_{21}(1,1)_{i-(\eS+1)+j}\nS_{21}(1,r)$ is an element of $M(1)^{+}_{\varnothing}$.
Since $\epsilon(\omega^{[i]},\lu)\leq 1$ and
	$\epsilon(\Har^{[i]},\lu)\leq 3$ for all $i=1,\ldots,\rankL$,
for any homogeneous element $a\in M(1)^{+}_{\varnothing}$, by using \eqref{eq:a(1)(-l1)ldots a(m)lm},
an inductive argument on $\wt a$ shows that
$a_{i}\lu=0$ for all $i>\wt(a)-1$. 
We see that for $j\in\Z$ and $r\in\Z_{>0}$,
\begin{align}
2\eS+r-j-(\wt((\nS_{21}(1,1)_{i-(\eS+1)+j}\nS_{21}(1,r)))-1)&=\eS+i-2.
\end{align}
Thus, if $\eS\geq 2$, then  by \eqref{eq:nS21(1,1)inS21(1,r)}, $\nS_{21}(1,1)_{i}\nS_{21}(1,r)_{\eS+r-1}\lu=0$ for all $i\geq \eS$
and hence 
\begin{align}
\label{eq:epsilon(S21(1,1),lwleq eS1}
\epsilon(S_{21}(1,1),\lw)&\leq \eS-1.
\end{align}
For any $\lw\in \mW$ and $i=1,\ldots,\rankL$, it follows from Lemma \ref{lemma:epsilon(omegailuleq1epsilonHariluleq3} (1) that
$\omega^{[i]}_{1}\lw, \Har^{[i]}_{3}\lw\in \mW$,
and $\epsilon(\omega^{[i]},\lw)\leq 1$,  $\epsilon(\Har^{[i]},\lw)\leq 3$.
Thus, if ${\eS}\geq 2$, then by 
\eqref{eq:epsilon(Sjk(1,1),lw)&leqepsilon(Sjk(1,1),lu)},
\eqref{eq:epsilon(Sjk(1,1),lw)&leqepsilon(Sjk(1,1),lu)-2},
\eqref{eq:epsilon(Sjk(1,1),lw)&leqepsilon(Sjk(1,1),lu)-3},
and \eqref{eq:epsilon(S21(1,1),lwleq eS1},  we can take a simultaneous eigenvector $\lv$ of 
$\{\omega^{[i]}_{1},\Har^{[i]}_{3}\}_{i=1}^{\rankL}$ in $\mW$
such that $\epsilon(S_{21}(1,1),\lv)<\epsilon(S_{21}(1,1),\lu)$
and $\epsilon(S_{ij}(1,1),\lv)\leq \max\{\epsilon(S_{ij}(1,1),\lu),1\}$
for any pair of distinct element $i,j\in\{1,\ldots,\rankL\}$ with $\{i,j\}\neq \{1,2\}$.
Replacing $\lu$ by this $\lv$ repeatedly, we get a non-zero element $\lu\in\Omega_{M(1)^{+}}(M)$.
Now, the result follows from Lemma 
\ref{lemma:epsilon(omegailuleq1epsilonHariluleq3} (2),
\end{proof}

By Lemma \ref{lemma:Zhu-Omega} and \cite[Theorem 6.2]{DLM1998t}, we have the following result,
which is already shown in \cite[Proposition 3.13]{Tanabe2021-1} when $\rank \lattice=1$:
\begin{proposition}
\label{proposition:M(1)plusmoduleinweak}
Let $\lattice$ be a non-degenerate even lattice of finite rank
and $\module$ a non-zero weak $V_{\lattice}^{+}$-module.
Then, there exists a non-zero $M(1)^{+}$-submodule of $\module$. 
\end{proposition}

\section{Extension groups for $M(1)^{+}$}
\label{section:Extension groups for M(1)+}
In this section we study some weak modules for $M(1)^{+}$ with rank $\rankL$. 
As stated in Section \ref{section:introduction},
the irreducible $M(1)^{+}$-modules 
are classified in \cite[Theorem 4.5]{DN1999-1} for the case of $\dim_{\C}\fh=1$ 
and \cite[Theorem 6.2.2]{DN2001} for the general case (see \eqref{eqn:classificationM(1)p}).
Results in this section will be used in Part $3$ of this series of papers
to show that every irreducible weak $V_{\lattice}^{+}$-module is 
a direct sum of irreducible $M(1)^{+}$-modules.
When $\rankL=1$, some of the results in this section have already 
been obtained in \cite[Section 5]{Abe2005}.
In some parts of the following argument, we shall use techniques in \cite[Section 5]{Abe2005}.
Throughout this section, $\module$ is an $M(1)^{+}$-module,
$W$ is an irreducible $M(1)^{+}$-module, 
and $N$ is a weak $M(1)^{+}$-module.
In this section,  we consider the following exact sequence
\begin{align}
	0\rightarrow W\overset{}{\rightarrow} N\overset{\pi}{\rightarrow} M\rightarrow 0
\label{eq:exact-seq}
	\end{align}
of weak $M(1)^{+}$-modules.
We shall use  the symbols in \eqref{eq:def-oega-i-H-i} and \eqref{eq:def-nS_ij(l,m)}.
 We note that
$[\omega^{[i]}_{1},\omega^{[j]}_{1}]=[\omega^{[i]}_{1},\Har^{[j]}_{3}]=[\Har^{[i]}_{3},\Har^{[j]}_{3}]=0$ for all $i,j=1,\ldots, \rankL$.
Let $B$ be an irreducible $A(M(1)^{+})$-submodule of $\module(0)$.
For $\zeta=(\zeta^{[1]},\ldots,\zeta^{[\rankL]}), \xi=(\xi^{[1]},\ldots,\xi^{[\rankL]})\in \C^{\rankL}$, 
let $\lv\in B$ such that 
\begin{align}
	\label{eq:eigenvalue-omega-H}
(\omega^{[i]}_{1}-\zeta^{[i]})\lv&=(\Har^{[i]}_{3}-\xi^{[i]})\lv=0
\end{align}
for all $i=1,\ldots,\rankL$
and 
we define 
\begin{align}
W_{\zeta,\xi}=\bigcap_{j=1}^{\rankL}\Ker (\omega^{[j]}_{1}-\zeta^{[j]})
\cap
\bigcap_{j=1}^{\rankL}\Ker (\Har^{[j]}_{3}-\xi^{[j]})\cap W.
\end{align}
\begin{lemma}
\label{lemma:preimageofv}
Under the setting above,
there exists $\lu\in N$ such that
\begin{align}
\label{eq:keromega-kerH}
\pi(\lu)=\lv,\quad (\omega^{[i]}_{1}-\zeta^{[i]})\lu, (\Har^{[i]}_{3}-\xi^{[i]})\lu\in W_{\zeta,\xi}
\end{align}
and
\begin{align}
\label{eq:pi(lu)lvmboxand} 
(\omega^{[i]}_{1}-\zeta^{[i]})^2\lu=(\Har^{[i]}_{3}-\xi^{[i]})^2\lu=0
\end{align}
for all $i=1,\ldots, \rankL$.
\end{lemma}
\begin{proof}
Let $\lu\in N$ such that $\pi(\lu)=\lv$.
Since $(\omega^{[i]}_{1}-\zeta^{[i]})\lu, (\Har^{[i]}_{3}-\xi^{[i]})\lu\in W$
and the actions of $\omega^{[i]}_{1}$ and $\Har^{[i]}_{3}$ on $W$ are semisimple for all $i=1,\ldots,\rankL$,
the subspace $U:=\Span_{\C}\{a_{(1)}\cdots a_{(n)}\lu\ |\ n\in\Z_{\geq 0}, a_{(1)},\ldots,a_{(n)}\in \{\omega^{[i]}_{1},\Har^{[i]}_{3}\}_{i=1}^{\rankL}\}$ of $\mN$
is  finite dimensional.
For $\rho=(\rho^{[i]})_{i=1}^{\rankL}, \sigma=(\sigma^{[i]})_{i=1}^{\rankL}\in \C^{\rankL}$, we 
define 
\begin{align}
U_{\rho,\sigma}:=\Big\{\lw\in U\ \Big|
\mbox{
\begin{tabular}{l}
there exists $n\in\Z_{>0}$ such that\\
$(\omega^{[i]}_{1}-\rho^{[i]})^{n}\lw=(\Har^{[i]}_{3}-\sigma^{[i]})^{n}\lw=0$
for all $i=1,\ldots,\rankL$.
\end{tabular}}
\Big\}
\end{align}
and we take a decomposition $U=\oplus_{\rho,\sigma\in\C^{\rankL}}U_{\rho,\sigma}$.
For any $\rho,\sigma\in\C^{\rankL}$, we also take a linear map $f^{\rho,\sigma}\in \Span_{\C}\{a_{(1)}\cdots a_{(n)}\ |\ n\in\Z_{\geq 0}, a_{(1)},\ldots,a_{(n)}\in \{\omega^{[i]}_{1},\Har^{[i]}_{3}\}_{i=1}^{\rankL}\}$  such that $f^{\rho,\sigma}|_{U_{\rho,\sigma}}=\id_{U_{\rho,\sigma}}$
and  $f^{\rho,\sigma}|_{U_{\mu,\nu}}=0$ for all $(\mu,\nu)\neq (\rho,\sigma)$.
We write $\lu=\sum_{\rho,\sigma\in\C^{\rankL}}\lu_{\rho,\sigma}$ where $\lu_{\rho,\sigma}\in U_{\rho,\sigma}$.
We fix $i\in \{1,\ldots,\rankL\}$. Since 
$(\omega^{[i]}_{1}-\zeta^{[i]})\lu\in W$, we have
\begin{align}
\label{eq:(omega[i]1-zeta[i])lurhosigma}
(\omega^{[i]}_{1}-\zeta^{[i]})\lu_{\rho,\sigma}
&=(\omega^{[i]}_{1}-\zeta^{[i]})f^{\rho,\sigma}\lu
=f^{\rho,\sigma}(\omega^{[i]}_{1}-\zeta^{[i]})\lu\in W\cap U_{\rho,\sigma}.
\end{align}
Since the action of $\omega^{[i]}_{1}$ on $W$ is semisimple,
\begin{align}
\label{eq:(omega1]1-rho[i])omega[i]1}
(\omega^{[i]}_1-\rho^{[i]})(\omega^{[i]}_{1}-\zeta^{[i]})\lu_{\rho,\sigma}&=0.
\end{align}
Since $\lu_{\rho,\sigma}\in U_{\rho,\sigma}$, there exists $k\in\Z_{>0}$ such that $(\omega^{[i]}_1-\rho^{[i]})^{k}\lu_{\rho,\sigma}=0$.
When $\rho^{[i]}\neq \zeta^{[i]}$, regarding $(\omega^{[i]}_1-\rho^{[i]})^{k}$ and the left-hand side of \eqref{eq:(omega1]1-rho[i])omega[i]1} 
as polynomials in $\omega^{[i]}_1-\rho^{[i]}$, and dividing the former by the latter, we get
$(\omega^{[i]}_1-\rho^{[i]})\lu_{\rho,\sigma}=0$ and hence
\begin{align}
\lu_{\rho,\sigma}&=\frac{1}{\rho^{[i]}-\zeta^{[i]}}(\omega^{[i]}_1-\zeta^{[i]})\lu_{\rho,\sigma}\in W
\end{align}
by \eqref{eq:(omega[i]1-zeta[i])lurhosigma}.
The same argument shows that
\begin{align}
\label{eq:(omega1]1-rho[i])omega[i]2}
(\Har^{[i]}_3-\sigma^{[i]})(\Har^{[i]}_{3}-\xi^{[i]})\lu_{\rho,\sigma}&=0
\end{align}
and if $\sigma^{[i]}\neq \xi^{[i]}$, then $\lu_{\rho,\sigma}\in W$.
Thus if $(\rho,\sigma)\neq (\zeta,\xi)$, then $\lu_{\rho,\sigma}\in W$
and hence we can take $\lu=\lu_{\zeta,\xi}\in U_{\zeta,\xi}$. In this case,
\eqref{eq:keromega-kerH} 
and 
\eqref{eq:pi(lu)lvmboxand} hold by 
\eqref{eq:(omega1]1-rho[i])omega[i]1}, and \eqref{eq:(omega1]1-rho[i])omega[i]2}.
\end{proof}
Let $\lu\in \mN$ that satisfies \eqref{eq:keromega-kerH} and \eqref{eq:pi(lu)lvmboxand}.
If $(W,B) \not\cong (M(1)^{+},M(1)^{-}(0))$, then \textcolor{black}{it follows from \cite[Lemma 4.8]{Abe2005} that}
\begin{align}
\label{eq:epsilon(omegailu)leq1}
\epsilon(\omega^{[i]},\lu)\leq 1\mbox{ and }\epsilon(\Har^{[i]},\lu)\leq 3 
\end{align}
for all $i=1,\ldots, \rankL$, where $\epsilon=\epsilon_{Y_N}$.

For a pair of distinct elements $i,j\in\{1,\ldots,\rankL\}$, a direct computation shows that
\begin{align}
0&=
6\omega^{[i]}_{-2 } S_{ij}(1,1)
+2\omega^{[j]}_{-2 } S_{ij}(1,1)
\nonumber\\&\quad{}
-4\omega_{0 } \omega^{[i]}_{-1 } S_{ij}(1,1)
+\omega_{0 } \omega_{0 } \omega_{0 } S_{ij}(1,1)
\nonumber\\&\quad{}
+4\omega^{[i]}_{-1 } S_{ij}(1,2)
-4\omega^{[j]}_{-1 } S_{ij}(1,2)
\nonumber\\&\quad{}
-3\omega_{0 } \omega_{0 } S_{ij}(1,2)
+6\omega_{0 } S_{ij}(1,3),\label{eq:s11-3}\\
0&=
32\omega^{[i]}_{-3 } S_{ij}(1,1)
-24\Har^{[i]}_{-1 } S_{ij}(1,1)
\nonumber\\&\quad{}
-8\omega^{[j]}_{-3 } S_{ij}(1,1)
+24\Har^{[j]}_{-1 } S_{ij}(1,1)
\nonumber\\&\quad{}
-120\omega_{0 } \omega^{[i]}_{-2 } S_{ij}(1,1)
+36\omega_{0 } \omega^{[j]}_{-2 } S_{ij}(1,1)
\nonumber\\&\quad{}
+72\omega_{0 } \omega_{0 } \omega^{[i]}_{-1 } S_{ij}(1,1)
-9\omega_{0 } \omega_{0 } \omega_{0 } \omega_{0 } S_{ij}(1,1)
\nonumber\\&\quad{}
+12\omega^{[i]}_{-2 } S_{ij}(1,2)
+12\omega^{[j]}_{-2 } S_{ij}(1,2)
\nonumber\\&\quad{}
-72\omega_{0 } \omega^{[i]}_{-1 } S_{ij}(1,2)
-72\omega_{0 } \omega^{[j]}_{-1 } S_{ij}(1,2)
\nonumber\\&\quad{}
+18\omega_{0 } \omega_{0 } \omega_{0 } S_{ij}(1,2),\label{eq:s11-4-1}\\
0&=
8\omega^{[j]}_{-3 } S_{ij}(1,1)
-24\Har^{[j]}_{-1 } S_{ij}(1,1)
\nonumber\\&\quad{}
+54\omega_{0 } \omega^{[i]}_{-2 } S_{ij}(1,1)
-36\omega_{0 } \omega^{[j]}_{-2 } S_{ij}(1,1)
\nonumber\\&\quad{}
-36\omega_{0 } \omega_{0 } \omega^{[i]}_{-1 } S_{ij}(1,1)
+9\omega_{0 } \omega_{0 } \omega_{0 } \omega_{0 } S_{ij}(1,1)
\nonumber\\&\quad{}
+54\omega^{[i]}_{-2 } S_{ij}(1,2)
-12\omega^{[j]}_{-2 } S_{ij}(1,2)
\nonumber\\&\quad{}
+72\omega_{0 } \omega^{[j]}_{-1 } S_{ij}(1,2)
-18\omega_{0 } \omega_{0 } \omega_{0 } S_{ij}(1,2)
\nonumber\\&\quad{}
+72\omega^{[i]}_{-1 } S_{ij}(1,3),\label{eq:s11-4-2}\\
0&=
14\omega^{[j]}_{-3 } S_{ij}(1,1)
+12\Har^{[j]}_{-1 } S_{ij}(1,1)
-3\omega^{[j]}_{-2 } S_{ij}(1,2)
-36\omega^{[j]}_{-1 } S_{ij}(1,3).\label{eq:s11-4-3}
\end{align}

The following result is a direct consequence of \eqref{eq:def-oega-i-H-i}, Lemma \ref{lemma:(M(1)GammaPM(1)Gammaprime}, and \eqref{eq:a(1)(-l1)ldots a(m)lm}:
\begin{lemma}
\label{lemma:m1plusmodulegenerated}
Let $\mK$ be an $M(1)^{+}$-module such that $\mK=M(1)^{+}\cdot \mK(0)$.
Then, $\mK$ is spanned 
by
$a^{(1)}_{i_1}\cdots a^{(n)}_{i_n}b$ where $n\in\Z_{\geq 0}$,
$b\in \mK(0)$, 
$a^{(j)}\in\{\omega^{[k]}, J^{[k]}\ |\ k=1,\ldots,\rankL\}\cup
\{S_{lm}(1,r)\ |\ 1\leq m<l\leq\rankL, r=1,2,3\}$ 
and $i_j\in\Z_{\leq \wt a^j-2}$
for $j=1,\ldots,n$.
\end{lemma}

\begin{lemma}
Let $U$ be a subspace of a weak $M(1)^{+}$-module which is stable under the actions of	$\{\omega^{[i]}_{1},\Har^{[i]}_{3}\}_{i=1}^{\rankL}$. 
Assume $\epsilon(\omega^{[i]},\lu)\leq 1$ and $\epsilon(\Har^{[i]},\lu)\leq 3$
	for all $\lu\in U$ and $i=1,\ldots,\rankL$.
	Let $i,j\in \{1,\ldots,\rankL\}$ with $i\neq j$
	and ${\eS}\in\Z$ such that
	${\eS}\geq \epsilon(S_{ij},\lu)$ for all $\lu\in U$.
Then, for $\lu\in U$
\begin{align}
\label{eq:(epsilonS-1)big((18zeta[i]+3)epsilonS5-0}
0&=
-{\eS} ({\eS}+1)^2S_{ij}(1,1)_{{\eS} }\lu
-({\eS}+2) (3 {\eS}+1)S_{ij}(1,2)_{{\eS}+1 }\lu
\nonumber\\&\quad{}
+4 {\eS}S_{ij}(1,1)_{{\eS} } \omega^{[i]}_{1 }\lu
-4S_{ij}(1,1)_{{\eS} } \omega^{[j]}_{1 }\lu
\nonumber\\&\quad{}
-2 (3 {\eS}+1)S_{ij}(1,3)_{{\eS}+2 }\lu
+4S_{ij}(1,2)_{{\eS}+1 } \omega^{[i]}_{1 }\lu
-4S_{ij}(1,2)_{{\eS}+1 } \omega^{[j]}_{1 }\lu,\\
\label{eq:(epsilonS-1)big((18zeta[i]+3)epsilonS5-1}
0&=
-\eS(\eS+1)(3\eS^2+27\eS+22)S_{ij}(1,1)_{\eS } \lu\nonumber\\&\quad{}
-2(3\eS^3+39\eS^2+82\eS+24)S_{ij}(1,2)_{\eS+1 } \lu
\nonumber\\&\quad{}
+8\eS(3\eS+11)S_{ij}(1,1)_{\eS } \omega^{[i]}_{1 } \lu
+8(3\eS-13)S_{ij}(1,1)_{\eS } \omega^{[j]}_{1 } \lu
\nonumber\\&\quad{}
-48(3\eS+1)S_{ij}(1,3)_{\eS+2 } \lu
+8(3\eS-13)S_{ij}(1,2)_{\eS+1 } \omega^{[j]}_{1 } \lu
\nonumber\\&\quad{}
-8S_{ij}(1,1)_{\eS } \Har^{[i]}_{3 } \lu
+8S_{ij}(1,1)_{\eS } \Har^{[j]}_{3 } \lu
\nonumber\\&\quad{}
+8(3\eS+11)S_{ij}(1,2)_{\eS+1 } \omega^{[i]}_{1 } \lu,\\
\label{eq:(epsilonS-1)big((18zeta[i]+3)epsilonS5-2}
0&=
2(3\eS^3+21\eS^2+42\eS+14)S_{ij}(1,2)_{\eS+1 } \lu
-8(3\eS-7)S_{ij}(1,1)_{\eS } \omega^{[j]}_{1 } \lu
\nonumber\\&\quad{}
+4(18\eS+7)S_{ij}(1,3)_{\eS+2 } \lu
-8(3\eS-7)S_{ij}(1,2)_{\eS+1 } \omega^{[j]}_{1 } \lu
\nonumber\\&\quad{}
-8S_{ij}(1,1)_{\eS } \Har^{[j]}_{3 } \lu
+3\eS(\eS+1)^2(\eS+4)S_{ij}(1,1)_{\eS } \lu
\nonumber\\&\quad{}
-12\eS(\eS+4)S_{ij}(1,1)_{\eS } \omega^{[i]}_{1 } \lu
-36S_{ij}(1,2)_{\eS+1 } \omega^{[i]}_{1 } \lu
\nonumber\\&\quad{}
+24S_{ij}(1,3)_{\eS+2 } \omega^{[i]}_{1 } \lu,\\
\label{eq:(epsilonS-1)big((18zeta[i]+3)epsilonS5-3}
0&=
-S_{ij}(1,2)_{{\eS}+1 }\lu
-5S_{ij}(1,1)_{{\eS} }\omega^{[j]}_{1 } \lu
-S_{ij}(1,3)_{{\eS}+2 }\lu
-11S_{ij}(1,2)_{{\eS}+1 }\omega^{[j]}_{1 } \lu
\nonumber\\&\quad{}
+2S_{ij}(1,1)_{{\eS} }\Har^{[j]}_{3 } \lu
-6S_{ij}(1,3)_{{\eS}+2 }\omega^{[j]}_{1 } \lu.
\end{align}
If $\lu$ is a simultaneous eigenvector of $\{ \omega^{[i]}_{1},\omega^{[j]}_{1},\Har^{[i]}_{3},\Har^{[j]}_{3}\}$
with eigenvalues $\{\zeta^{[i]},\zeta^{[j]},\xi^{[i]},\xi^{[j]}\}$:
\begin{align}
	(\omega^{[i]}_{1}-\zeta^{[i]})\lu&=
	(\omega^{[j]}_{1}-\zeta^{[j]})\lu
=(\Har^{[i]}_{3}-\xi^{[i]})\lu	=(\Har^{[j]}_{3}-\xi^{[j]})\lu=0,
\label{eq:omega-zeta-H4-0}
\end{align}
then
\begin{align}
0&=-(\eS-1)\big((18\zeta^{[i]}+3)\eS^5+(-54\zeta^{[i]}+6)\eS^4\nonumber\\
&\qquad{}+(1-36 \zeta^{[j]}-78 \zeta^{[i]}-216 \zeta^{[j]} \zeta^{[i]}+216 (\zeta^{[i]})^2)\eS^3\nonumber\\
&\qquad{}+(-2 + 4\zeta^{[j]} + 22 \zeta^{[i]} + 744 \zeta^{[j]} \zeta^{[i]} + 24 (\zeta^{[i]})^2)\eS^2\nonumber\\
&\qquad{}+(12 \zeta^{[i]} - 192 \zeta^{[j]} \zeta^{[i]} - 48 (\zeta^{[i]})^2 - 1152 \zeta^{[j]} (\zeta^{[i]})^2)\eS\nonumber\\
&\qquad{}+384\zeta^{[j]}(\zeta^{[i]})^2-16\zeta^{[j]}\zeta^{[i]}\big)\nonumber\\
&\quad{}+8\big((9\eS^4+12\eS^3+(-18\zeta^{[i]}-36\zeta^{[j]})\eS^2+(-24\zeta^{[i]}-1)\eS\nonumber\\
&\qquad{}-4 \zeta^{[j]} - 6 \zeta^{[i]} - 24 \zeta^{[j]} \zeta^{[i]} + 24 (\zeta^{[i]})^2\big)\xi^{[i]}\nonumber\\
 	&\quad-8\big(((18\zeta^{[i]}+3)\eS^2+
(-24\zeta^{[i]}+1)\eS+24(\zeta^{[i]})^2+(-24\zeta^{[j]}-6)\zeta^{[i]}-4\zeta^{[j]})\big)\xi^{[j]},\label{eq:s11-zeta-1}\\
0&=-(\eS-1)\big((18\zeta^{[j]}+3)\eS^5+(-54\zeta^{[j]}+6)\eS^4\nonumber\\
&\qquad{}+(1-36 \zeta^{[i]}-78 \zeta^{[j]}-216 \zeta^{[i]} \zeta^{[j]}+216 (\zeta^{[j]})^2)\eS^3\nonumber\\
&\qquad{}+(-2 + 4\zeta^{[i]} + 22 \zeta^{[j]} + 744 \zeta^{[i]} \zeta^{[j]} + 24 (\zeta^{[j]})^2)\eS^2\nonumber\\
&\qquad{}+(12 \zeta^{[j]} - 192 \zeta^{[i]} \zeta^{[j]} - 48 (\zeta^{[j]})^2 - 1152 \zeta^{[i]} (\zeta^{[j]})^2)\eS\nonumber\\
&\qquad{}+384\zeta^{[i]}(\zeta^{[j]})^2-16\zeta^{[i]}\zeta^{[j]}\big)\nonumber\\
&\quad{}+8\big((9\eS^4+12\eS^3+(-18\zeta^{[j]}-36\zeta^{[i]})\eS^2+(-24\zeta^{[j]}-1)\eS\nonumber\\
&\qquad{}-4 \zeta^{[i]} - 6 \zeta^{[j]} - 24 \zeta^{[i]} \zeta^{[j]} + 24 (\zeta^{[j]})^2\big)\xi^{[j]}\nonumber\\
 	&\quad-8\big(((18\zeta^{[j]}+3)\eS^2+
(-24\zeta^{[j]}+1)\eS+24(\zeta^{[j]})^2+(-24\zeta^{[i]}-6)\zeta^{[j]}-4\zeta^{[i]})\big)\xi^{[i]}.
\label{eq:s11-zeta-2}
\end{align}
\end{lemma}
\begin{proof}
We first note that \textcolor{black}{interchanging} the positions of  $\zeta^{[i]}$ and $\zeta^{[j]}$, and $\xi^{[i]}$ and $\xi^{[j]}$ in \eqref{eq:s11-zeta-1},
we get \eqref{eq:s11-zeta-2}. Thus, these two equations \eqref{eq:s11-zeta-1} and \eqref{eq:s11-zeta-2}
are essentially the same, however,
we put them here  because they are convenient for later use.
	By Lemma \ref{lemma:bound-index-S}, $\epsilon(S_{ij}(1,r),\lu)\leq {\eS}+r-1$ for all
	$r=1,2,\ldots$ and $\lu\in U$.
We shall apply Lemma \ref{lemma:ksumklimitsm=1n(kepsilon(a(m))+im)kk+(kepsilon(b)+j)-1+kr} (2) to \eqref{eq:s11-3}
with $A:=\{\omega^{[i]},\omega^{[j]},\Har^{[i]},\Har^{[j]}\}$, $B:=\{S_{ij}(1,r)\ |\ r=1,2,3\}$,
$\epsilon(\omega^{[k]}):=\wt(\omega^{[k]})-1=1$, $\epsilon(\Har^{[k]}):=\wt(\Har^{[k]})-1=3$ for $k=i,j$, and 
$\epsilon(S_{ij}(1,r))=\epsilon(S)+r-1$ for $r=1,2,3$.
By Lemma \ref{lemma:(M(1)GammaPM(1)Gammaprime}, we have
$
A\cdot \langle A_{-}\rangle\vac \subset \langle A_{-}\rangle\vac,
A\cdot (\langle A_{-}\rangle B_{-}\vac)\subset \langle A_{-}\rangle B_{-}\vac,
B\cdot \langle A_{-}\rangle\vac \subset B_{-}\langle A_{-}\rangle\vac$,
and 
$A\cdot (B_{-}\langle A_{-}\rangle\vac) \subset B_{-}\langle A_{-}\rangle\vac$,
 where the symbols $\langle A_{-}\rangle\vac, \langle A_{-}\rangle B_{-}\vac$, and $ B_{-}\langle A_{-}\rangle\vac$
are defined in \eqref{eqn:A-B:=Span} and \eqref{eq:a(1)-i1cdotsa(n)-inb}.
The weight of each term in \eqref{eq:s11-3} is $5$. The same argument as in \eqref{eq:delta(wt(a)+wt(S211,j))-1)-epsilon(S21(1,j))}
shows $\delta(5)=\epsilon(S)+3$, where $\delta$ is defined in \eqref{eq:gamma(k):=maxBiggsum}.
By Lemma \ref{lemma:ksumklimitsm=1n(kepsilon(a(m))+im)kk+(kepsilon(b)+j)-1+kr} (2), the $({\eS}+3)$-th action of \eqref{eq:s11-3} on $\lu$ is a linear combination of  elements of the form
\begin{align}
q_{\epsilon(q)}p^{(1)}_{\epsilon(p^{(1)})}\cdots p^{(m)}_{\epsilon(p^{(m)})}\lu
\end{align}
where $m\in\Z_{\geq 0}$, $p^{(1)},\ldots,p^{(m)}\in A$, and $q\in B$.
To obtain the explicit expression of the result \eqref{eq:(epsilonS-1)big((18zeta[i]+3)epsilonS5-0} we use computer algebra system Risa/Asir\cite{Risa/Asir}.
By	taking the $({\eS}+4)$-th actions of \eqref{eq:s11-4-1}--\eqref{eq:s11-4-3} on $\lu$,
the same argument shows \eqref{eq:(epsilonS-1)big((18zeta[i]+3)epsilonS5-1}--\eqref{eq:(epsilonS-1)big((18zeta[i]+3)epsilonS5-3}.
Deleting the terms including $S_{ij}(1,3)_{{\eS}+2}\lu$ and $S_{ij}(1,2)_{{\eS}+2}\lu$
from \eqref{eq:(epsilonS-1)big((18zeta[i]+3)epsilonS5-0}--\eqref{eq:(epsilonS-1)big((18zeta[i]+3)epsilonS5-3},  we have 
	\eqref{eq:s11-zeta-1} and \eqref{eq:s11-zeta-2}.
		\end{proof}

\begin{lemma}
\label{lemma:M1lambda-submodule}
Let 
\begin{align}
	0\rightarrow \mW\overset{}{\rightarrow} N\overset{\pi}{\rightarrow} \module\rightarrow 0
\label{eq:exact-seq-lambda}
	\end{align}
be an exact sequence of weak $M(1)^{+}$-modules
where $\mW$ is an irreducible $M(1)^{+}$-module,
$N$ is a weak $M(1)^{+}$-module and $\module=\oplus_{i\in \gamma+\Z_{\geq 0}}\module_{i}$ is an $M(1)^{+}$-module.
Let $B$ be an irreducible $A(M(1)^{+})$-submodule of $\module_{\gamma}$
that is not isomorphic to $\mW(0)$ 
and $\lv$
 a simultaneous eigenvector in $B$ of $\{\omega^{[i]}_{1},\Har^{[i]}_{3}\}_{i=1}^{\rankL}$
with eigenvalues $\{\zeta^{[i]},\xi^{[i]}\}_{i=1}^{\rankL}$:
\begin{align}
(\omega^{[i]}_1-\zeta^{[i]})\lv=(\Har^{[i]}_3-\xi^{[i]})\lv=0.
\end{align}
If $(W,B)\not\cong (M(1)^{+},M(1)^{-}(0))$,
then there exists a preimage $\lu\in N_{\gamma}$ of $\lv$ under the canonical projection $\mN_{\gamma}\rightarrow \module_{\gamma}$
such that
\begin{align}
(\omega^{[i]}_1-\zeta^{[i]})\lu=(\Har^{[i]}_3-\xi^{[i]})\lu=0
\end{align}
for all $i=1,\ldots,\rankL$.
\end{lemma}
\begin{proof}
Using \cite[Proposition 4.3]{Abe2005} and eigenvalues of $\omega^{[i]}_{1}$ and $\Har^{[i]}_{3}$  for $i=1,\ldots,\rankL$ on 
irreducible $M(1)^{+}$-modules in \cite[Table 1]{AD2004},
we see that the result holds
if $\mW=M(1)(\theta)^{\pm}$ or $B=M(1)(\theta)^{\pm}(0)$.
We discuss the other cases.
For $v\in B$ with \eqref{eq:eigenvalue-omega-H}
, we take $\lu\in\mN$ that satisfies \eqref{eq:keromega-kerH} and \eqref{eq:pi(lu)lvmboxand}.

Let $B=\C e^{\lambda}$ for some $\lambda\in\fh\setminus\{0\}$.
In this case $\zeta^{[i]}=\langle \lambda,h^{[i]}\rangle^2/2$ and
$\xi^{[i]}=0$ for $i=1,\ldots,\rankL$.
We note that at least one of $\zeta^{[1]},\ldots,\zeta^{[\rankL]}$ is not zero.
Let $W=M(1,\mu)$.  Since $B\not\cong W(0)$ as $A(M(1)^{+})$-modules, we have $\mu\in\fh\setminus\{0,\pm\lambda\}$.
Since $\cap_{j=1}^{\rankL}\Ker \Har^{[j]}_{3}\cap M(1,\mu)=\C e^{\mu}$ by \cite[Proposition 4.3]{Abe2005}, 
		\begin{align}
			M(1,\mu)_{\zeta,(0,\ldots,0)}&=\bigcap_{j=1}^{\rankL}\Ker (\omega^{[j]}_{1}-\zeta^{[j]})\cap\C e^{\mu}.
			\label{eq:cap-omega-mu}
		\end{align}
Assume  
\begin{align}
\label{eq:omega-i-zeta-or-H-0}
(\omega^{[i]}_{1}-\zeta^{[i]})\lu\neq 0\mbox{ or }\Har^{[i]}_{3}\lu\neq 0
\mbox{ for some }i\in \{1,\ldots,\rankL\}.
\end{align}
It follows from  \eqref{eq:keromega-kerH} that $M(1,\mu)_{\zeta,(0,\ldots,0)}\neq 0$ 
		and hence 
		$\langle\lambda,h^{[j]}\rangle=\pm\langle\mu,h^{[j]}\rangle$ for all $j=1,\ldots,\rankL$ by \eqref{eq:cap-omega-mu}.
Thus, $\gamma=\langle\lambda,\lambda\rangle/2=\langle\mu,\mu\rangle/2$.
By this and $\lambda\neq\pm\mu$, 
we see that there exists an $A(M(1)^{+})$-submodule of $N(0)$ which is isomorphic to 
$B\oplus M(1,\mu)(0)\cong M(1,\lambda)(0)\oplus M(1,\mu)(0)$. Thus, we have the result.
If $W=M(1)^{\pm}$, then the result follows from the fact that 
$M(1)^{\pm}_{\zeta,(0,\ldots,0)}=0$.

If $B=\C \vac=M(1)^{+}(0)$, then the same argument as above shows the result.

Let $B=M(1)^{-}(0)$,
$W=M(1,\lambda)$ such that $\lambda\in\fh\setminus\{0\}$, and $\lv=h^{[j]}(-1)\vac$ for some $j\in \{1,\ldots,\rankL\}$.
Since $\xi^{[i]}=\delta_{ij}$ for all $i=1,\ldots,\rankL$, it follows from \cite[Proposition 4.3]{Abe2005} that
\begin{align}
M(1,\lambda)_{\zeta,\xi}&\subset \C h^{[j]}(-1)e^{\lambda}.
\end{align}
Suppose there exists $i\in\{1,\ldots,\rankL\}$ such that
$(\omega^{[i]}-\delta_{ij})\lu\neq 0$ or $(\Har^{[i]}-\delta_{ij})\lu\neq 0$.
Then, $M(1,\lambda)_{\zeta,\xi}\neq 0$ and hence $\delta_{jk}=\langle \lambda, h^{[k]}\rangle^2/2+\delta_{jk}$
for all $k=1,\ldots,\rankL$, which contradicts that $\lambda\neq 0$.
The proof is complete.
\end{proof}

We will prepare the following symbol for Lemmas \ref{lemma:M1lambda-submodule-2} and \ref{lemma:Ext-M-M}:
\begin{definition}
Let $R[x]$ be a polynomial ring over a commutative ring $R$.
For two polynomials $A_1=\sum_{i=0}^{\deg A_1}A_{1,i}x^i,A_2=\sum_{i=0}^{\deg A_2}A_{2,i}x^i\in R[x]$ with $A_{ki}\in R$,
we define a polynomial $G(A_1,A_2)\in R[x]$ as follows.
We first prepare indeterminates $\hat{A}_{1,0},\ldots,\hat{A}_{1,\deg A_1}, \hat{A}_{2,0},\ldots,\hat{A}_{2,\deg A_2}$
over $\C$.
We define $\hat{R}:=\C[\hat{A}_{1,0},\ldots,\hat{A}_{1,\deg A_1}, \hat{A}_{2,0},\ldots,\hat{A}_{2,\deg A_2}]$
and two polynomials $\hat{A_1}:=
\sum_{i=0}^{\deg A_1}\hat{A}_{1,i}x^i, \hat{A}_2:=\sum_{i=0}^{\deg A_2}\hat{A}_{2,i}x^i\in 
\hat{R}[x]$. 
In the following,  $\deg \hat{P}$ is the degree of $\hat{P}\in \hat{R}[x]$ with respect to $x$. 
If $\deg \hat{A}_1\geq \deg \hat{A}_2$, then we define $\hat{A}_3\in \hat{R}[x]$ by the remainder  after dividing $\hat{A}_{2,\deg \hat{A}_2}^{\deg \hat{A}_1-\deg \hat{A}_2+1}\hat{A}_1$ by $\hat{A}_2$:
 \begin{align}
\hat{A}_{2,\deg A_2}^{\deg \hat{A}_1-\deg \hat{A_2}+1}\hat{A}_1&=\hat{B}_{2}\hat{A}_2+\hat{A}_3,\quad \hat{B}_{2},\hat{A}_{3}\in \hat{R}[x], \deg  \hat{A}_3<\deg  \hat{A}_2.
\end{align}
Here we note that such an $\hat{A}_3$ exists uniquely since we took $\hat{A}_{2,\deg \hat{A}_2}^{\deg \hat{A}_1-\deg \hat{A}_2+1}\hat{A}_1$ instead of $\hat{A}_1$.
If $\deg \hat{A}_2>\deg\hat{A}_1$, then we define $\hat{A}_3$ by $\hat{A}_1$.
Defining
$A_3:=\hat{A}_{3}|_{\hat{A}_{1,0}=A_{1,0},\hat{A}_{1,1}=A_{1,1},\ldots,\hat{A}_{2,0}=A_{2,0},\ldots}$,
$B_2:=\hat{B}_{2}|_{\hat{A}_{1,0}=A_{1,0},\hat{A}_{1,1}=A_{1,1},\ldots,\hat{A}_{2,0}=A_{2,0},\ldots}\in R[x]$,
we have 
 \begin{align}
{A}_{2,\deg A_2}^{\deg {A}_1-\deg {A_2}+1}{A}_1&={B}_{2}{A}_2+{A}_3,\quad {B}_{2},{A}_{3}\in {R}[x], \deg  {A}_3<\deg  {A}_2.
\end{align}
We replace $A_1$ by $A_2$ and $A_2$ by $A_3$ and repeat this operation as many times as possible, which is the essentially Euclidean algorithm:
\begin{align*}
\begin{array}{c}
\vdots\\
{A}_{k+1,\deg  {A}_{k+1}}^{\deg  {A}_k-\deg  {A}_{k+1}+1}{A}_k={B}_{k+1}{A}_{k+1}+{A}_{k+2},\quad {B}_{k+1}, {A}_{k+2}\in {R}[x], \deg  {A}_{k+2}<\deg  {A}_{k+1}\ (k=0,1,\ldots)\\
\vdots\\
{A}_{d+1,\deg  {A}_{d+1}}^{\deg  {A}_d-\deg  {A}_{d+1}+1}{A}_d={B}_{d+1}{A}_{d+1},\quad {B}_{d+1}\in {R}[x]. 
\end{array}
\end{align*}
Then, we define 
\begin{align}
\label{eq:G(A1,A2):=Ad+1inR[x]}
G(A_1,A_2)&:={A}_{d+1}\in R[x].
\end{align}
\end{definition}
\begin{lemma}
\label{lemma:M1lambda-submodule-2}
Let 
\begin{align}
	0\rightarrow \mW\overset{}{\rightarrow} N\overset{\pi}{\rightarrow} \module\rightarrow 0
\label{eq:exact-seq-lambda-2}
	\end{align}
be an exact sequence of weak $M(1)^{+}$-modules
where $\mW$ is an irreducible $M(1)^{+}$-module,
$N$ is a weak $M(1)^{+}$-module and $\module=\oplus_{i\in \gamma+\Z_{\geq 0}}\module_{i}$ is an $M(1)^{+}$-module.
Let $B$ be an irreducible $A(M(1)^{+})$-submodule of $\module_{\gamma}$
and $\lv$
a simultaneous eigenvector in $B$ for $\{\omega^{[i]}_{1},\Har^{[i]}_{3}\}_{i=1}^{\rankL}$ 
with eigenvalues $\{\zeta^{[i]},\xi^{[i]}\}_{i=1}^{\rankL}$:
\begin{align}
(\omega^{[i]}_1-\zeta^{[i]})\lv=(\Har^{[i]}_3-\xi^{[i]})\lv=0.
\end{align}
Let $\lw\in \mN_{\gamma}$
 such that 
$(\omega^{[i]}_1-\zeta^{[i]})\lw=(\Har^{[i]}_3-\xi^{[i]})\lw=0$
 for all $i=1,\ldots,\rankL$.
If $(W,B)\not\cong (M(1)^{+}, M(1)^{-}(0))$, then $\lw \in \Omega_{M(1)^{+}}(N_{\gamma})$.
\end{lemma}
\begin{proof}
Assume $(W,B)\not\cong (M(1)^{+}, M(1)^{-}(0))$.
By Lemma \ref{lemma:bound-index-S}
and \eqref{eq:epsilon(omegailu)leq1}, 
it is enough to show that $\epsilon(S_{ij}(1,1),\lw)\leq 1$
for any pair of distinct elements $i,j\in\{1,\ldots, \rankL\}$.
For such a pair $i,j$, we write $\epsilon(S_{ij})=\epsilon(S_{ij}(1,1),\lw)$	for simplicity.

\begin{enumerate}
\item
Let $B\cong M(1,\lambda)(0)$ for some $\lambda\in\fh\setminus\{0\}$. In this case $\xi^{[i]}=0$ for all $i=1,\ldots,\rankL$.
		Assume $\langle\lambda,\lambda\rangle\neq 0$.
		Then, we may assume $\lambda\in \C h^{[1]}$ and hence
		$\langle\lambda,h^{[i]}\rangle=\zeta^{[i]}=0$ for all $i=2,\ldots,\rankL$.
For $i=2,\ldots,\rankL$, 
substituting $\zeta^{[i]}=0$ and $\xi^{[1]}=\xi^{[i]}=0$
into \eqref{eq:s11-zeta-1} and \eqref{eq:s11-zeta-2} with $j=1$,  we have
\begin{align}
0&=\epsilon(S_{i1})^2 (\epsilon(S_{i1})-1) (4(-9\epsilon(S_{i1})+1)\zeta^{[1]}\textcolor{black}{+}(\epsilon(S_{i1})+1)(3 \epsilon(S_{i1})^2+3 \epsilon(S_{i1})-2) )\quad\mbox{and} \label{eq:s11-k-11}\\
0&=\epsilon(S_{i1}) (\epsilon(S_{i1})-1) 
\Big((216\epsilon(S_{i1})^2+24\epsilon(S_{i1})-48)(\zeta^{[1]})^2\nonumber\\
&\qquad{}+(18\epsilon(S_{i1})^4-54\epsilon(S_{i1})^3-78\epsilon(S_{i1})^2+22\epsilon(S_{i1})+12)\zeta^{[1]}\nonumber\\
&\qquad{}+3\epsilon(S_{i1})^4+6\epsilon(S_{i1})^3+\epsilon(S_{i1})^2-2\epsilon(S_{i1})\Big).
 \label{eq:s11-k-2}
\end{align}
If we take $A_1$ to the right-hand side of \eqref{eq:s11-k-2} and $A_2$ to the right-hand side of \eqref{eq:s11-k-11}
and if we regard $A_1$ and $A_2$ as polynomials in the variable $\zeta^{[1]}$, 
then $G(A_1,A_2)$ in \eqref{eq:G(A1,A2):=Ad+1inR[x]} is a non-zero scalar multiple of 
\begin{align}
\label{eqn:epsilon(Si1)5epsilon(Si1-1)4}
&\epsilon(S_{i1})^5\epsilon(S_{i1}-1)^4(\epsilon(S_{i1}\textcolor{black}{)}+1)(2\epsilon(S_{i1})+1)\nonumber\\
&\quad{}\times(3\epsilon(S_{i1})-2)(3\epsilon(S_{i1})+1)^2
(3\epsilon(S_{i1})^2+3\epsilon(S_{i1})-2).
\end{align}
Since $G(A_1,A_2)=0$, we have $\epsilon(S_{i1})\leq 1$.
	
Assume $\langle\lambda,\lambda\rangle=0$. Then, we may assume $0\neq \langle \lambda,h^{[1]}\rangle^2
		=-\langle \lambda,h^{[2]}\rangle^2$ 
		and $\langle\lambda,h^{[j]}\rangle=0$ for all $j=3,4,\ldots,\rankL$. 
By substituting $\zeta^{[2]}=-\zeta^{[1]}$ into \eqref{eq:s11-zeta-1} and \eqref{eq:s11-zeta-2},
the same argument as above shows
that 
$\epsilon(S_{21})\leq 1$.

In both the cases of $\langle\lambda,\lambda\rangle\neq0$ and $\langle\lambda,\lambda\rangle=0$,
for the other $i,j$, since one of $\zeta^{[i]}$ or $\zeta^{[j]}$ is $0$,
		the same argument as above also shows that $\epsilon(S_{ij})\leq 1$.
\item
Let $B\cong M(1)^{-}(0)$.  Assume $\mW\not\cong M(1)^{+}$.
If $\rankL=1$, then the result is shown in \textcolor{black}{\cite[Theorem 5.5]{Abe2005}}.
Assume $\rankL\geq 2$.
Let $i,j$ be a pair of distinct elements in $\{1,\ldots, \rankL\}$.
If $(\zeta^{[i]},\xi^{[i]})=(\zeta^{[j]},\xi^{[j]})=(0,0)$,
then it follows from  \eqref{eq:epsilon(omegailu)leq1} and \eqref{eq:s11-zeta-1} that
\begin{align}
	0&=\epsilon(S_{ij})^2(\epsilon(S_{ij})-1)(\epsilon(S_{ij})+1)(3\epsilon(S_{ij})^2+3\epsilon(S_{ij})-2)
	\end{align}
and hence $\epsilon(S_{ij})\leq 1$.

Assume $(\zeta^{[i]},\xi^{[i]})=(0,0)$ and $(\zeta^{[j]},\xi^{[j]})=(1,1)$. It follows from 
\eqref{eq:epsilon(omegailu)leq1} 
and
\eqref{eq:s11-zeta-1} that
\begin{align}
\label{eq:epsilonS13epsilonS615}
0=(\epsilon(S_{ij})-2) (\epsilon(S_{ij})-1) (3 \epsilon(S_{ij})^4+12 \epsilon(S_{ij})^3-11 \epsilon(S_{ij})^2-20 \epsilon(S_{ij})-16)
\end{align}
and hence $\epsilon(S_{ij})=1$ or $2$.
We further assume that $\epsilon(S_{ij})=2$. 
By \eqref{eq:(epsilonS-1)big((18zeta[i]+3)epsilonS5-0}--\eqref{eq:(epsilonS-1)big((18zeta[i]+3)epsilonS5-3},
\begin{align}
\label{eq:Sij(1,2)3textcolorredlw}
S_{ij}(1,2)_{3}\textcolor{black}{\lw}&=-2S_{ij}(1,1)_{2}\lw
\mbox{ and }S_{ij}(1,3)_{4}\lw=3S_{ij}(1,1)_{2}\lw.
\end{align}
We note that \textcolor{black}{$S_{ij}(1,1)_{2}\lw\in W$.}
Using commutation relations (see Remark \ref{remark:com-appendix}) and \eqref{eq:Sij(1,2)3textcolorredlw}, we have 
\begin{align}
\label{eq:}
\omega^{[k]}_{1}S_{ij}(1,1)_{2}\lw&=\Har^{[k]}_{3}S_{ij}(1,1)_{2}\lw=0
\end{align}
for all $k=1,\ldots,\rankL$.
It follows from \cite[Proposition 4.3]{Abe2005} that
there is no non-zero element $\lv$ in any irreducible $M(1)^{+}$-module except $\vac\in M(1)^{+}$  that satisfies 
$\omega^{[k]}_{1}\lv=\Har^{[k]}_{3}\lv=0$ for all $k=1,\ldots,\rankL$.
This is a contradiction.

The same argument as above shows the results for the case that $B\cong M(1)^{+}(0)$ or $M(1)(\theta)^{\pm}(0)$.
\end{enumerate}
\end{proof}
	\begin{lemma}
\label{lemma:Ext-split-lambda-mu}
For any pair of non-isomorphic irreducible $M(1)^{+}$-modules $\module,\mW$
such that $(M,W)\not\cong (M(1)^{+},M(1)^{-})$ and $(M(1)^{-},M(1)^{+})$,
$\Ext^{1}_{M(1)^{+}}(\module,W)=0$.
\end{lemma}
	\begin{proof}
Let $\mN$ be a weak $M(1)^{+}$-module and 
\begin{align}
\label{eq:rightarrowM(1mu)oversetiotarightarrowN}
	0\rightarrow \mW\overset{}{\rightarrow} N\overset{\pi}{\rightarrow} \module\rightarrow 0
	\end{align}
an exact sequence of weak $M(1)^{+}$-modules.
We write $\module =\oplus_{i\in\gamma+\Z_{\geq 0}}\module_{i}$ with $\module_{\gamma}\neq 0$
and $W=\oplus_{i\in \delta+\Z_{\geq 0}}\mW_{i}$ with $\mW_{\delta}\neq 0$.
By Lemmas \ref{lemma:M1lambda-submodule} and \ref{lemma:M1lambda-submodule-2}, there exists $\lu\in \Omega_{M(1)^{+}}(\mN)$ such that
$0\neq \pi(\lu)\in \module_{\gamma}$.

Assume $\mW\cap (M(1)^{+}\cdot \lu)\neq 0$.
Since $\mW_{\delta}\subset \mW\cap (M(1)^{+}\cdot \lu)$,
$\delta\in \gamma+\Z_{\geq 0}$.
Since $\mW\cap (M(1)^{+}\cdot \lu)\neq 0$, we have $\Ext^{1}_{M(1)^{+}}(\module, \mW)\neq 0$
and hence
$\Ext^{1}_{M(1)^{+}}(\mW,\module)\neq 0$ by \cite[Proposition 2.5]{Abe2005}
and \cite[Proposition 3.5]{ADL2005}.
Thus,  there exists a 
non-split exact sequence
$0\rightarrow \module\overset{}{\rightarrow} N\overset{}{\rightarrow} \mW\rightarrow 0
$ of weak $M(1)^{+}$-modules.
The same argument as above shows that $\gamma\in \delta+\Z_{\geq 0}$
and hence $\gamma=\delta$.
Since $M\not\cong W$, $N(0)=N_{\gamma}\cong \module_{\gamma}\oplus \mW_{\delta}$ as $A(M(1)^{+})$-modules.
Thus, the sequence \eqref{eq:rightarrowM(1mu)oversetiotarightarrowN} splits, a contradiction.

\end{proof}
	
\begin{lemma}
\label{lemma:Ext-M-M}
For $M=M(1)^{+},M(1)^{-},M(1)(\theta)^{+}$, and $M(1)(\theta)^{-}$,
$\Ext^{1}_{M(1)^{+}}(M,M)=0$.
\end{lemma}
\begin{proof}
	Let $\mW$ be an $M(1)^{+}$-module such that $\mW\cong \module$,
	$\mN$ a weak $M(1)^{+}$-module, and 
	\begin{align}
		\label{eq:rightarrowM(1mu)oversetiotarightarrowN-1}
		0\rightarrow \mW\overset{}{\rightarrow} N\overset{\pi}{\rightarrow} \module\rightarrow 0
	\end{align}
	an exact sequence of weak $M(1)^{+}$-modules.
	We take $\lv\in\module(0)$ and $\lu\in \mN$ as in \eqref{eq:eigenvalue-omega-H}, \eqref{eq:keromega-kerH},
and	\eqref{eq:pi(lu)lvmboxand}.	
In the case of $M=M(1)^{+}$, the same argument as in the proof of 
\cite[Proposition 5.1]{Abe2005} shows that $\Ext^{1}_{M(1)^{+}}(M(1)^{+},M(1)^{+})=0$.

For $M=M(1)^{-}$ or $M(1)(\theta)^{\pm}$, it is enough to show that
$N(0)\cong \mW(0)\oplus M(0)\cong M(0)\oplus M(0)$ as $A(M(1)^{+})$-modules.
In the Zhu algebra $A(M(1)^{+})$, we have
\begin{align}
\label{eq:omega-H-commute}
\omega^{[i]}*\Har^{[i]}\equiv \Har^{[i]}*\omega^{[i]}
\end{align}
and recall that the following relations from \cite[(6.1.11) and (6.1.10)]{DN2001}:
\begin{align}
\label{eq:omega-1-omega-1/16-omega-9/16-2}
(\omega^{[i]}-\vac)*(\omega^{[i]}-\dfrac{1}{16}\vac)*(\omega^{[i]}-\dfrac{9}{16}\vac)*\Har^{[i]}&\equiv 0,\\
\label{eq:omega-1-omega-1/16-omega-9/16-1}
(132(\omega^{[i]})^2-65\omega^{[i]}-70\Har^{[i]}+3)*\Har^{[i]}&\equiv 0
\end{align}  
 for $i=1,\ldots, \rankL$. Here, we note that $\Har_{a}$ in \cite[Section 6]{DN2001}
is equal to the image of $-9\Har^{[a]}$ under the projection $M(1)^{+}\rightarrow A(M(1)^{+})$
for $a=1,\ldots,\rankL$. Let $A_1$ be the quotient of the right-hand side of \eqref{eq:omega-1-omega-1/16-omega-9/16-2} by $\Har^{[i]}$ 
and $A_2$ the quotient of the right-hand side of \eqref{eq:omega-1-omega-1/16-omega-9/16-1} by 
$\Har^{[i]}$: $A_1:=(\omega^{[i]}-\vac)*(\omega^{[i]}-(1/16)\vac)*(\omega^{[i]}-(9/16)\vac)$
and $A_2:=132(\omega^{[i]})^2-65\omega^{[i]}-70\Har^{[i]}+3$.
If we regard $A_1$ and $A_2$ as polynomials in $\omega^{[i]}$, then $G(A_1,A_2)$ in \eqref{eq:G(A1,A2):=Ad+1inR[x]} is a non-zero scalar multiple of 
\begin{align}
(\Har^{[i]}-1)*(\Har^{[i]}-\dfrac{-1}{128})*(\Har^{[i]}-\frac{15}{128})
\end{align}
and hence
\begin{align}
\label{eq:\Har[i]*(Har[i]-1)*(Har[i]-dfrac-1128)}
\Har^{[i]}*(\Har^{[i]}-1)*(\Har^{[i]}-\dfrac{-1}{128})*(\Har^{[i]}-\frac{15}{128})
&\equiv 0
\end{align}
for all $i=1,\ldots,\rankL$.
\begin{enumerate}
\item Let $M=M(1)(\theta)^{+}$.
Since $S_{ij}(1,1)_{1}\vac_{\tw}=0$ for any pair of distinct elements $i,j\in\{1,\ldots,\rankL\}$,
$S_{ij}(1,1)_{1}\lu\in \C v$ in $W$.
We note that $\omega^{[i]}_{1}\vac_{\tw}=(1/16)\vac_{\tw}$ and 
$\Har^{[i]}_{3}\vac_{\tw}=(-1/128)\vac_{\tw}$. 
By \eqref{eq:omega-1-omega-1/16-omega-9/16-2} and \eqref{eq:\Har[i]*(Har[i]-1)*(Har[i]-dfrac-1128)}, $\omega^{[i]}_{1}\lw=(1/16)\lw$ and 
$\Har^{[i]}_{3}\lw=(-1/128)\lw$ for all $\lw \in N(0)$.
We denote $\epsilon(S_{ij}(1,1), \lu)$ by $\eS$ for simplicity.
By \eqref{eq:s11-zeta-1},
\begin{align}
0&=\eS(11\eS^2-15\eS+6)(6\eS^3+6\eS^2-7\eS+1)
\label{eq:tw-0-S11}
\end{align}
and hence ${\eS}=0$. Thus, $S_{ij}(1,k)_{k}\lu=0$ for all $k\in \Z_{\geq 1}$ by
Lemma \ref{lemma:bound-index-S} (1)
and hence $N(0)\cong M(1)(\theta)^{+}(0)\oplus M(1)(\theta)^{+}(0)$ as $A(M(1)^{+})$-modules.

\item Let $M=M(1)^{-}$. We consider $N(0)$.
Since $N(0)/W(0)\cong M(1)^{-}(0)$, $A^{u}\cdot N(0)\neq 0$.
Since $A^{t}\cdot N(0)\subset W(0), A^{u}*A^{t}=0$, and $A^{u}\cdot \lw\neq 0$ for any non-zero $\lw\in W(0)$,
we  have $A^{t}\cdot N(0)=0$.
For any pair of distinct elements $i,j\in\{1,\ldots,\rankL\}$,
since  $A^{u}*\Lambda_{ij}=0$ by \cite[Proposition 5.3.12]{DN2001}, 
the same argument shows that $\Lambda_{ij}\cdot N(0)=0$.
We note that the eigenvalues for $\omega^{[i]}|_{N(0)}$ are $0$ or $1$,
and  those for $\Har^{[i]}|_{N(0)}$ are also $0$ or $1$.
We take a non-zero $\lv\in\module(0)$ and $\lu\in \mN$ as in \eqref{eq:eigenvalue-omega-H}, \eqref{eq:keromega-kerH},
and	\eqref{eq:pi(lu)lvmboxand}.	
We fix $i=1,\ldots,\rankL$.
By \eqref{eq:\Har[i]*(Har[i]-1)*(Har[i]-dfrac-1128)},
\begin{align}
\Har^{[i]}_{3}\lu=\lu\mbox{ or }\Har^{[i]}_{3}\lu=0.
\end{align}
\textcolor{black}{We study the following three cases:}
\begin{enumerate}
\item[(2-1)]
If $\Har^{[i]}_{3}\lu=\lu$, then it follows from \eqref{eq:omega-1-omega-1/16-omega-9/16-2} that
$\omega^{[i]}_{1}\lu=\lu$. 
\item[(2-2)]
The case that $\Har^{[i]}_{3}\lu=0$ and $(\omega^{[i]}_{1}-1)^2\lu=0$. Since
$(\omega^{[i]}_{1}-1)\lu\in \mW(0)$ and
there is no non-zero vector $\lw\in M(1)^{-}(0)$ such that $\omega^{[i]}_{1}\lw=\lw$ and $\Har^{[i]}_{3}\lw=0$, we have
$(\omega^{[i]}_{1}-1)\lu=0$.
\item[(2-3)]
The case that $\Har^{[i]}_{3}\lu=0$ and $(\omega^{[i]}_{1})^2\lu=0$.
Since $0\neq \lu\in M(0)\cong M(1)^{-}(0)$, there exists $k$ such that $\Har^{[k]}_{3}\lu\neq 0$. 
The argument (2-1) above shows that $\Har^{[k]}_{3}\lu=\lu$ and  $\omega^{[k]}_{1}\lu=\lu$. 
Since $\omega^{[k]}_{1}\omega^{[i]}_{1}\lu=\omega^{[i]}_{1}\lu$, we have $E^{u}_{kk}\omega^{[i]}_{1}\lu=\omega^{[i]}_{1}\lu$ in $\mW(0)$.
By \cite[Lemma 5.2.2]{DN2001}, $\omega^{[i]}_{1}\lu=0$. 
\end{enumerate}
Thus $A(M(1)^{+})\cdot \lu=A^{u}\cdot \lu$. 
Since $A^{u}$ is isomorphic to the matrix algebra, 
$A^{u}\cdot \lu$ is an irreducible  $A(M(1)^{+})$-module. Thus $N(0)\cong M(1)^{-}(0)\oplus M(1)^{-}(0)$.

\item In the case of $M=M(1)(\theta)^{-}$, the same argument as in (2) above shows that 
$N(0)\cong M(1)(\theta)^{-}(0)\oplus M(1)(\theta)^{-}(0)$.

\end{enumerate}
\end{proof}

By Lemmas \ref{lemma:Ext-split-lambda-mu}, \ref{lemma:Ext-M-M}, \cite[Proposition 2.5]{Abe2005}, and  \cite[Proposition 3.5]{ADL2005}, 
we have the following result:
\begin{proposition}
\label{proposition:Ext-total}
If a pair $(M,W)$ of irreducible $M(1)^{+}$-modules satisfies one of the following conditions,
then $\Ext^{1}_{M(1)^{+}}(M,W)=\Ext^{1}_{M(1)^{+}}(W,M)=0$.
 
\begin{enumerate}
\item $M\cong M(1,\lambda)$ with $\lambda\in\fh\setminus\{0\}$ and $W\not\cong M(1,\lambda)$.
\item $M\cong M(1)(\theta)^{\pm}$.
\item $M\cong M(1)^{+}$ and $W\not\cong M(1)^{-}$.
\item $M\cong M(1)^{-}$ and $W\not\cong M(1)^{+}$.
\end{enumerate}
\end{proposition}

The following result is a direct consequence of Lemmas \ref{lemma:M1lambda-submodule} and \ref{lemma:M1lambda-submodule-2}.
Here we call the $\N$-graded module $\bar{M}(U)$ in \cite[Theorem 6.2]{DLM1998t} the generalized Verma module associated with a module $U$ for the Zhu algebra.
\begin{corollary}\label{corollary:verma-irreducible}
	Let $\Omega$ be an irreducible $A(M(1)^{+})$-module such that 
	$\Omega\not\cong M(1)^{+}_{0}=\C\vac$.
Then the generalized Verma module for $M(1)^{+}$ associated
with  $\Omega$ is irreducible.
\end{corollary}
\begin{proof}
Let $\mN=\oplus_{i\in\delta+\Z_{\geq 0}}\mN_{i}$ with $\mN_{\delta}=\Omega$ be the generalized Verma module for $M(1)^{+}$ associated
with $\Omega$ and $\mW=\oplus_{i\in \gamma+\Z_{\geq 0}}\mW_{i}$ the maximal submodule of $\module$ such that $\Omega\cap \mW=0$.
We take $\gamma$ so that $W_{\gamma}\neq 0$ if $W\neq 0$.
We note that $\gamma-\delta\in\Z_{>0}$.
Taking the restricted dual of the exact sequence 
$	0\rightarrow \mW\overset{}{\rightarrow} \mN\overset{}{\rightarrow} \mN/\mW\rightarrow 0$,
we have the following exact sequence
\begin{align}
0\rightarrow (\mN/\mW)^{\prime}\overset{}{\rightarrow} 
\mN^{\prime}\overset{}{\rightarrow}\mW^{\prime} \rightarrow 0.
\end{align}
We note that $(\mN/\mW)^{\prime}\not\cong M(1)^{+}$ by \cite[Proposition 3.5]{ADL2005}.
Assume $\mW_{\gamma}\neq 0$ and let $B$ be an irreducible $A(M(1)^{+})$-submodule  
of $W_{\gamma}^{\prime}$.
By Lemmas \ref{lemma:M1lambda-submodule} and \ref{lemma:M1lambda-submodule-2},
there exists a non-zero
 $\lu^{\prime}\in \Omega_{M(1)^{+}}(N_{\gamma}^{\prime})$.
For any homogeneous element $a\in M(1)^{+}$ such that $\omega_{2}a=0$
and $i\in\Z_{\geq \wt a}$,
it follows from \cite[5.2.4]{FHL} that
\begin{align}
\label{eq:langleailuprime}
0&=\langle a_{i}\lu^{\prime},\lw\rangle=(-1)^{\wt a}\langle \lu^{\prime},a_{2\wt a-i-2}\lw\rangle
\end{align}
for all $\lw\in \mN$.
Since $\omega_2\omega^{[i]}=\omega_2J^{[i]}=0$ for all $i=1,\ldots,\rankL$,
it follows from Lemma \ref{lemma:m1plusmodulegenerated} and \eqref{eq:langleailuprime} that
$\lu^{\prime}=0$, a contradiction.
\end{proof}

\begin{lemma}
\label{lemma:generalizedVermaModule-M(1)-}
Let $\mW$ be the generalized Verma module associated to the $A(M(1)^{+})$-module $\C \vac$
and $\pi : \mW\rightarrow M(1)^{+}$ the canonical projection.
Then, $\Ker \pi\cong (M(1)^{-})^{\oplus k}$
as $M(1)^{+}$-modules  for some $k\in \{1,\ldots,\rankL\}$.
\end{lemma}
\begin{proof}
The same argument as in \cite[(6.1)]{Abe2005}
shows that there is a non-split exact sequence
\begin{align}
\label{eq:non-split-verma-M1plus}
0\rightarrow M(1)^{-}\overset{}{\rightarrow} \mN\overset{}{\rightarrow} M(1)^{+}\rightarrow 0
\end{align}
of $M(1)^{+}$-modules. Thus, $\Ker \pi \neq 0$.
Let $\lu\in\mW$ such that $\pi(\lu)=\vac$.
Note that $u\in \Omega_{M(1)^{+}}(N)$. 
Since $M(1)^{-}_{0}=0$, we have $\omega^{[k]}_1\lu=\Har^{[k]}_{3}\lu=S_{ij}(1,r)_{r}\lu=0$ 
for all $k=1,\ldots,\rankL$, pairs of distinct elements $i,j\in \{1,\ldots,\rankL\}$, and $r=1,2,3$.
For $i=1,\ldots,\rankL$, $P^{(8),\Har,i}$ denotes the element obtained by replacing $\omega$ by 
$\omega^{[i]}$ and $\Har$ by $\Har^{[i]}$ in $P^{(8),\Har}$ in \cite[(3.27)]{Tanabe2021-1}.
We  have shown in \cite[Lemma 3.5]{Tanabe2021-1} that
$\sv^{(8),\Har}=0$. A direct computation shows that
\begin{align}
\label{eq:P(8)H6lu=144}
0&=P^{(8),H,i}_{6}\lu=144(\omega^{[i]}_{0}-3\Har^{[i]}_{2})\lu
\end{align}
for all $i=1,\ldots,\rankL$.
Taking the $3$rd action of \eqref{eq:s11-3} on $\lu$, we have
\begin{align}
\label{eq:Si1121luSi1132}
0&=S_{ij}(1,2)_{1}\lu+S_{ij}(1,3)_{2}\lu
\end{align}
for any pair of distinct elements $i,j\in\{1,\ldots,\rankL\}$. 
By \eqref{eq:P(8)H6lu=144}, \eqref{eq:Si1121luSi1132}, 
and 
\begin{align}
S_{ij}(1,1)_{0}\lu&=-(\omega_{0}S_{ij}(1,1))_{1}\lu=-S_{ji}(1,2)_{1}\lu-S_{ij}(1,2)_{1}\lu
\end{align}
for any pair of distinct elements $i,j\in\{1,\ldots,\rankL\}$,
$N_1$ is spanned by $\{\omega^{[j]}_{0}\lu,S_{ij}(1,2)_{1}\lu\ |\ i,j=1,\ldots,\rankL, i\neq j\}$. 
For distinct $i,j,k\in\{1,\ldots,\rankL\}$, by \eqref{eq:P(8)H6lu=144}, \eqref{eq:Si1121luSi1132},
and commutation relations (see Remark \ref{remark:com-appendix}), a direct computation shows that
\begin{align}
\label{eqn:omega1[j]omega[j]0}
\omega_{1}^{[j]}\omega^{[j]}_{0}\lu&=\omega^{[j]}_{0}\lu,\nonumber\\
\omega_{1}^{[i]}\omega^{[j]}_{0}\lu&=0,\nonumber\\
\Har_{3}^{[j]}\omega^{[j]}_{0}\lu&=\omega^{[j]}_{0}\lu,\nonumber\\
\Har_{3}^{[i]}\omega^{[j]}_{0}\lu&=0,\nonumber\\
S_{ij}(1,1)_{1}\omega^{[j]}_{0}\lu&=-S_{ij}(1,2)_{1}\lu,\nonumber\\
S_{ij}(1,2)_{2}\omega^{[j]}_{0}\lu&=2S_{ij}(1,2)_{1}\lu,\nonumber\\
S_{ij}(1,3)_{3}\omega^{[j]}_{0}\lu&=-3S_{ij}(1,2)_{1}\lu,\nonumber\\
\omega^{[i]}_{1}S_{ij}(1,2)_{1}\lu&=S_{ij}(1,2)_{1}\lu,\nonumber\\
\omega^{[j]}_{1}S_{ij}(1,2)_{1}\lu&=0,\nonumber\\
\Har_{3}^{[i]}S_{ij}(1,2)_{1}\lu&=S_{ij}(1,2)_{1 }\lu,\nonumber\\
\Har_{3}^{[j]}S_{ij}(1,2)_{1}\lu&=0,\nonumber\\
S_{ij}(1,1)_{1}S_{ij}(1,2)_{1}\lu&=\textcolor{black}{-\omega^{[j]}_{0}\lu},\nonumber\\
S_{ij}(1,2)_{2}S_{ij}(1,2)_{1}\lu&=0,\nonumber\\
S_{ij}(1,3)_{3}S_{ij}(1,2)_{1}\lu&=0,\nonumber\\
S_{kj}(1,1)_{1}S_{ij}(1,2)_{1}\lu&=0,\nonumber\\
S_{kj}(1,2)_{2}S_{ij}(1,2)_{1}\lu&=0,\nonumber\\
S_{kj}(1,3)_{3}S_{ij}(1,2)_{1}\lu&=0,\nonumber\\
S_{ki}(1,1)_{1}S_{ij}(1,2)_{1}\lu&=S_{kj}(1,2)_{1}\lu,\nonumber\\
S_{ki}(1,2)_{2}S_{ij}(1,2)_{1}\lu&=-2S_{kj}(1,2)_{1}\lu,\nonumber\\
S_{ki}(1,3)_{3}S_{ij}(1,2)_{1}\lu&=3S_{kj}(1,2)_{1}\lu.
\end{align}
Thus, by \eqref{eq:norm2Si1(1,1)1lu=-Si1(1,2)2lu-2}, for each $j=1,\ldots,\rankL$,
the linear subjective map \\$M(1)^{-}(0)\rightarrow 
U^{\langle j\rangle}:=\Span_{\C}\{\omega^{[j]}_{0}\lu,S_{ij}(1,2)_{1}\lu\ |\ i\neq j\}$
sending \textcolor{black}{$h^{[j]}(-1)\vac$ to $\omega^{[j]}_{0}\lu$ and $h^{[i]}(-1)\vac$ to $-S_{ij}(1,2)_{1}\lu$ }for $i\neq j$
is an $A(M(1)^{+})$-homomorphism.
Since $M(1)^{-}(0)$ is an irreducible $A(M(1)^{+})$-module, if $U^{\langle j\rangle}\neq 0$, then $U^{\langle j\rangle}\cong M(1)^{-}(0)$
as $A(M(1)^{+})$-modules.
Since $\sum_{j=1}^{\rankL}U^{\langle j\rangle}=N_1$, $(W/(M(1)^{+}\cdot (\sum_{j=1}^{\rankL}U^{\langle j\rangle}))_{1}=0$.
Thus $(W/(M(1)^{+}\cdot (\sum_{j=1}^{\rankL}U^{\langle j\rangle}))\cong M(1)^{+}$ and hence
$\Ker \pi=M(1)^{+}\cdot (\sum_{j=1}^{\rankL}U^{\langle j\rangle})$.
Now the result follows from Corollary \ref{corollary:verma-irreducible}.
\end{proof}

By Lemma \ref{lemma:Zhu-Omega}, Proposition \ref{proposition:M(1)plusmoduleinweak}, 
Corollary \ref{corollary:verma-irreducible}, and Lemma \ref{lemma:generalizedVermaModule-M(1)-},
we have the following result.
\begin{corollary}
\label{corollary:M(1)plusirreduciblemoduleinweak}
Let $\lattice$ be a non-degenerate even lattice of finite rank
and $\module$ a non-zero weak $V_{\lattice}^{+}$-module.
Then, there exists an irreducible $M(1)^{+}$-submodule of $\module$. 
\end{corollary}
\begin{proof}
By  Lemma \ref{lemma:Zhu-Omega}, there exists an irreducible $A(M(1)^{+})$-submodule $\Omega$
of $\Omega_{M(1)^{+}}(M)$. 
Let $\mN$ be the generalized Verma module for $M(1)^{+}$ associated
with $\Omega$ (Proposition \ref{proposition:M(1)plusmoduleinweak})
and $f : \mN\rightarrow \module$ the associated $M(1)^{+}$-homomorphism.
If $\Omega\not\cong \C\vac$, then  by Corollary \ref{corollary:verma-irreducible},
$N$ is irreducible and hence so is $f(N)$.
If $\Omega\cong \C\vac$, then  by  Lemma \ref{lemma:generalizedVermaModule-M(1)-},
$M(1)^{-}\subset f(N)$ or $f(N)\cong M(1)^{+}$.
This completes the proof.
\end{proof}

\appendix
\renewcommand{\thesection}{Appendix A\arabic{section}}
\setcounter{section}{1}
\renewcommand{\thesubsection}{A\arabic{section}-\arabic{subsection}}
\section{}
\label{section:appendix}
\renewcommand{\thesection}{A\arabic{section}}

In this appendix, for some $a,b\in M(1)^{+}$, we put the computations of 
$a_{k}b$ for $k\in\Z_{\geq 0}$.
For $k\in\Z_{\geq 0}$ not listed below, $a_{k}b=0$.
Using these results, we can compute the commutation relation $[a_i,b_j]=\sum_{k=0}^{\infty}\binom{i}{k}(a_kb)_{i+j-k}$.
Let $h^{[1]},\ldots,h^{[\rankL]}$ be an orthonormal basis of $\fh$.
The rest of this appendix, $i,j,k,l$ are distinct elements of $\{1,\ldots,\rankL\}$.
\label{section:appendix-The general case}

\begin{align}
\label{eq:omega[j]0Sij(1,1)=Sij(1,2)-1vac}
	\omega^{[j]}_{0}S_{ij}(1,1)&=
	S_{ij}(1,2)_{-1 } \vac
	,&
	\omega^{[j]}_{1}S_{ij}(1,1)&=
	S_{ij}(1,1)_{-1 } \vac
	,
\end{align}

\begin{align}
	\omega^{[j]}_{0}S_{ij}(1,2)&=
	2S_{ij}(1,3)_{-1 } \vac
	,
&
	\omega^{[j]}_{1}S_{ij}(1,2)&=
	2S_{ij}(1,2)_{-1 } \vac
	,
&
	\omega^{[j]}_{2}S_{ij}(1,2)&=
	2S_{ij}(1,1)_{-1 } \vac
	,
\end{align}

\begin{align}
	\omega^{[j]}_{0}S_{ij}(1,3)&=
	-\omega^{[j]}_{-2 } S_{ij}(1,1)_{-1 } \vac
	+2\omega^{[j]}_{-1 } S_{ij}(1,2)_{-1 } \vac
	,
	\nonumber\\
	\omega^{[j]}_{1}S_{ij}(1,3)&=
	3S_{ij}(1,3)_{-1 } \vac
	,
	\nonumber\\
	\omega^{[j]}_{2}S_{ij}(1,3)&=
	3S_{ij}(1,2)_{-1 } \vac
	,
	\nonumber\\
	\omega^{[j]}_{3}S_{ij}(1,3)&=
	3S_{ij}(1,1)_{-1 } \vac
	,
\end{align}

\begin{align}
\label{eq:Har[j]0Sij(1,1)=-2omega[j]-2Sij(1,1)-1vac}
	\Har^{[j]}_{0}S_{ij}(1,1)&=
	-2\omega^{[j]}_{-2 } S_{ij}(1,1)_{-1 } \vac
	+4\omega^{[j]}_{-1 } S_{ij}(1,2)_{-1 } \vac
	,
	\nonumber\\
	\Har^{[j]}_{1}S_{ij}(1,1)&=
	4S_{ij}(1,3)_{-1 } \vac
	,
	\nonumber\\
	\Har^{[j]}_{2}S_{ij}(1,1)&=
	\frac{7}{3}S_{ij}(1,2)_{-1 } \vac
	,
	\nonumber\\
	\Har^{[j]}_{3}S_{ij}(1,1)&=
	S_{ij}(1,1)_{-1 } \vac
	,
\end{align}

\begin{align}
	\Har^{[j]}_{0}S_{ij}(1,2)&=
	-6\omega_{0 } \omega^{[j]}_{-2 } S_{ij}(1,1)_{-1 } \vac
	+6\omega^{[i]}_{-2 } S_{ij}(1,2)_{-1 } \vac
	\nonumber\\&\quad{}
	-4\omega_{0 } \omega^{[i]}_{-1 } S_{ij}(1,2)_{-1 } \vac
	+12\omega_{0 } \omega^{[j]}_{-1 } S_{ij}(1,2)_{-1 } \vac
	\nonumber\\&\quad{}
	+\omega_{0 } \omega_{0 } \omega_{0 } S_{ij}(1,2)_{-1 } \vac
	+8\omega^{[i]}_{-1 } S_{ij}(1,3)_{-1 } \vac
	\nonumber\\&\quad{}
	-6\omega_{0 } \omega_{0 } S_{ij}(1,3)_{-1 } \vac
	,
	\nonumber\\
	\Har^{[j]}_{1}S_{ij}(1,2)&=
	-6\omega^{[j]}_{-2 } S_{ij}(1,1)_{-1 } \vac
	+12\omega^{[j]}_{-1 } S_{ij}(1,2)_{-1 } \vac
	,
	\nonumber\\
	\Har^{[j]}_{2}S_{ij}(1,2)&=
	\frac{38}{3}S_{ij}(1,3)_{-1 } \vac
	,
	\nonumber\\
	\Har^{[j]}_{3}S_{ij}(1,2)&=
	8S_{ij}(1,2)_{-1 } \vac
	,
	\nonumber\\
	\Har^{[j]}_{4}S_{ij}(1,2)&=
	4S_{ij}(1,1)_{-1 } \vac
	,
\end{align}

\begin{align}
	\Har^{[j]}_{0}S_{ij}(1,3)&=
	\frac{40}{29}\omega^{[j]}_{-3 } S_{ij}(1,2)_{-1 } \vac
	+\frac{60}{29}\Har^{[j]}_{-1 } S_{ij}(1,2)_{-1 } \vac
	\nonumber\\&\quad{}
	+\frac{-60}{29}\omega^{[j]}_{-2 } S_{ij}(1,3)_{-1 } \vac
	,
	\nonumber\\
	\Har^{[j]}_{1}S_{ij}(1,3)&=
	-12\omega_{0 } \omega^{[j]}_{-2 } S_{ij}(1,1)_{-1 } \vac
	+12\omega^{[i]}_{-2 } S_{ij}(1,2)_{-1 } \vac
	\nonumber\\&\quad{}
	-8\omega_{0 } \omega^{[i]}_{-1 } S_{ij}(1,2)_{-1 } \vac
	+24\omega_{0 } \omega^{[j]}_{-1 } S_{ij}(1,2)_{-1 } \vac	\nonumber\\&\quad{}
	+2\omega_{0 } \omega_{0 } \omega_{0 } S_{ij}(1,2)_{-1 } \vac
	+16\omega^{[i]}_{-1 } S_{ij}(1,3)_{-1 } \vac	\nonumber\\&\quad{}
	-12\omega_{0 } \omega_{0 } S_{ij}(1,3)_{-1 } \vac
	,
	\nonumber\\
	\Har^{[j]}_{2}S_{ij}(1,3)&=
	\frac{-37}{3}\omega^{[j]}_{-2 } S_{ij}(1,1)_{-1 } \vac
	+\frac{74}{3}\omega^{[j]}_{-1 } S_{ij}(1,2)_{-1 } \vac
	,
	\nonumber\\
	\Har^{[j]}_{3}S_{ij}(1,3)&=
	27S_{ij}(1,3)_{-1 } \vac
	,
	\nonumber\\
	\Har^{[j]}_{4}S_{ij}(1,3)&=
	18S_{ij}(1,2)_{-1 } \vac
	,
	\nonumber\\
	\Har^{[j]}_{5}S_{ij}(1,3)&=
	10S_{ij}(1,1)_{-1 } \vac
	,
\end{align}

\begin{align}
	\omega^{[i]}_{0}S_{ij}(1,1)&=
	\omega_{0 } \textcolor{black}{S_{ij}}(1,1)_{-1 } \vac
	-S_{ij}(1,2)_{-1 } \vac
	,
	\nonumber\\
	\omega^{[i]}_{1}S_{ij}(1,1)&=
	S_{ij}(1,1)_{-1 } \vac
	,
\end{align}

\begin{align}
	\omega^{[i]}_{0}S_{ij}(1,2)&=
	\omega_{0 } S_{ij}(1,2)_{-1 } \vac
	-2S_{ij}(1,3)_{-1 } \vac
	,
	\nonumber\\
	\omega^{[i]}_{1}S_{ij}(1,2)&=
	S_{ij}(1,2)_{-1 } \vac
	,
\end{align}

\begin{align}
	\omega^{[i]}_{0}S_{ij}(1,3)&=
	\omega^{[j]}_{-2 } S_{ij}(1,1)_{-1 } \vac
	-2\omega^{[j]}_{-1 } S_{ij}(1,2)_{-1 } \vac
	+\omega_{0 } S_{ij}(1,3)_{-1 } \vac
	,
	\nonumber\\
	\omega^{[i]}_{1}S_{ij}(1,3)&=
	S_{ij}(1,3)_{-1 } \vac
	,
\end{align}

\begin{align}
	\Har^{[i]}_{0}S_{ij}(1,1)&=
	2\omega^{[j]}_{-2 } S_{ij}(1,1)_{-1 } \vac
	+\omega_{0 } \omega_{0 } \omega_{0 } S_{ij}(1,1)_{-1 } \vac
	\nonumber\\&\quad{}
	-4\omega^{[j]}_{-1 } S_{ij}(1,2)_{-1 } \vac
	-3\omega_{0 } \omega_{0 } S_{ij}(1,2)_{-1 } \vac
	\nonumber\\&\quad{}
	+6\omega_{0 } S_{ij}(1,3)_{-1 } \vac
	,
	\nonumber\\
	\Har^{[i]}_{1}S_{ij}(1,1)&=
	2\omega_{0 } \omega_{0 } S_{ij}(1,1)_{-1 } \vac
	-4\omega_{0 } S_{ij}(1,2)_{-1 } \vac
	+4S_{ij}(1,3)_{-1 } \vac
	,
	\nonumber\\
	\Har^{[i]}_{2}S_{ij}(1,1)&=
	\frac{7}{3}\omega_{0 } S_{ij}(1,1)_{-1 } \vac
	+\frac{-7}{3}S_{ij}(1,2)_{-1 } \vac
	,
	\nonumber\\
	\Har^{[i]}_{3}S_{ij}(1,1)&=
	S_{ij}(1,1)_{-1 } \vac
	,
\end{align}

\begin{align}
	\Har^{[i]}_{0}S_{ij}(1,2)&=
	-6\omega^{[i]}_{-2 } S_{ij}(1,2)_{-1 } \vac
	+4\omega_{0 } \omega^{[i]}_{-1 } S_{ij}(1,2)_{-1 } \vac
	-8\omega^{[i]}_{-1 } S_{ij}(1,3)_{-1 } \vac
	,
	\nonumber\\
	\Har^{[i]}_{1}S_{ij}(1,2)&=
	-4\omega^{[j]}_{-2 } S_{ij}(1,1)_{-1 } \vac
	+8\omega^{[j]}_{-1 } S_{ij}(1,2)_{-1 } \vac
	\nonumber\\&\quad{}
	+2\omega_{0 } \omega_{0 } S_{ij}(1,2)_{-1 } \vac
	-8\omega_{0 } S_{ij}(1,3)_{-1 } \vac
	,
	\nonumber\\
	\Har^{[i]}_{2}S_{ij}(1,2)&=
	\frac{7}{3}\omega_{0 } S_{ij}(1,2)_{-1 } \vac
	+\frac{-14}{3}S_{ij}(1,3)_{-1 } \vac
	,
	\nonumber\\
	\Har^{[i]}_{3}S_{ij}(1,2)&=
	S_{ij}(1,2)_{-1 } \vac
	,
\end{align}

\begin{align}
	\Har^{[i]}_{0}S_{ij}(1,3)&=
	\frac{-4}{3}\Har^{[i]}_{-1 } S_{ij}(1,2)_{-1 } \vac
	+\frac{16}{9}\omega^{[i]}_{-3 } S_{ij}(1,2)_{-1 } \vac
	\nonumber\\&\quad{}
	+\frac{-11}{3}\omega_{0 } \omega^{[i]}_{-2 } S_{ij}(1,2)_{-1 } \vac
	+2\omega_{0 } \omega_{0 } \omega^{[i]}_{-1 } S_{ij}(1,2)_{-1 } \vac
	\nonumber\\&\quad{}
	+\frac{4}{3}\omega^{[i]}_{-2 } S_{ij}(1,3)_{-1 } \vac
	-4\omega_{0 } \omega^{[i]}_{-1 } S_{ij}(1,3)_{-1 } \vac
	,
	\nonumber\\
	\Har^{[i]}_{1}S_{ij}(1,3)&=
	-2\omega_{0 } \omega^{[j]}_{-2 } S_{ij}(1,1)_{-1 } \vac
	+6\omega^{[i]}_{-2 } S_{ij}(1,2)_{-1 } \vac
	\nonumber\\&\quad{}
	-4\omega_{0 } \omega^{[i]}_{-1 } S_{ij}(1,2)_{-1 } \vac
	+4\omega_{0 } \omega^{[j]}_{-1 } S_{ij}(1,2)_{-1 } \vac
	\nonumber\\&\quad{}
	+\omega_{0 } \omega_{0 } \omega_{0 } S_{ij}(1,2)_{-1 } \vac
	+8\omega^{[i]}_{-1 } S_{ij}(1,3)_{-1 } \vac
	\nonumber\\&\quad{}
	-4\omega_{0 } \omega_{0 } S_{ij}(1,3)_{-1 } \vac
	,
	\nonumber\\
	\Har^{[i]}_{2}S_{ij}(1,3)&=
	\frac{7}{3}\omega^{[j]}_{-2 } S_{ij}(1,1)_{-1 } \vac
	+\frac{-14}{3}\omega^{[j]}_{-1 } S_{ij}(1,2)_{-1 } \vac
	+\frac{7}{3}\omega_{0 } S_{ij}(1,3)_{-1 } \vac
	,
	\nonumber\\
	\Har^{[i]}_{3}S_{ij}(1,3)&=
	S_{ij}(1,3)_{-1 } \vac
	,
\end{align}

\begin{align}
	\label{eq:Sij(1,1)0Sij(1,1)-1}
S_{ij}(1,1)_{0}S_{ij}(1,1)&=
\omega_{0 } \omega^{[i]}_{-1 } \vac
+\omega_{0 } \omega^{[j]}_{-1 } \vac
,
\nonumber\\
S_{ij}(1,1)_{1}S_{ij}(1,1)&=
2\omega^{[i]}_{-1 } \vac
+2\omega^{[j]}_{-1 } \vac
,
\nonumber\\
S_{ij}(1,1)_{2}S_{ij}(1,1)&=0,\nonumber\\
S_{ij}(1,1)_{3}S_{ij}(1,1)&=
\vac
,
\end{align}

\begin{align}
	\label{eq:Sij(1,1)0Sij(1,1)-2}
S_{ij}(1,1)_{0}S_{ij}(1,2)&=
2\Har^{[i]}_{-1 } \vac
-2\Har^{[j]}_{-1 } \vac
+2\omega_{0 }^{2}\omega^{[i]}_{-1 } \vac
+\omega_{0 }^{2}\omega^{[j]}_{-1 } \vac
,
\nonumber\\
S_{ij}(1,1)_{1}S_{ij}(1,2)&=
2\omega_{0 } \omega^{[i]}_{-1 } \vac
+\omega_{0 } \omega^{[j]}_{-1 } \vac
,
\nonumber\\
S_{ij}(1,1)_{2}S_{ij}(1,2)&=
4\omega^{[i]}_{-1 } \vac
,
\nonumber\\
S_{ij}(1,1)_{3}S_{ij}(1,2)&=0,\nonumber\\
S_{ij}(1,1)_{4}S_{ij}(1,2)&=
2\vac
,
\end{align}

\begin{align}
	\label{eq:Sij(1,1)0Sij(1,1)-3}
S_{ij}(1,1)_{0}S_{ij}(1,3)&=
3\omega_{0 } \Har^{[i]}_{-1 } \vac
-\omega_{0 } \Har^{[j]}_{-1 } \vac
+\omega_{0 }^{3}\omega^{[i]}_{-1 } \vac
+\omega_{0 }^{3}\omega^{[j]}_{-1 } \vac
,
\nonumber\\
S_{ij}(1,1)_{1}S_{ij}(1,3)&=
3\Har^{[i]}_{-1 } \vac
+\Har^{[j]}_{-1 } \vac
+\omega_{0 }^{2}\omega^{[i]}_{-1 } \vac
+\omega_{0 }^{2}\omega^{[j]}_{-1 } \vac
,
\nonumber\\
S_{ij}(1,1)_{2}S_{ij}(1,3)&=
3\omega_{0 } \omega^{[i]}_{-1 } \vac
,
\nonumber\\
S_{ij}(1,1)_{3}S_{ij}(1,3)&=
6\omega^{[i]}_{-1 } \vac
,
\nonumber\\
S_{ij}(1,1)_{4}S_{ij}(1,3)&=0,\nonumber\\
S_{ij}(1,1)_{5}S_{ij}(1,3)&=
3\vac
,
\end{align}

\begin{align}
S_{ij}(1,1)_{0}S_{k j}(1,1)&=
S_{k i}(1,2)_{-1 } \vac
,
\nonumber\\
S_{ij}(1,1)_{1}S_{k j}(1,1)&=
S_{k i}(1,1)_{-1 } \vac
,
\end{align}

\begin{align}
S_{ij}(1,1)_{0}S_{k j}(1,2)&=
2S_{k i}(1,3)_{-1 } \vac
,
\nonumber\\
S_{ij}(1,1)_{1}S_{k j}(1,2)&=
2S_{k i}(1,2)_{-1 } \vac
,
\nonumber\\
S_{ij}(1,1)_{2}S_{k j}(1,2)&=
2S_{k i}(1,1)_{-1 } \vac
,
\end{align}

\begin{align}
S_{ij}(1,1)_{0}S_{k j}(1,3)&=
-3\omega^{[k]}_{-1 } S_{k i}(1,1)_{-2 } \vac
+\omega_{0 } \omega^{[k]}_{-1 } S_{k i}(1,1)_{-1 } \vac
\nonumber\\&\quad{}
+\omega_{0 }^{3}S_{k i}(1,1)_{-1 } \vac
+2\omega^{[k]}_{-1 } S_{k i}(1,2)_{-1 } \vac
\nonumber\\&\quad{}
-3S_{k i}(1,2)_{-3 } \vac
+3\omega_{0 } S_{k i}(1,3)_{-1 } \vac
,
\nonumber\\
S_{ij}(1,1)_{1}S_{k j}(1,3)&=
3S_{k i}(1,3)_{-1 } \vac
,
\nonumber\\
S_{ij}(1,1)_{2}S_{k j}(1,3)&=
3S_{k i}(1,2)_{-1 } \vac
,
\nonumber\\
S_{ij}(1,1)_{3}S_{k j}(1,3)&=
3S_{k i}(1,1)_{-1 } \vac
,
\end{align}

\begin{align}
S_{ij}(1,1)_{0}S_{k i}(1,1)&=
S_{k j}(1,2)_{-1 } \vac
,
&
S_{ij}(1,1)_{1}S_{k i}(1,1)&=
S_{k j}(1,1)_{-1 } \vac
,
\end{align}

\begin{align}
S_{ij}(1,1)_{0}S_{k i}(1,2)&=
2S_{k j}(1,3)_{-1 } \vac
,
\nonumber\\
S_{ij}(1,1)_{1}S_{k i}(1,2)&=
2S_{k j}(1,2)_{-1 } \vac
,
\nonumber\\
S_{ij}(1,1)_{2}S_{k i}(1,2)&=
2S_{k j}(1,1)_{-1 } \vac
,
\end{align}

\begin{align}
S_{ij}(1,1)_{0}S_{k i}(1,3)&=
3\omega^{[k]}_{-2 } S_{k j}(1,1)_{-1 } \vac
-2\omega_{0 } \omega^{[k]}_{-1 } S_{k j}(1,1)_{-1 } \vac
\nonumber\\&\quad{}
+\omega_{0 }^{3}S_{k j}(1,1)_{-1 } \vac
+2\omega^{[k]}_{-1 } S_{k j}(1,2)_{-1 } \vac
\nonumber\\&\quad{}
-3\omega_{0 } S_{k j}(1,2)_{-2 } \vac
+3\omega_{0 } S_{k j}(1,3)_{-1 } \vac
,
\nonumber\\
S_{ij}(1,1)_{1}S_{k i}(1,3)&=
3S_{k j}(1,3)_{-1 } \vac
,
\nonumber\\
S_{ij}(1,1)_{2}S_{k i}(1,3)&=
3S_{k j}(1,2)_{-1 } \vac
,
\nonumber\\
S_{ij}(1,1)_{3}S_{k i}(1,3)&=
3S_{k j}(1,1)_{-1 } \vac
,
\end{align}

\begin{align}
S_{ij}(1,2)_{0}S_{ij}(1,2)&=
-3\omega_{0 } \Har^{[i]}_{-1 } \vac
-\omega_{0 } \Har^{[j]}_{-1 } \vac
-\omega_{0 }^{3}\omega^{[i]}_{-1 } \vac
+\omega_{0 }^{3}\omega^{[j]}_{-1 } \vac
,
\nonumber\\
S_{ij}(1,2)_{1}S_{ij}(1,2)&=
-6\Har^{[i]}_{-1 } \vac
-2\Har^{[j]}_{-1 } \vac
-2\omega_{0 }^{2}\omega^{[i]}_{-1 } \vac
+\omega_{0 }^{2}\omega^{[j]}_{-1 } \vac
,
\nonumber\\
S_{ij}(1,2)_{2}S_{ij}(1,2)&=
-6\omega_{0 } \omega^{[i]}_{-1 } \vac
,
\nonumber\\
S_{ij}(1,2)_{3}S_{ij}(1,2)&=
-12\omega^{[i]}_{-1 } \vac
,
\nonumber\\
S_{ij}(1,2)_{4}S_{ij}(1,2)&=0,\nonumber\\
S_{ij}(1,2)_{5}S_{ij}(1,2)&=
-6\vac
,
\end{align}

\begin{align}
S_{ij}(1,2)_{0}S_{ij}(1,3)&=
-8\omega^{[i]}_{-1 } \Har^{[i]}_{-1 } \vac
-2\omega^{[i]}_{-2 } \omega^{[i]}_{-2 } \vac
+2\omega^{[j]}_{-2 } \omega^{[j]}_{-2 } \vac
+8\omega^{[j]}_{-1 } \Har^{[j]}_{-1 } \vac
\nonumber\\&\quad{}
-7\omega_{0 }^{2}\Har^{[i]}_{-1 } \vac
+2\omega_{0 }^{2}\omega^{[i]}_{-1 } \omega^{[i]}_{-1 } \vac
-2\omega_{0 }^{2}\omega^{[j]}_{-1 } \omega^{[j]}_{-1 } \vac
-3\omega_{0 }^{2}\Har^{[j]}_{-1 } \vac
\nonumber\\&\quad{}
-19\omega_{0 }^{4}\omega^{[i]}_{-1 } \vac
+2\omega_{0 }^{4}\omega^{[j]}_{-1 } \vac
,
\nonumber\\
S_{ij}(1,2)_{1}S_{ij}(1,3)&=
-6\omega_{0 } \Har^{[i]}_{-1 } \vac
-\omega_{0 } \Har^{[j]}_{-1 } \vac
-\omega_{0 }^{3}\omega^{[i]}_{-1 } \vac
+\omega_{0 }^{3}\omega^{[j]}_{-1 } \vac
,
\nonumber\\
S_{ij}(1,2)_{2}S_{ij}(1,3)&=
-12\Har^{[i]}_{-1 } \vac
-4\omega_{0 }^{2}\omega^{[i]}_{-1 } \vac
,
\nonumber\\
S_{ij}(1,2)_{3}S_{ij}(1,3)&=
-12\omega_{0 } \omega^{[i]}_{-1 } \vac
,
\nonumber\\
S_{ij}(1,2)_{4}S_{ij}(1,3)&=
-24\omega^{[i]}_{-1 } \vac
,
\nonumber\\
S_{ij}(1,2)_{5}S_{ij}(1,3)&=0,\nonumber\\
S_{ij}(1,2)_{6}S_{ij}(1,3)&=
-12\vac
,
\end{align}

\begin{align}
S_{ij}(1,2)_{0}S_{k j}(1,1)&=
-2S_{k i}(1,3)_{-1 } \vac
,
\nonumber\\
S_{ij}(1,2)_{1}S_{k j}(1,1)&=
-2S_{k i}(1,2)_{-1 } \vac
,
\nonumber\\
S_{ij}(1,2)_{2}S_{k j}(1,1)&=
-2S_{k i}(1,1)_{-1 } \vac
,
\end{align}

\begin{align}
S_{ij}(1,2)_{0}S_{k j}(1,2)&=
6\omega^{[k]}_{-1 } S_{k i}(1,1)_{-2 } \vac
-2\omega_{0 } \omega^{[k]}_{-1 } S_{k i}(1,1)_{-1 } \vac
\nonumber\\&\quad{}
 -\omega_{0 }^{3}S_{k i}(1,1)_{-1 } \vac
-4\omega^{[k]}_{-1 } S_{k i}(1,2)_{-1 } \vac
\nonumber\\&\quad{}
+6S_{k i}(1,2)_{-3 } \vac
-6\omega_{0 } S_{k i}(1,3)_{-1 } \vac
,
\nonumber\\
S_{ij}(1,2)_{1}S_{k j}(1,2)&=
-6S_{k i}(1,3)_{-1 } \vac
,
\nonumber\\
S_{ij}(1,2)_{2}S_{k j}(1,2)&=
-6S_{k i}(1,2)_{-1 } \vac
,
\nonumber\\
S_{ij}(1,2)_{3}S_{k j}(1,2)&=
-6S_{k i}(1,1)_{-1 } \vac
,
\end{align}

\begin{align}
S_{ij}(1,2)_{0}S_{k j}(1,3)&=
16\omega^{[k]}_{-1 } S_{k i}(1,1)_{-3 } \vac
-4\Har^{[k]}_{-1 } S_{k i}(1,1)_{-1 } \vac
\nonumber\\&\quad{}
+44\omega_{0 } \omega^{[k]}_{-1 } S_{k i}(1,1)_{-2 } \vac
+2\omega_{0 } \omega^{[i]}_{-2 } S_{k i}(1,1)_{-1 } \vac
\nonumber\\&\quad{}
-16\omega_{0 }^{2}\omega^{[k]}_{-1 } S_{k i}(1,1)_{-1 } \vac
 -\omega_{0 }^{4}S_{k i}(1,1)_{-1 } \vac
\nonumber\\&\quad{}
-10\omega^{[k]}_{-2 } S_{k i}(1,2)_{-1 } \vac
-4\omega^{[k]}_{-1 } S_{k i}(1,2)_{-2 } \vac
\nonumber\\&\quad{}
+2\omega_{0 } S_{k i}(1,2)_{-3 } \vac
-4\omega_{0 } \omega^{[i]}_{-1 } S_{k i}(1,2)_{-1 } \vac
,
\nonumber\\
S_{ij}(1,2)_{1}S_{k j}(1,3)&=
12\omega^{[k]}_{-1 } S_{k i}(1,1)_{-2 } \vac
-4\omega_{0 } \omega^{[k]}_{-1 } S_{k i}(1,1)_{-1 } \vac
\nonumber\\&\quad{}
-2\omega_{0 }^{3}S_{k i}(1,1)_{-1 } \vac
-8\omega^{[k]}_{-1 } S_{k i}(1,2)_{-1 } \vac
\nonumber\\&\quad{}
+12S_{k i}(1,2)_{-3 } \vac
-12\omega_{0 } S_{k i}(1,3)_{-1 } \vac
,
\nonumber\\
S_{ij}(1,2)_{2}S_{k j}(1,3)&=
-12S_{k i}(1,3)_{-1 } \vac
,
\nonumber\\
S_{ij}(1,2)_{3}S_{k j}(1,3)&=
-12S_{k i}(1,2)_{-1 } \vac
,
\nonumber\\
S_{ij}(1,2)_{4}S_{k j}(1,3)&=
-12S_{k i}(1,1)_{-1 } \vac
,
\end{align}

\begin{align}
S_{ij}(1,2)_{0}S_{k i}(1,1)&=
2S_{k j}(1,3)_{-1 } \vac
,
&
S_{ij}(1,2)_{1}S_{k i}(1,1)&=
S_{k j}(1,2)_{-1 } \vac
,
\end{align}

\begin{align}
S_{ij}(1,2)_{0}S_{k i}(1,2)&=
6\omega^{[k]}_{-2 } S_{k j}(1,1)_{-1 } \vac
-4\omega_{0 } \omega^{[k]}_{-1 } S_{k j}(1,1)_{-1 } \vac
\nonumber\\&\quad{}
+\omega_{0 }^{3}S_{k j}(1,1)_{-1 } \vac
+4\omega^{[k]}_{-1 } S_{k j}(1,2)_{-1 } \vac
\nonumber\\&\quad{}
-3\omega_{0 } S_{k j}(1,2)_{-2 } \vac
+6\omega_{0 } S_{k j}(1,3)_{-1 } \vac
,
\nonumber\\
S_{ij}(1,2)_{1}S_{k i}(1,2)&=
4S_{k j}(1,3)_{-1 } \vac
,
\nonumber\\
S_{ij}(1,2)_{2}S_{k i}(1,2)&=
2S_{k j}(1,2)_{-1 } \vac
,
\end{align}

\begin{align}
S_{ij}(1,2)_{0}S_{k i}(1,3)&=
-16\omega^{[k]}_{-1 } S_{k j}(1,1)_{-3 } \vac
+4\Har^{[k]}_{-1 } S_{k j}(1,1)_{-1 } \vac
\nonumber\\&\quad{}
-98\omega_{0 } \omega^{[k]}_{-1 } S_{k j}(1,1)_{-2 } \vac
+34\omega_{0 }^{2}\omega^{[k]}_{-1 } S_{k j}(1,1)_{-1 } \vac
\nonumber\\&\quad{}
+3\omega_{0 }^{4}S_{k j}(1,1)_{-1 } \vac
+2\omega^{[k]}_{-1 } S_{k j}(1,2)_{-2 } \vac
\nonumber\\&\quad{}
+22\omega_{0 } \omega^{[k]}_{-1 } S_{k j}(1,2)_{-1 } \vac
-4\omega_{0 }^{3}S_{k j}(1,2)_{-1 } \vac
\nonumber\\&\quad{}
+6\omega_{0 } S_{k j}(1,3)_{-2 } \vac
,
\nonumber\\
S_{ij}(1,2)_{1}S_{k i}(1,3)&=
9\omega^{[k]}_{-2 } S_{k j}(1,1)_{-1 } \vac
-6\omega_{0 } \omega^{[k]}_{-1 } S_{k j}(1,1)_{-1 } \vac
\nonumber\\&\quad{}
+3\omega_{0 }^{3}S_{k j}(1,1)_{-1 } \vac
+6\omega^{[k]}_{-1 } S_{k j}(1,2)_{-1 } \vac
\nonumber\\&\quad{}
-9\omega_{0 } S_{k j}(1,2)_{-2 } \vac
+9\omega_{0 } S_{k j}(1,3)_{-1 } \vac
,
\nonumber\\
S_{ij}(1,2)_{2}S_{k i}(1,3)&=
6S_{k j}(1,3)_{-1 } \vac
,
\nonumber\\
S_{ij}(1,2)_{3}S_{k i}(1,3)&=
3S_{k j}(1,2)_{-1 } \vac
,
\end{align}

\begin{align}
S_{ij}(1,3)_{0}S_{ij}(1,3)&=
2\omega_{0 } \omega^{[i]}_{-1 } \Har^{[i]}_{-1 } \vac
+5\omega_{0 } \omega^{[i]}_{-2 } \omega^{[i]}_{-2 } \vac
+\omega_{0 } \omega^{[j]}_{-2 } \omega^{[j]}_{-2 } \vac
+2\omega_{0 } \omega^{[j]}_{-1 } \Har^{[j]}_{-1 } \vac
\nonumber\\&\quad{}
+\omega_{0 }^{3}\Har^{[i]}_{-1 } \vac
 -\omega_{0 }^{3}\omega^{[i]}_{-1 } \omega^{[i]}_{-1 } \vac
 -\omega_{0 }^{3}\omega^{[j]}_{-1 } \omega^{[j]}_{-1 } \vac
-3\omega_{0 }^{3}\Har^{[j]}_{-1 } \vac
\nonumber\\&\quad{}
+7\omega_{0 }^{5} \omega^{[i]}_{-1 } \vac
+\omega_{0 }^{5} \omega^{[j]}_{-1 } \vac
,
\nonumber\\
S_{ij}(1,3)_{1}S_{ij}(1,3)&=
4\omega^{[i]}_{-1 } \Har^{[i]}_{-1 } \vac
+5\omega^{[i]}_{-2 } \omega^{[i]}_{-2 } \vac
+\omega^{[j]}_{-2 } \omega^{[j]}_{-2 } \vac
+4\omega^{[j]}_{-1 } \Har^{[j]}_{-1 } \vac
\nonumber\\&\quad{}
+7\omega_{0 }^{2}\Har^{[i]}_{-1 } \vac
 -\omega_{0 }^{2}\omega^{[i]}_{-1 } \omega^{[i]}_{-1 } \vac
 -\omega_{0 }^{2}\omega^{[j]}_{-1 } \omega^{[j]}_{-1 } \vac
-3\omega_{0 }^{2}\Har^{[j]}_{-1 } \vac
\nonumber\\&\quad{}
+19\omega_{0 }^{4}\omega^{[i]}_{-1 } \vac
+\omega_{0 }^{4}\omega^{[j]}_{-1 } \vac
,
\nonumber\\
S_{ij}(1,3)_{2}S_{ij}(1,3)&=
15\omega_{0 } \Har^{[i]}_{-1 } \vac
+5\omega_{0 }^{3}\omega^{[i]}_{-1 } \vac
,
\nonumber\\
S_{ij}(1,3)_{3}S_{ij}(1,3)&=
30\Har^{[i]}_{-1 } \vac
+10\omega_{0 }^{2}\omega^{[i]}_{-1 } \vac
,
\nonumber\\
S_{ij}(1,3)_{4}S_{ij}(1,3)&=
30\omega_{0 } \omega^{[i]}_{-1 } \vac
,
\nonumber\\
S_{ij}(1,3)_{5}S_{ij}(1,3)&=
60\omega^{[i]}_{-1 } \vac
,
\nonumber\\
S_{ij}(1,3)_{6}S_{ij}(1,3)&=0,\nonumber\\
S_{ij}(1,3)_{7}S_{ij}(1,3)&=
30\vac
,
\end{align}

\begin{align}
S_{ij}(1,3)_{0}S_{k j}(1,1)&=
-3\omega^{[k]}_{-1 } S_{k i}(1,1)_{-2 } \vac
+\omega_{0 } \omega^{[k]}_{-1 } S_{k i}(1,1)_{-1 } \vac
\nonumber\\&\quad{}
+\omega_{0 }^{3}S_{k i}(1,1)_{-1 } \vac
+2\omega^{[k]}_{-1 } S_{k i}(1,2)_{-1 } \vac
\nonumber\\&\quad{}
-3S_{k i}(1,2)_{-3 } \vac
+3\omega_{0 } S_{k i}(1,3)_{-1 } \vac
,
\nonumber\\
S_{ij}(1,3)_{1}S_{k j}(1,1)&=
3S_{k i}(1,3)_{-1 } \vac
,
\nonumber\\
S_{ij}(1,3)_{2}S_{k j}(1,1)&=
3S_{k i}(1,2)_{-1 } \vac
,
\nonumber\\
S_{ij}(1,3)_{3}S_{k j}(1,1)&=
3S_{k i}(1,1)_{-1 } \vac
,
\end{align}

\begin{align}
S_{ij}(1,3)_{0}S_{k j}(1,2)&=
-16\omega^{[k]}_{-1 } S_{k i}(1,1)_{-3 } \vac
+4\Har^{[k]}_{-1 } S_{k i}(1,1)_{-1 } \vac
\nonumber\\&\quad{}
-44\omega_{0 } \omega^{[k]}_{-1 } S_{k i}(1,1)_{-2 } \vac
-2\omega_{0 } \omega^{[i]}_{-2 } S_{k i}(1,1)_{-1 } \vac
\nonumber\\&\quad{}
+16\omega_{0 }^{2}\omega^{[k]}_{-1 } S_{k i}(1,1)_{-1 } \vac
+\omega_{0 }^{4}S_{k i}(1,1)_{-1 } \vac
\nonumber\\&\quad{}
+10\omega^{[k]}_{-2 } S_{k i}(1,2)_{-1 } \vac
+4\omega^{[k]}_{-1 } S_{k i}(1,2)_{-2 } \vac
\nonumber\\&\quad{}
-2\omega_{0 } S_{k i}(1,2)_{-3 } \vac
+4\omega_{0 } \omega^{[i]}_{-1 } S_{k i}(1,2)_{-1 } \vac
,
\nonumber\\
S_{ij}(1,3)_{1}S_{k j}(1,2)&=
-12\omega^{[k]}_{-1 } S_{k i}(1,1)_{-2 } \vac
+4\omega_{0 } \omega^{[k]}_{-1 } S_{k i}(1,1)_{-1 } \vac
\nonumber\\&\quad{}
+2\omega_{0 }^{3}S_{k i}(1,1)_{-1 } \vac
+8\omega^{[k]}_{-1 } S_{k i}(1,2)_{-1 } \vac
\nonumber\\&\quad{}
-12S_{k i}(1,2)_{-3 } \vac
+12\omega_{0 } S_{k i}(1,3)_{-1 } \vac
,
\nonumber\\
S_{ij}(1,3)_{2}S_{k j}(1,2)&=
12S_{k i}(1,3)_{-1 } \vac
,
\nonumber\\
S_{ij}(1,3)_{3}S_{k j}(1,2)&=
12S_{k i}(1,2)_{-1 } \vac
,
\nonumber\\
S_{ij}(1,3)_{4}S_{k j}(1,2)&=
12S_{k i}(1,1)_{-1 } \vac
,
\end{align}

\begin{align}
S_{ij}(1,3)_{0}S_{k j}(1,3)&=
30\Har^{[i]}_{-1 } S_{k i}(1,2)_{-1 } \vac
+20\omega^{[i]}_{-3 } S_{k i}(1,2)_{-1 } \vac
-30\omega^{[i]}_{-2 } S_{k i}(1,3)_{-1 } \vac
,
\nonumber\\
S_{ij}(1,3)_{1}S_{k j}(1,3)&=
-40\omega^{[k]}_{-1 } S_{k i}(1,1)_{-3 } \vac
+10\Har^{[k]}_{-1 } S_{k i}(1,1)_{-1 } \vac
\nonumber\\&\quad{}
-110\omega_{0 } \omega^{[k]}_{-1 } S_{k i}(1,1)_{-2 } \vac
-5\omega_{0 } \omega^{[i]}_{-2 } S_{k i}(1,1)_{-1 } \vac
\nonumber\\&\quad{}
+40\omega_{0 }^{2}\omega^{[k]}_{-1 } S_{k i}(1,1)_{-1 } \vac
+5\omega_{0 }^{4}S_{k i}(1,1)_{-1 } \vac
\nonumber\\&\quad{}
+25\omega^{[k]}_{-2 } S_{k i}(1,2)_{-1 } \vac
+10\omega^{[k]}_{-1 } S_{k i}(1,2)_{-2 } \vac
\nonumber\\&\quad{}
-5\omega_{0 } S_{k i}(1,2)_{-3 } \vac
+10\omega_{0 } \omega^{[i]}_{-1 } S_{k i}(1,2)_{-1 } \vac
,
\nonumber\\
S_{ij}(1,3)_{2}S_{k j}(1,3)&=
-30\omega^{[k]}_{-1 } S_{k i}(1,1)_{-2 } \vac
+10\omega_{0 } \omega^{[k]}_{-1 } S_{k i}(1,1)_{-1 } \vac
\nonumber\\&\quad{}
+5\omega_{0 }^{3}S_{k i}(1,1)_{-1 } \vac
+20\omega^{[k]}_{-1 } S_{k i}(1,2)_{-1 } \vac
\nonumber\\&\quad{}
-30S_{k i}(1,2)_{-3 } \vac
+30\omega_{0 } S_{k i}(1,3)_{-1 } \vac
,
\nonumber\\
S_{ij}(1,3)_{3}S_{k j}(1,3)&=
30S_{k i}(1,3)_{-1 } \vac
,
\nonumber\\
S_{ij}(1,3)_{4}S_{k j}(1,3)&=
30S_{k i}(1,2)_{-1 } \vac
,
\nonumber\\
S_{ij}(1,3)_{5}S_{k j}(1,3)&=
30S_{k i}(1,1)_{-1 } \vac
,
\end{align}

\begin{align}
S_{ij}(1,3)_{0}S_{k i}(1,1)&=
3\omega^{[k]}_{-2 } S_{k j}(1,1)_{-1 } \vac
-2\omega_{0 } \omega^{[k]}_{-1 } S_{k j}(1,1)_{-1 } \vac
\nonumber\\&\quad{}
+\omega_{0 }^{3}S_{k j}(1,1)_{-1 } \vac
+2\omega^{[k]}_{-1 } S_{k j}(1,2)_{-1 } \vac
\nonumber\\&\quad{}
-3\omega_{0 } S_{k j}(1,2)_{-2 } \vac
+3\omega_{0 } S_{k j}(1,3)_{-1 } \vac
,
\nonumber\\
S_{ij}(1,3)_{1}S_{k i}(1,1)&=
S_{k j}(1,3)_{-1 } \vac
,
\end{align}

\begin{align}
S_{ij}(1,3)_{0}S_{k i}(1,2)&=
-16\omega^{[k]}_{-1 } S_{k j}(1,1)_{-3 } \vac
+4\Har^{[k]}_{-1 } S_{k j}(1,1)_{-1 } \vac
\nonumber\\&\quad{}
-98\omega_{0 } \omega^{[k]}_{-1 } S_{k j}(1,1)_{-2 } \vac
+34\omega_{0 }^{2}\omega^{[k]}_{-1 } S_{k j}(1,1)_{-1 } \vac
\nonumber\\&\quad{}
+3\omega_{0 }^{4}S_{k j}(1,1)_{-1 } \vac
+2\omega^{[k]}_{-1 } S_{k j}(1,2)_{-2 } \vac
\nonumber\\&\quad{}
+22\omega_{0 } \omega^{[k]}_{-1 } S_{k j}(1,2)_{-1 } \vac
-4\omega_{0 }^{3}S_{k j}(1,2)_{-1 } \vac
\nonumber\\&\quad{}
+6\omega_{0 } S_{k j}(1,3)_{-2 } \vac
,
\nonumber\\
S_{ij}(1,3)_{1}S_{k i}(1,2)&=
6\omega^{[k]}_{-2 } S_{k j}(1,1)_{-1 } \vac
-4\omega_{0 } \omega^{[k]}_{-1 } S_{k j}(1,1)_{-1 } \vac
\nonumber\\&\quad{}
+\omega_{0 }^{3}S_{k j}(1,1)_{-1 } \vac
+4\omega^{[k]}_{-1 } S_{k j}(1,2)_{-1 } \vac
\nonumber\\&\quad{}
-3\omega_{0 } S_{k j}(1,2)_{-2 } \vac
+6\omega_{0 } S_{k j}(1,3)_{-1 } \vac
,
\nonumber\\
S_{ij}(1,3)_{2}S_{k i}(1,2)&=
2S_{k j}(1,3)_{-1 } \vac
,
\end{align}

\begin{align}
S_{ij}(1,3)_{0}S_{k i}(1,3)&=
-120\omega^{[k]}_{-1 } \omega^{[j]}_{-1 } S_{k j}(1,1)_{-2 } \vac
-300\omega^{[j]}_{-4 } S_{k j}(1,1)_{-1 } \vac
\nonumber\\&\quad{}
-180\Har^{[j]}_{-1 } S_{k j}(1,1)_{-2 } \vac
-60\Har^{[j]}_{-2 } S_{k j}(1,1)_{-1 } \vac
\nonumber\\&\quad{}
+40\omega_{0 } \omega^{[k]}_{-1 } \omega^{[j]}_{-1 } S_{k j}(1,1)_{-1 } \vac
+240\omega_{0 } \omega^{[j]}_{-3 } S_{k j}(1,1)_{-1 } \vac
\nonumber\\&\quad{}
-90\omega_{0 }^{2}\omega^{[j]}_{-2 } S_{k j}(1,1)_{-1 } \vac
+20\omega_{0 }^{3}\omega^{[j]}_{-1 } S_{k j}(1,1)_{-1 } \vac
\nonumber\\&\quad{}
+330\omega^{[k]}_{-2 } S_{k j}(1,2)_{-2 } \vac
+780\omega^{[k]}_{-1 } S_{k j}(1,2)_{-3 } \vac
\nonumber\\&\quad{}
-90\Har^{[k]}_{-1 } S_{k j}(1,2)_{-1 } \vac
+180\Har^{[j]}_{-1 } S_{k j}(1,2)_{-1 } \vac
\nonumber\\&\quad{}
-120\omega_{0 }^{2}\omega^{[k]}_{-1 } S_{k j}(1,2)_{-1 } \vac
-135\omega_{0 }^{4}S_{k j}(1,2)_{-1 } \vac
\nonumber\\&\quad{}
-470\omega^{[k]}_{-2 } S_{k j}(1,3)_{-1 } \vac
-500\omega^{[k]}_{-1 } S_{k j}(1,3)_{-2 } \vac
\nonumber\\&\quad{}
+750S_{k j}(1,3)_{-4 } \vac
,
\nonumber\\
S_{ij}(1,3)_{1}S_{k i}(1,3)&=
-8\omega^{[k]}_{-1 } S_{k j}(1,1)_{-3 } \vac
+2\Har^{[k]}_{-1 } S_{k j}(1,1)_{-1 } \vac
\nonumber\\&\quad{}
-49\omega_{0 } \omega^{[k]}_{-1 } S_{k j}(1,1)_{-2 } \vac
+17\omega_{0 }^{2}\omega^{[k]}_{-1 } S_{k j}(1,1)_{-1 } \vac
\nonumber\\&\quad{}
+9\omega_{0 }^{4}S_{k j}(1,1)_{-1 } \vac
+\omega^{[k]}_{-1 } S_{k j}(1,2)_{-2 } \vac
\nonumber\\&\quad{}
+11\omega_{0 } \omega^{[k]}_{-1 } S_{k j}(1,2)_{-1 } \vac
-6\omega_{0 }^{3}S_{k j}(1,2)_{-1 } \vac
\nonumber\\&\quad{}
+9\omega_{0 } S_{k j}(1,3)_{-2 } \vac
,
\nonumber\\
S_{ij}(1,3)_{2}S_{k i}(1,3)&=
9\omega^{[k]}_{-2 } S_{k j}(1,1)_{-1 } \vac
-6\omega_{0 } \omega^{[k]}_{-1 } S_{k j}(1,1)_{-1 } \vac
\nonumber\\&\quad{}
+3\omega_{0 }^{3}S_{k j}(1,1)_{-1 } \vac
+6\omega^{[k]}_{-1 } S_{k j}(1,2)_{-1 } \vac
\nonumber\\&\quad{}
-9\omega_{0 } S_{k j}(1,2)_{-2 } \vac
+9\omega_{0 } S_{k j}(1,3)_{-1 } \vac
,
\nonumber\\
S_{ij}(1,3)_{3}S_{k i}(1,3)&=
3S_{k j}(1,3)_{-1 } \vac
,\label{eq:last-appendix}
\end{align}

\renewcommand{\thesection}{Notation}

\section{}\label{section:notation}

\begin{tabular}{lp{13cm}}
$V$ & a vertex algebra.\\
$U$ & a subspace of a weak $V$-module.\\
$\Omega_{V}(U)$&$
=\{\lu\in U\ \Big|\ 
a_{i}\lu=0\ \mbox{for all homogeneous }a\in V\mbox{and }i>\wt a-1\}$.\\
$\hei$ & a finite dimensional vector space equipped with a nondegenerate symmetric bilinear form
$\langle \mbox{ }, \mbox{ }\rangle$.\\
$h^{[1]},\ldots,h^{[\rankL]}$ & an orthonormal basis of $\fh$.\\
$M(1)$ & the vertex operator algebra associated to the Heisenberg algebra.\\
$\lattice$ & a non-degenerate even lattice of finite rank.\\
$\rankL$ & the rank of $\lattice$.\\
$V_{\lattice}$ & the vertex algebra associated to $\lattice$.\\
$\theta$ & the automorphism of $V_{\lattice}$ induced from the $-1$ symmetry of $\lattice$.\\
$M(1)^{+}$ & the fixed point subalgbra of $M(1)$ under the action of $\theta$.\\
$V_{\lattice}^{+}$ & the fixed point subalgbra of $V_{\lattice}$ under the action of $\theta$.\\
$I(\mbox{ },x)$ & an intertwining operator for $M(1)^{+}$.\\
$\epsilon(\lu,\lv)$ & 
$\lu_{\epsilon(\lu,\lv)}\lv\neq 0\mbox{ and }\lu_{i}\lv=0\mbox{ for all }i>\epsilon(\lu,\lv)$
if $I(\lu,x)\lv\neq 0$ and $\epsilon(\lu,\lv)=-\infty$ if $I(\lu,x)\lv= 0$,
where $I : \module\times\mW\rightarrow \mN\db{x}$ is an intertwining operator and 
$\lu\in\module$, $\lv\in\mW$ (see \eqref{eqn:max-vanish}).\\
$A(V)$ & the Zhu algebra of a vertex operator algebra $V$.\\
$A_{-}B$&$:=\Span_{\C}\{a_{-i}b\ |\ a\in A, b\in B,\mbox{ and }i\in\Z_{>0}\}$ (see \eqref{eqn:A-B:=Span}).\\
$\langle A_{-}\rangle B$ & see \eqref{eq:a(1)-i1cdotsa(n)-inb}.\\
$\omega$&$=(1/2)\sum_{i=1}^{\rankL}h^{[i]}(-1)^2\vac$.\nonumber\\
$\ExB(\alpha)$&$=e^{\alpha}+\theta(e^{\alpha})$ where $\alpha\in\fh$.\\
$\omega^{[i]}$&$=(1/2)h^{[i]}(-1)^2$.\nonumber\\
$\Har^{[i]}$&$=(1/3)(h^{[i]}(-3)h^{[i]}(-1)\vac-h^{[i]}(-2)^2\vac)$.
\end{tabular}


\bigskip
\noindent\textbf{Conflict of Interest}\\
The author declares that he has no conflict of interest.

\bigskip
\noindent\textbf{Acknowledgments}\\
The author thanks the anonymous referee for helpful comments. 

\providecommand{\MR}{\relax\ifhmode\unskip\space\fi MR }
\providecommand{\MRhref}[2]{%
  \href{http://www.ams.org/mathscinet-getitem?mr=#1}{#2}
}
\providecommand{\href}[2]{#2}

\bibliographystyle{ijmart}

\end{document}